\documentclass[12pt]{amsart}

\usepackage{amsmath,amsfonts,amssymb,amscd,verbatim,delarray}
\newcommand{\A}{{\mathbb A}}
\newcommand{\F}{{\mathbb F}}
\renewcommand{\P}{{\mathbb P}}
\newcommand{\Q}{{\mathbb Q}}
\newcommand{\Z}{{\mathbb Z}}

\newcommand{\kf}{{\mathcal F}}
\newcommand{\kc}{{\mathcal C}}
\newcommand{\lra}{\longrightarrow}

 \DeclareMathOperator{\spoly}{spoly}
 \DeclareMathOperator{\Fix}{Fix}
\DeclareMathOperator{\Fr}{Fr}

\DeclareMathOperator{\const}{const}

\DeclareMathOperator{\Spec}{Spec}

\begin{document}
\numberwithin{equation}{section}

\newtheorem{theorem}{Theorem}[section]
\newtheorem{lemma}[theorem]{Lemma}
\newtheorem{prop}[theorem]{Proposition}
\newtheorem{proposition}[theorem]{Proposition}
\newtheorem{corollary}[theorem]{Corollary}
\newtheorem{corol}[theorem]{Corollary}
\newtheorem{conj}[theorem]{Conjecture}

\theoremstyle{definition}
\newtheorem{defn}[theorem]{Definition}
\newtheorem{example}[theorem]{Example}
\newtheorem{examples}[theorem]{Examples}
\newtheorem{remarks}[theorem]{Remarks}
\newtheorem{remark}[theorem]{Remark}
\newtheorem{algorithm}[theorem]{Algorithm}
\newtheorem{question}[theorem]{Question}
\newtheorem{problem}[theorem]{Problem}
\newtheorem{subsec}[theorem]{}


\def\toeq{{\stackrel{\sim}{\longrightarrow}}}
\def\into{{\hookrightarrow}}


\def\alp{{\alpha}}  \def\bet{{\beta}} \def\gam{{\gamma}}
 \def\del{{\delta}}
\def\eps{{\varepsilon}}
\def\kap{{\kappa}}                   \def\Chi{\text{X}}
\def\lam{{\lambda}}
 \def\sig{{\sigma}}  \def\vphi{{\varphi}} \def\om{{\omega}}
\def\Gam{{\Gamma}}   \def\Del{{\Delta}}
\def\Sig{{\Sigma}}   \def\Om{{\Omega}}
\def\ups{{\upsilon}}


\def\F{{\mathbb{F}}}
\def\Q{{\mathbb{Q}}}
\def\Ql{{\overline{\Q }_{\ell }}}
\def\CC{{\mathbb{C}}}
\def\R{{\mathbb R}}
\def\V{{\mathbf V}}
\def\D{{\mathbf D}}

\def\XX{\mathbf{X}^*}
\def\xx{\mathbf{X}_*}

\def\AA{\Bbb A}
\def\HH{\mathbb H}
\def\PP{\Bbb P}

\def\Gm{{{\mathbb G}_{\textrm{m}}}}
\def\Gmk{{{\mathbb G}_{\textrm m,k}}}
\def\GmL{{\mathbb G_{{\textrm m},L}}}
\def\Ga{{{\mathbb G}_a}}

\def\Fb{{\overline{\F }}}
\def\Kb{{\overline K}}
\def\Yb{{\overline Y}}
\def\Xb{{\overline X}}
\def\Tb{{\overline T}}
\def\Bb{{\overline B}}
\def\Gb{{\bar{G}}}
\def\Ub{{\overline U}}
\def\Vb{{\overline V}}
\def\Hb{{\bar{H}}}
\def\kb{{\bar{k}}}

\def\Th{{\hat T}}
\def\Bh{{\hat B}}
\def\Gh{{\hat G}}

\def\cF{{\mathfrak{F}}}
\def\cC{{\mathcal C}}
\def\cU{{\mathcal U}}

\def\Xt{{\widetilde X}}
\def\Gt{{\widetilde G}}

\def\gg{{\mathfrak g}}
\def\hh{{\mathfrak h}}

\def\min{^{-1}}

\def\textrm#1{\text{\textnormal{#1}}}

\def\GL{\textrm{GL}}            \def\Stab{\textrm{Stab}}
\def\Gal{\textrm{Gal}}          \def\Aut{\textrm{Aut\,}}
\def\Lie{\textrm{Lie\,}}        \def\Ext{\textrm{Ext}}
\def\PSL{\textrm{PSL}}          \def\SL{\textrm{SL}}
\def\loc{\textrm{loc}}
\def\coker{\textrm{coker\,}}    \def\Hom{\textrm{Hom}}
\def\im{\textrm{im\,}}           \def\int{\textrm{int}}
\def\inv{\textrm{inv}}           \def\can{\textrm{can}}
\def\id{\textrm{id}}
\def\Cl{\textrm{Cl}}
\def\Sz{\textrm{Sz}}

\def\tors{_{\textrm{tors}}}      \def\tor{^{\textrm{tor}}}
\def\red{^{\textrm{red}}}         \def\nt{^{\textrm{ssu}}}

\def\sss{^{\textrm{ss}}}          \def\uu{^{\textrm{u}}}
\def\ad{^{\textrm{ad}}}           \def\mm{^{\textrm{m}}}
\def\tm{^\times}                  \def\mult{^{\textrm{mult}}}

\def\uss{^{\textrm{ssu}}}         \def\ssu{^{\textrm{ssu}}}
\def\comp{_{\textrm{c}}}
\def\ab{_{\textrm{ab}}}

\def\et{_{\textrm{\'et}}}
\def\nr{_{\textrm{nr}}}

\def\nil{_{\textrm{nil}}}
\def\sol{_{\textrm{sol}}}

\def\til{\;\widetilde{}\;}


\title{{\bf Engel-like Identities \\ Characterizing Finite Solvable Groups
}}
\author[Bandman, Greuel, Grunewald, Kunyavski\u\i , Pfister, Plotkin]{
Tatiana Bandman, Gert-Martin Greuel, Fritz Grunewald,\\
Boris Kunyavski\u\i , Gerhard Pfister, and Eugene Plotkin }

\address{Bandman, Kunyavski\u\i \ and Plotkin: Department of
Mathematics and Statistics, Bar-Ilan University, 52900 Ramat Gan,
ISRAEL}

\email{bandman@macs.biu.ac.il , kunyav@macs.biu.ac.il,
plotkin@macs.biu.ac.il}

\address{Greuel and Pfister: Fachbereich Mathematik, Universit\"at
Kaiserslautern, Postfach 3049, 67653 Kaisers\-lautern, GERMANY}

\email{greuel@mathematik.uni-kl.de, pfister@mathematik.uni-kl.de}

\address{Grunewald: Mathematisches Institut der Universit\"at
Heinrich Heine D\"usseldorf, Universit\"atsstr. 1, 40225
D\"usseldorf, GERMANY}

\email{grunewald@math.uni-duesseldorf.de}

\date{\today}
\maketitle

\thispagestyle{empty} \vspace{1.0cm} \setcounter{page}{0}

\tableofcontents
\newpage

\section{Introduction} \label{sec:intro}

\subsection{Statement of the problem and main results}

The starting point for this research is the following classical
fact: the class of finite nilpotent groups is characterized by
Engel identities. To be more precise, Zorn's theorem \cite{Zo}
says that a finite group $G$ is nilpotent if and only if it
satisfies one of the identities $e_n(x,y)=[y,x,x,\dots ,x]=1$
(here $[x,y]=xyx\min y\min$, $[y,x,x]=[[y,x],x]$, etc.).


Our goal is to characterize the class of finite solvable groups by
two-variable identities in a similar way.
More precisely, a sequence of words $u_1,\ldots,u_n,\ldots $ is
called correct if $u_k\equiv 1$ in a group $G$ implies $u_m\equiv
1$ in a group $G$ for all $m>k$. We are looking for an explicit
correct sequence of words $u_1(x,y),\ldots,u_n(x,y),\ldots$ such that a
group $G$ is solvable if and only if for some $n$ the word $u_n$
is an identity in $G$.

B.~Plotkin suggested some Engel-like identities which could
characterize finite solvable groups (see \cite{PPT}, \cite{GKNP}).
In the present paper we establish B.~Plotkin's conjecture (in a
slightly modified form).

Let $w$ denote a word in $x$, $y$, $x\min$, $y\min$, and let
${}^w\!u_n(x,y)$ be an infinite sequence defined by the rule
\begin{equation}
\begin{array}{ccl}
{}^w\!u_1 & = & w, \\
{}^w\!u_{n+1} & = & [x\,{}^w\!u_n\,x\min,y\,{}^w\!u_n\,y\min],
\dots
\end{array}
\label{seq:gen}
\end{equation}

Here is our main result.

\begin{theorem}   \label{conj:P}
There exists $w$ such that a finite group $G$ is solvable if and
only if for some $n$ the identity ${}^w\!u_n(x,y)\equiv 1$ holds
in $G$.
\end{theorem}

In fact, we exhibit an explicit initial term $w=x^{-2}y\min x$ for
which the statement of the theorem holds. Note two obvious
properties of $w$: 1) if a group $G$ satisfies the identity
$w\equiv 1$, then $G=\{1\}$; 2) the words $w$ and $x$ generate the
free group $F=\left<x,y\right>$. Thus $w$ can be used also as the
initial term of a sequence characterizing finite nilpotent groups,
see Proposition \ref{prop:Engel} below.

\begin{remark}   \label{rem:gen}
We believe that in the statement of Theorem \ref{conj:P} the
initial word $w$ can be chosen in the most natural way: $w=[x,y]$,
as in the nilpotent case.
\end{remark}


A theorem of J.~Thompson (\cite{Th}, \cite{Fl}) states that if $G$
is a finite group in which every two elements generate a solvable
subgroup then $G$ is solvable. As mentioned in \cite{BW}, together
with \cite[Satz 2.12]{Br} this implies that finite solvable groups
can be characterized by a countable set of two-variable
identities. (Note that this fact also follows from Lemma 16.1 and
Theorem 16.21 from \cite{Ne} saying that an $n$-generator group
$G$ belongs to a variety $V$ if and only if all $n$-variable
identities from $V$ are fulfilled in $G$.) However, this does not
provide  explicit two-variable identities for finite solvable
groups. Furthermore, in the above cited paper R.~Brandl and
J.~S.~Wilson construct a countable set of words $w_n(x,y)$ with
the property that a finite group $G$ is solvable if and only if
for almost all $n$ the identity $w_n(x,y)\equiv 1$ holds in $G$.
Since in their construction there is no easily described
relationship between terms of $w_n(x,y)$, they raise the question
whether one can characterize finite solvable groups by {\bf
sequences} of identities fitting into a simple recursive
definition.

Recently A.~Lubotzky proved that for any integer $d\geq 2$ the
free prosolvable group $\hat F_d(S)$ can be defined by a {\it
single} profinite relation \cite[Prop. 3.4]{Lu}. Using this
proposition and Thompson's theorem, one can derive the existence
of a needed sequence of identities characterizing finite solvable
groups (Lubotzky's result does not give, however, any candidate
for such a sequence).

The sequence constructed in our Theorem \ref{conj:P} answers the
question of Brandl--Wilson and fits very well into profinite
setting (see Subsection \ref{app:pro}).

One can mention here some more cases where certain interesting
classes of finite groups were characterized by two-variable
commutator identities \cite{Br}, \cite{BP}, \cite{BN}, \cite{Gu},
\cite{GH}, \cite{Ni1}, \cite{Ni2}; see \cite{GKNP} or the above
cited papers for more details.

\bigskip


Although Theorem \ref{conj:P} is a purely group-theoretic result,
its proof is surprisingly diverse involving a good bunch of
algebraic geometry and computer algebra (note, however, a paper of
Bombieri \cite{Bo} which served for us as an inspiring example of
such an approach). We want to emphasize a special role played by
problem oriented software (particularly, the packages \textsc{Singular} and
MAGMA): not only proofs but even precise statements of our results
would hardly be possible without extensive computer experiments.

We shall outline the general strategy of the proof in Section
\ref{sec:gen} and give the details in Sections \ref{sec:PSL} and
\ref{sec:Suzuki}. Before that we present some digressions related
to the main theorem (motivations, possible applications, open
problems, etc.)

\subsection{Analogues, problems, and generalizations}

Theorem \ref{conj:P} admits some natural analogues in
Lie-algebraic and group-schematic settings \cite{GKNP}. In particular, the
following analogue of the classical Engel theorem on nilpotent Lie
algebras is true.

\begin{theorem} \textrm{\cite{GKNP}}
Let $L$ be a finite dimensional Lie algebra defined over an
infinite field $k$ of characteristic different from $2$, $3$, $5$.
Define
\begin{equation}
v_1=[x,y], \quad v_{n+1}=[[v_n,x],[v_n,y]] \quad (n>1).
\label{eq:solv}
\end{equation}
Then $L$ is solvable if and only if for some $n$ one of the
identities $v_n(x,y)\equiv 0$ holds in $L$.
\end{theorem}

(Here [ , ] are Lie brackets.)

A much more challenging question is related to the
infinite-dimensional case. Namely, the remarkable
Kostrikin--Zelmanov theorem on locally nilpotent Lie algebras
\cite{Ko}, \cite{Ze2}, \cite{Ze3} and Zelmanov's theorem
\cite{Ze1} lead to the following

\begin{problem} \label{quest:Z}
Suppose that $L$ is a Lie algebra over a field $k$, the $v_n$'s
are defined by formulas $(\ref{eq:solv})$, and there is $n$ such
that the identity $v_n(x,y)\equiv 0$ holds in $L$. Is it true that
$L$ is locally solvable? If $k$ is of characteristic $0$, is it
true that $L$ is solvable?
\end{problem}

(A property is said to hold locally if it holds for all finitely
generated Lie subalgebras.)

\medskip

Of course, it would be of significant interest to consider similar
questions for arbitrary groups.

We call $G$ an {\it Engel} group if there is an integer $n$ such
that the Engel identity $e_n(x,y)\equiv 1$ holds in $G$.

We call $G$ an {\it unbounded Engel} group if for every $x,y\in G$
there is an integer $n=n(x,y)$ such that $e_n(x,y)=1$.

We introduce the following

\begin{defn}
We call $G$ a {\it $w$-quasi-Engel} group if there is an integer
$n$ such that the identity ${}^w\!u_n(x,y)\equiv 1$ holds in $G$.
\end{defn}

\begin{defn}
We call $G$ an {\it unbounded $w$-quasi-Engel} group if for every
$x,y\in G$ there is an integer $n=n(x,y)$ such that
${}^w\!u_n(x,y)=1$.
\end{defn}

Here and throughout below we assume that the initial word $w$ is
chosen so that Theorem \ref{conj:P} holds, and fix such a $w$.
Thus we shall drop $w$ in the above definitions and use the terms
`quasi-Engel group', `unbounded quasi-Engel group'.

\begin{problem}\label{prob:Engel}
Is every Engel group locally nilpotent?
\end{problem}

\begin{problem}\label{prob:q-Engel}
Is every quasi-Engel group locally solvable?
\end{problem}

(A property is said to hold locally if it holds for all finitely
generated subgroups.)

Problem \ref{prob:Engel} remains open for a long time, cf.
\cite{Plo3}. The answer in general is most likely negative,
however some positive results are known \cite{BM}, \cite{Gr},
\cite{Plo1}, \cite{Plo2}, \cite{Wi}, \cite{WZ}, etc. In the
solvable case the situation is even less clear. We dare to state
the following

\begin{conj}
Every residually finite, quasi-Engel group is locally solvable.
\end{conj}

(A group is said to be residually finite if the intersection of
all its normal subgroups of finite index is trivial.)

For profinite groups the situation looks more promising.

\begin{theorem} \textrm{\cite[Th. 5]{WZ}}
Every profinite, unbounded Engel group is locally nilpotent.
\end{theorem}

\begin{conj} \label{conj:pro}
Every profinite, unbounded quasi-Engel group is locally solvable.
\end{conj}

It is quite natural to consider restricted versions of Problems
\ref{prob:Engel} and \ref{prob:q-Engel} like it is considered for
the Burnside problem. Let $E_n$ be the Engel variety defined by
the identity $e_n\equiv 1$. Let $F=F_{k,n}$ be the free group with
$k$ generators in the variety $E_n$. One can prove that the
intersection of all conilpotent normal subgroups $H_\alpha$ in $F$
is also conilpotent. Hence there exists a group $F^0_{n,k}$ in
$E_n$ such that every nilpotent group $G\in E_n$ with $k$
generators is a homomorphic image of  $F^0_{n,k}$. This implies
that all locally nilpotent groups from $E_n$ form a variety. In
other words, the restricted Engel problem has a positive solution.
The situation with the restricted quasi-Engel problem is unclear.

\begin{problem}
Let $F=F_{k,n}$ be the free group with $k$ generators in the
variety of all quasi-Engel groups with fixed $n$. Is it true that
the intersection of all cosolvable normal subgroups in $F=F_{k,n}$
is also cosolvable?
\end{problem}
Our main theorem can be reformulated in profinite terms.

\begin{theorem} \label{qe0}
Let $F=F(x,y)$ denote the free group in two variables, and let
$\widehat F$ be its profinite completion. Let $v_1, v_2, \ldots,
v_m , \ldots$ be any convergent subsequence of $(\ref{seq:gen})$
with limit $f$ from $\widehat F$. Then the identity $f\equiv 1$
defines the profinite variety of prosolvable groups.
\end{theorem}

(See Section \ref{app:pro} for more details.)

It would be of great interest to consider the restricted
quasi-Engel problem for profinite groups.

\begin{remark}
There is no sense in generalizing Conjecture \ref{conj:pro} too
far: from the Golod--Shafarevich counterexamples one can deduce an
example of an unbounded quasi-Engel group which is not locally
nilpotent (and hence not locally solvable). We thank B.~Plotkin
for this observation.
\end{remark}

Note that our results can be viewed as a natural development of
the classical Thompson--Flavell theorem. Indeed, the main Theorem
\ref{conj:P} immediately implies that a finite group $G$ is
solvable if and only if every two-generated subgroup is solvable
\cite{Th}, \cite{Fl}.

See Corollary \ref{cor:two-gen} for the profinite setting.

Finally, consider an interesting particular case of linear groups.

\begin{corol} \label{cor:lin}
Suppose that $G\subset\GL (n,K)$ where $K$ is a field. Then $G$ is
solvable if and only if it is quasi-Engel.
\end{corol}

\begin{proof}
The ``only if'' part is obvious. The ``if'' part is an immediate
consequence of Theorem \ref{conj:P} and Platonov's theorem
\cite{Pla} stating that every linear group over a field satisfying
a non-trivial identity has a solvable subgroup of finite index.
(Of course, if $K$ is of characteristic zero, the assertion
follows from the Tits alternative \cite{Ti}.)
\end{proof}

\subsection{Related graphs}

Let us return to finite groups. Recall that to every finite group
$G$ one can attach its commuting graph $\Gam (G)$ with vertices
the non-identity elements of $G$ and edges joining commuting
elements. This graph plays a crucial role in recent works on the
group structure of the group of rational points $H(k)$ of a simple
linear algebraic group $H$ defined over a number field $k$ (to be
more precise, in establishing important cases of the
Margulis--Platonov conjecture on projective simplicity of $H(k)$
\cite{Se}, \cite{SS}, and of the Prasad--Rapinchuk conjecture on
solvability of the finite quotients of $H(k)$ \cite{RSS}).

Zorn's theorem and our Theorem \ref{conj:P} allow one to define
the nilpotency graph $\Gam\nil (G)$ and the solvability graph
$\Gam\sol (G)$ in a similar way: the vertices are the non-identity
elements of $G$ and two vertices $x,y$ are joined by an edge if
for some $n$ we have $e_n(x,y)=1$ (resp. $u_n(x,y)=1$); see
\cite{GKNP}. We believe that these graphs may provide an
additional helpful tool in studying properties of arithmetic
groups.

\bigskip

\noindent {\it Acknowledgements}. Bandman, Kunyavski\u\i , and
Plotkin were partially supported by the Ministry of Absorption
(Israel), the Israeli Science Foundation founded by the Israeli
Academy of Sciences --- Center of Excellence Program, and the
Minerva Foundation through the Emmy Noether Research Institute of
Mathematics. Kunyavski\u\i \ and Plotkin were also supported by
the RTN network HPRN-CT-2002-00287 and INTAS 00-566.

We are grateful to N.~Gordeev, D.~Grayson, L.~Illusie,
A.~Lubotzky, A.~Mann, S.~Margolis, R.~Pink, L.~Rowen, Y.~Segev,
J.-P.~Serre, Y.~Varshavsky, and N.~Vavilov for useful comments and
advise. We thank D.~Nikolova and R.~Shklyar for help in computer
experiments. Our special thanks go to B.~Plotkin for numerous
enlightening, encouraging, and inspiring discussions.

\bigskip

\noindent {\it Convention}. Throughout below we assume that
$w=x^{-2}y\min x$ is chosen as the initial term of the sequence
${}^w\!u_n$ and shorten ${}^w\!u_i$ to $u_i$.

\medskip
\noindent {\it Notation}. Because of extensive use of the
\textsc{Singular} package our notation sometimes differs from the
standard one: say, in the output of computer sessions, powers like
$a^{12}$ are denoted as {\tt a12} . We refer the reader to
\cite{GP3}, \cite{GP4}, \cite{GPS} for definitions of
\textsc{Singular} commands and their usage, and to \cite{Bu},
\cite{GP1}--\cite{GP3} for details on Gr\"obner bases.

All other notation is more or less standard.

{\it Rings and fields}: All rings are assumed commutative with 1;
$\mathbb Z$, $\mathbb Q$, $\mathbb F_q$ denote the ring of
integers, the field of rational numbers, the field of $q$
elements, respectively. $\overline k$ denotes a (fixed) algebraic
closure of a field $k$.

{\it Ideals and varieties}: If $I$ is an ideal in $R$ and $i\colon
R\to S$ is a ring homomorphism, $IS$ stands for the image of $I$
under $i$. The ideal generated by $f_1,\dots f_k$ is denoted
$\langle f_1,\dots , f_k\rangle$.

For $f\in R$ we denote $I:f^{\infty}=\cup_{n=1}^{\infty }I:f^n$.
If $R$ is noetherian, the chain of ideals $I:f\subseteq
I:f^2\subseteq\dots$ stabilizes, and we have $I:f^{\infty }=I:f^n$
for some $n$.

$\mathbb A^n$ and $\mathbb P^n$ denote affine and projective
spaces. $\overline{C}$ denotes the projective closure of an affine
set $C\subset\mathbb A^n$, and $I_h$ stands for the homogenization
of an ideal $I$. $\V (J)$ denotes the affine variety defined as
the set of common zeros of the functions from an ideal $J$.  If
$\V (J)\subset\mathbb A^n$, we denote $\D (J)=\mathbb
A^n\smallsetminus\V(J)$. We shorten $\V (\langle f_1,\dots
,f_k\rangle)$ to $\V (f_1,\dots ,f_k)$, and $\D (\langle f_1,\dots
,f_k\rangle)$ to $\D (f_1,\dots ,f_k)$. We denote by $V(k)$ the
set of rational points of a $k$-variety $V$.

$\chi(V)$ denotes the Euler characteristic of (the tangent bundle
of) a variety $V.$

If $C$ is a projective curve (maybe singular), $p_a(C)$ is the
arithmetic genus of $C$, and $g(C)$ denotes the arithmetic genus
of the normalization of $C.$


{\it Groups}: $\PSL (n,q)$ denotes the projective special linear
group of degree $n$ over $\mathbb F_q$. For $q=2^m$ we denote by
$\Sz (q)$ the Suzuki group (the twisted form of ${}^2\!B_2$, see
\cite{HB}).

All other notation will be explained when needed.

\section{General strategy} \label{sec:gen}

To prove Theorem \ref{conj:P}, we proceed as follows (cf.
\cite{GKNP}). We derive Theorem \ref{conj:P} from the following

\begin{theorem} \label{conj:main}
Let $G$ be one of the following groups:

\begin{enumerate}
\item $\PSL (2,p)$ $(p=5 \textrm{ or } p=\pm 2 \pmod 5, p\neq 3)$,
\item $\PSL (2,2^p)$,
\item $\PSL (2,3^p)$ $(p \textrm{ odd })$,
\item $\PSL (3,3)$,
\item $\Sz (2^p)$ $(p \textrm{ odd })$.
\end{enumerate}

Then there exists a word $w$ in $x$, $y$, $x\min$, $y\min$,
independent of $G$, such that none of the identities
$u_n(x,y)\equiv 1$ holds in $G$.
\end{theorem}

\begin{prop}
Theorem $\ref{conj:main}$ implies Theorem $\ref{conj:P}$.
\end{prop}

\begin{proof} First note that the ``only if'' part of the statement
of Theorem \ref{conj:P} is obvious. Indeed, if $G$ is solvable of
class $n$ then the identity $u_n\equiv 1$ holds in $G$ for any $w$
since the value $u_n(x,y)$ belongs to the corresponding term of
the derived series. Thus we only have to prove that the ``if" part
of the statement of Theorem \ref{conj:P} follows from Theorem
\ref{conj:main}. So let us assume that Theorem \ref{conj:main}
holds, take $w$ as in its statement, and suppose that there exists
a non-solvable finite group in which the identity $u_n\equiv 1$
holds. Denote by $G$ a minimal counter-example, i.e. a finite
non-solvable group of the smallest order with identity $u_n\equiv
1$. Such $G$ must be simple. Indeed, if $H$ is a proper normal
subgroup of $G$, then both $H$ and $G/H$ are solvable (because any
identity remains true in the subgroups and the quotients). But the
list of groups in Theorem \ref{conj:main} is none other than the
list of finite simple groups all of whose subgroups are solvable
\cite{Th}. Thus for any $G$ from this list the identity $u_n\equiv
1$ does not hold in $G$, contradiction.
\end{proof}

To prove Theorem \ref{conj:main}, it is enough to find a word $w$
and integers $i$ and $j$ such that the equation
\begin{equation}
u_i(x,y)=u_j(x,y) \label{eq:ui=uj}
\end{equation}
has a non-trivial solution in every $G$ from the above list
(non-trivial means that $u_i(x,y)\neq 1$). In the next sections we
explain how this can be done.

\subsection{First screening} \label{subsec:screen}

A possible attempt to find numerical evidence in support of the
main theorem could be as follows: pick up a word $w$ (say, take
$w=[x,y]$, as in the classical Engel sequence), and consider
equation (\ref{eq:ui=uj}) for small $i,j$ in each group $G$ from
the list of Theorem \ref{conj:main}. Let us focus on the case
$G=\PSL (2,p)$, and consider $1\le i,j\le 4$. Computer experiments
(with the help of MAPLE) immediately show arising difficulties.
Although the number of solutions to the above equations has a
tendency to grow with growth of $p$, for each pair $i,j$ there is
$p$ such that equation (\ref{eq:ui=uj}) has no non-trivial
solutions in $\PSL (2,p)$.

Here is a way out (cf.\ Proposition \ref{prop:Engel}): we vary the initial
word of the sequence. For simplicity, we limit ourselves to the equation
\begin{equation}
u_1(x,y)=u_2(x,y). \label{eq:u1=u2}
\end{equation}
The result may seem unexpected enough: there are certain words
(less than 0.1\% of the total number of words of given length)
such that equation (\ref{eq:u1=u2}) has a non-trivial solution for
all $p<1000$; moreover, for such initial words the rate of growth
of the number of solutions is significantly higher than for others
(the shortest words of this type are 1) $w=x\min yxy\min x$; 2)
$w=x^{-2}y\min x$; 3) $w=y^{-2}x\min y$).

This purely experimental numerical phenomenon allows us to reveal
even deeper properties of the equations under consideration. These
properties are of algebraic-geometric nature and are of key
importance for further investigation.


\subsection{Algebraic-geometric view} \label{subsec:ag}

The general idea can be described as follows. For a group $G$ in
the list of Theorem \ref{conj:main}, we fix its standard linear
representation (over the corresponding finite field $\F _q$). Then
the equation $u_1(x,y)=u_2(x,y)$ can be viewed as a matrix
equation. To be more precise, we regard the entries of the
matrices corresponding to $x$ and $y$ in this representation as
variables, and thus the above matrix equation becomes a system of
polynomial equations defining an algebraic variety over $\F _q$.
Our goal is to apply to this variety estimates of Lang--Weil type
which guarantee the existence of a solution for $q$ big enough
(see \cite{LW}). Small values of $q$ are checked case by case.

Here is how this strategy looks like in our setting.

\subsection{$\PSL(2)$ case} As mentioned at the end of
Subsection \ref{subsec:screen}, our experimental data can only be
explained by some algebraic-geometric phenomena. So, following the
general strategy described above, we fix an initial word $w$ and
represent equation (\ref{eq:u1=u2}) as an algebraic variety over
$\F _q$.

We restrict ourselves to looking for solutions among the matrices
$x,y$ of the special form:
$$
x=\left( \array{cr} t & -1 \\ 1 & 0
\endarray
\right)\,, \qquad y=\left( \array{cc} 1 & b \\ c & 1+bc
\endarray
\right)\,.
$$

We then study the arising variety $C_w\subset \A ^3$ (with affine
coordinates $b,c,t$), defined by the matrix equation $u_1=u_2$,
with the help of \textsc{Singular} package. The first striking observation
is the following {\em dimension jump}: there are four initial
words $w$ among about 10000 shortest ones such that the dimension
of $C_w$ is one (and not zero as one might expect and as it occurs
for most words $w$). Here are these four words: $w_1=x\min yxy\min
x, w_2=x^{-2}y\min x, w_3=y^{-2}x\min y, w_4=xy^{-2}x\min yx\min
$. Note that if $w=w_i$ ($i=1,2,3$), any solution to $u_1=u_2$ is
automatically non-trivial, thus we only consider these three words
(cf. hypothesis 1 of Proposition \ref{prop:Engel}).

For these ``good'' words, we proceed as follows. As explained
above, the main idea is to apply the Lang--Weil bound for the
number of rational points on a variety defined over a finite
field. It turns out that in the $PSL(2)$ case for our purposes it
is enough to use the classical Hasse--Weil bound (in a slightly
modified form adapted for singular curves, cf.
\cite[Th.~3.14]{FJ}, \cite{AP}, \cite{LY}).

\begin{lemma} \label{lem:HW}
Let $C$ be an absolutely irreducible projective algebraic curve
defined over a finite field $\F _q$, and let $N_q=\#C(\F _q)$
denote the number of its rational points. Then $|N_q-(q+1)|\leq
2p_a\sqrt{q}$, where $p_a$ stands for the arithmetic genus of $C$
(in particular, if $C$ is a plane curve of degree $d$,
$p_a=(d-1)(d-2)/2$).
\end{lemma}

In fact, we need an affine version of the lower estimate of Lemma
\ref{lem:HW} (cf. \cite[Th.~4.9, Cor.~4.10]{FJ}) based on the fact
that the affine curve $C$ has at most $\deg (\overline C)$
rational points less than the projective closure $\overline{C}$.

\begin{corol} \label{cor:HW}
Let $C \subset \A^n$ be an absolutely irreducible affine curve
defined over the finite field $\F_q$ and $\overline{C} \subset
\P^n$ the projective closure. Then the number of $\F_q$-rational
points of $C$ is at least $q+1 - 2p_a \sqrt{q}-d$ where $d$ is the
degree and $p_a$ the arithmetic genus of $\overline{C}$.
\end{corol}

To apply Lemma \ref{lem:HW} (or Corollary \ref{cor:HW}) we have to
compute the arithmetic genus of the curve $C_w$ (or the degree of
some plane projection of $C_w$) and to prove that the curve is
absolutely irreducible.


The case $G=\PSL (3,3)$ is easily settled by full search. For
example, for $w=w_2$ we find a solution to $u_1=u_2$ given by the
images in $G$ of the following matrices:
$$
x=\left( \array{ccc} 0 & 0 & 1 \\ 0 & 1 & 0 \\ 1 & 0 & 1
\endarray \right)\;, \qquad
y=\left( \array{ccc} 2 & 0 & 2 \\ 0 & 1 & 1 \\ 2 & 1 & 1
\endarray
\right)\,.
$$

\subsection{Suzuki case} \label{subsec:Suzuki}
The last remaining case $G=\Sz (q)$ is the most complicated one, particularly
from the computational side. Large group orders require heavy computations (we
used MAGMA for the group-theoretic part and \textsc{Singular} for the
algebraic-geometric one). Moreover, there are even deeper reasons making the
Suzuki case especially difficult. Although both $\PSL (2,q)$ and $\Sz (q)$ are
groups of Lie type of rank 1, and their algebraic structure is very similar,
geometric properties of equations under consideration are significantly
different. Namely, in the $\PSL (2,p)$ case the algebraic variety given by
equation (\ref{eq:u1=u2}) is in fact defined over the ring of integers $\Z$,
and the corresponding variety over $\F _p$ is obtained by reducing modulo $p$;
in particular, the degree is the same for all $p$, and we thus are able to
apply the Lang--Weil estimates for $p$ big enough.

In the Suzuki case the situation is quite different. The group $\Sz (q)$ is
defined with the help of a Frobenius-like automorphism, and hence the standard
matrix representation for $\Sz (q)$ (see below) contains entries depending on
$q$. Therefore, the degree of the resulting variety $u_1=u_2$ depends on $q$
(and grows with growing $q$) which prevents from direct application of the
Lang--Weil estimates. Fortunately, there is a way out described below.

Our strategy is essentially the same: we start with screening for
``good'' initial words $w$ such that the equation $u_1=u_2$ has a
solution in $\Sz (q)$ for $q=8, 32, 128, \dots$, and the number of
solutions grows with growth of $q$ (this last condition should be
emphasized because it gives hope for using algebraic-geometric
machinery). To be more precise, we use the standard embedding of
$\Sz (q)$ into $\GL (4,q)$ (see \cite{HB}) and look for solutions
among the matrices of the special form:

$$
x=\left( \array{clll}
a^{2+\theta }+ab+b^{\theta } & b & a & 1 \\
a^{1+\theta }+b & a^{\theta } & 1 & 0 \\
a & 1 & 0 & 0 \\
1 & 0 & 0 & 0
\endarray
\right)\,,
$$
$$
 y=\left( \array{clll}
c^{2+\theta }+cd+d^{\theta } & d & c & 1 \\
c^{1+\theta }+d & c^{\theta } & 1 & 0 \\
c & 1 & 0 & 0 \\
1 & 0 & 0 & 0
\endarray
\right)\,.
$$
Here $a,b,c,d\in\F _q$, $q=2^n$, $n$ odd, and $\theta$ stands for
the automorphism of $\F _q$ with $\theta ^2=2$. The following
amazing fact is crucial for us: ``good'' initial words $w$ are
those for which the $\F _q$-variety $V_n$, corresponding to the
equation $u_1(x,y)=u_2(x,y)$ with $x,y$ chosen as above, is, in a
certain sense, $\theta$-invariant. To be more precise, if $w$ is a
``good'' word, we can construct a ``universal'' model $V$ for all
varieties $V_n$ in the following sense: $V$ is an $\F _2$-variety
carrying an operator $\alp\colon V\to V$ such that fixed points of
$\alp ^n$ correspond to rational points of $V_n$ (note that a
similar operator appears in \cite[Section 11]{DL}). We are thus
reduced to the proof of the existence of a fixed point of $\alp
^n$ for every odd $n$. To prove that, we use a Lefschetz trace
formula resulting from Deligne's conjecture (established by Pink
and Fujiwara) in order to get an estimate of Lang--Weil type for
the number of fixed points (see Section \ref{sec:Suzuki} for more
details). These estimates guarantee the existence of a solution
for $q$ big enough (small values of $q$ are checked directly).
Thus our final step here is another screening for initial words
satisfying the above invariance condition. Luckily enough, among
these words we find the word $w_2=x^{-2}y\min x$ which is also
good for the $\PSL (2)$ case and for the $\PSL (3,3)$ case. This
establishes Theorem \ref{conj:main}.

\begin{remark} \label{rem:deep} We have no conceptual
explanation of the computational phenomenon of ``good'' words
(which yield the dimension jump and, in the Suzuki case,
additional symmetries). Here is another related observation:
``good'' initial words $w$ correspond to ``deep minima'' of the
length function
$$
l(n)=\textrm{min} \, length(w(x,y)\cdot [u_2(x,y)]\min ).
$$
where minimum is taken over all words $w(x,y)$ of length $n$.
\end{remark}

\begin{remark} Note that the word $w=w_2=x^{-2}y\min x$ satisfies the
hypotheses of Proposition \ref{prop:Engel} and thus this
proposition is true with $w_2$ chosen as the initial term of the
Engel--like sequence.
Note also that with $u_1 = w_2$ we have $u_1 = 1$ if and only if $y =
x^{-1}$.  Hence, for $w = w_2$, $u_1(x,y) = u_2(x,y)$ has a non--trivial
solution if and only if it has a solution with $y \not= x^{-1}$.
\end{remark}

\section{$\PSL (2)$ case: details} \label{sec:PSL}

We shall prove in this section:

\begin{proposition}\label{prop1.1}
If $q = p^k$ for a prime $p$ and $q \neq 2,3$, then there are
$x,y$ in $\PSL(2,\F_q)$ with $y \neq x^{-1}$ and $u_1(x,y) =
u_2(x,y)$.
\end{proposition}

The proof will use some explicit computations with the following
matrices. Let $R$ be a commutative ring with identity and define
\[
x(t) =
\begin{pmatrix}
  t & -1\\1 & 0
\end{pmatrix}, \qquad y(b,c) =
\begin{pmatrix}
  1 & b\\
c & 1+bc
\end{pmatrix} \in \SL(2,R)
\]
for $t,b,c \in R$.

\begin{remark} \label{rem:PSL}\mbox \\

\begin{enumerate}
\item We have
\[
x(t)^{-1} =
\begin{pmatrix}
  0 & 1\\-1 & t
\end{pmatrix}, \qquad y(b,c)^{-1} =
\begin{pmatrix}
  1+bc & -b\\
-c & 1
\end{pmatrix}
\]
for $t,b,c$.

\item For any $t,b,c \in R$ we have $y(b,c) \neq x(t)^{-1}$,
even for the images of $x(t)$ and $y(b,c)$ in $\PSL (2,R)$.

\item The equation $u_1=u_2$ is equivalent to
$x^{-1} yx^{-1}y^{-1} x^2 = yx^{-2}y^{-1}xy^{-1}$; for $t,b,c \in
R$ we put $x = x(t)$, $y = y(b,c)$, and write
\[
x^{-1} yx^{-1}y^{-1} x^2 - yx^{-2}y^{-1}xy^{-1} =
\begin{pmatrix}
n_1(t,b,c) & n_2(t,b,c)\\
n_3(t,b,c) & n_4(t,b,c)
\end{pmatrix}\,.
\]
\end{enumerate}

\end{remark}

Let $I = \langle n_1, n_2, n_3, n_4\rangle \subseteq \Z[b,c,t]$ be
the ideal generated by the entries of the matrix.

Using \textsc{Singular} we can obtain $I$ as follows:\footnote{A
  file with all \textsc{Singular} computations can be found at \\
  \texttt{http://www.mathematik.uni-kl.de/\symbol{126}pfister/SolubleGroups}.}
\begin{verbatim}
LIB"linalg.lib";  option(redSB);
ring R = 0,(c,b,t),(c,lp);
matrix X[2][2] = t, -1,
                 1,  0;
matrix Y[2][2] = 1, b,
                 c, 1+bc;
matrix iX = inverse(X);     matrix iY = inverse(Y);
matrix M=iX*Y*iX*iY*X*X-Y*iX*iX*iY*X*iY;  ideal I=flatten(M);   I;
I[1]=c2b3t2+c2b2t3-c2b2t2+c2b2t+c2b2-c2bt3+2c2bt2+c2bt-c2t2+c2t+c2-cb3t
     +cb2t2+cb2t+cbt3-cbt2+cbt+2cb-ct3+ct2+2ct+c-b2t+bt+1
I[2]=c2b2t+c2bt2+c2t-cb3t2-cb2t3-cb2t-cb2-2cbt2+cbt+ct2-ct-c+b3t-bt-b-1
I[3]=c3b3t2+c3b2t3+c3b2t+2c3bt2+c3t-c2b3t-c2b2t3+2c2b2t2+c2b2t-c2bt4
     +2c2bt3+c2bt2+c2bt-c2t3+2c2t2+c2t+2cb2t2-2cb2t-cb2+cbt2+cbt+cb-ct4
     +ct3+3ct2-c-b2t+bt2-bt-b+1
I[4]=-c2b3t2-c2b2t3+c2b2t2-c2b2t+c2bt3-2c2bt2+c2t2-c2t+cb3t-cb2t2-2cb2t
     -cbt3-cbt+ct3-ct2-2ct+b2t+b2-bt-b-t+1
\end{verbatim}

Denote by $C$ the $\F_q$-variety defined by the ideal
$I\F_q[b,c,t]$.

To prove Proposition \ref{prop1.1}, it is enough to prove

\begin{proposition}\label{prop1.2}
Let $q$ be as in Proposition $\ref{prop1.1}$, then the set
$C(\F_q)$ of rational points of $C$ is not empty.
\end{proposition}

The proof is based on the Hasse--Weil estimate (see Corollary
\ref{cor:HW}).

Note that the Hilbert function of $\overline{C}$, $H(t) = dt - p_a
+1$, can be computed from the homogeneous ideal $I_h$ of
$\overline{C}$, hence we can compute $d$ and $p_a$ without any
knowledge about the singularities of $\overline{C}$.  The ideal $I_h$ can be
computed by homogenising the elements of a Gr\"obner basis of $I$ with respect
to a degree ordering (cf.\ \cite{GP3}).

In the following let $q = p^k$ be an arbitrary, fixed prime power and $L$ the
algebraic closure of $\F_q$. To apply Corollary \ref{cor:HW}, we have to prove

\begin{prop} \label{prop:PSL-irr}
$IL[b,c,t]$ is a prime ideal.
\end{prop}

We start with the following

\begin{lemma}\label{lemma1.4}

The following polynomials form a Gr\"obner basis of $IL[c,b,t]$ with respect
to the lexicographical ordering $c > b>t$,
\begin{verbatim}
J[1]=(t2)*b4+(-t4+2t3)*b3+(-t5+3t4-2t3+2t+1)*b2+(t5-4t4+3t3+2t2)*b
      +(t4-4t3+2t2+4t+1)
J[2]=(t3-2t2-t)*c+(t2)*b3+(-t4+2t3)*b2+(-t5+3t4-2t3+2t+1)*b+(t5-4t4+3t3+2t2)
J[3]=(t)*cb+(-t2+2t+1)
J[4]=cb2+(-t2+2t+1)*c+(-t)*b3+(t3-2t2)*b2+(t4-3t3+2t2-t)*b+(-t4+4t3-3t2-2t)
J[5]=(t)*c2-cb+(t)*c+(-t2)*b3+(t4-2t3+t)*b2+(t5-3t4+t3+2t2-2t-1)*b
      +(-t5+3t4-4t2+t)
\end{verbatim}
\end{lemma}

\begin{proof}
The Gr\"obner basis can be computed in \textsc{Singular} as follows (in
characteristic $0$):

\begin{verbatim}
ideal J=std(I);
\end{verbatim}

We want to verify ``by hand'' that $J$ is indeed a Gr\"obner basis
for each $q$. Indeed, given some intermediate data obtained with
the help of a computer, the truth of the lemma can be verified
without computer. We first show that $I$ and $J$ generate the same
ideal.

\begin{verbatim}
matrix M=lift(I,J); M;

M[1,1]=b2t4+2bt3-t4-t3+3t2+2t            M[1,2]=bt4-t4+2t3+t2+t
M[1,3]=t M[1,4]=-bt3+t3-2t2-t
M[1,5]=-cb2t3-cbt4-ct3+ct+b2t2-bt4-bt+t6-t5-t4+t3-2t2-t
M[2,1]=-cbt5+cbt4+2cbt3-cbt2-cbt-ct4+ct2+ct-bt3+bt-t5+t4+2t3-2t2-2t
M[2,2]=cbt4-cbt2-cbt+ct2+ct+t4-t3-t2
M[2,3]=-cbt-t M[2,4]=-cbt3+2cbt-ct-t3+t2+2t-1
M[2,5]=-cbt4+2cbt2-cbt+ct5-2ct4+ct3+ct2-3ct-bt+t5-3t4+t3+3t2-3t+1
M[3,1]=bt4-bt3-2bt2+bt+b+t3-t-1             M[3,2]=-bt3+bt+b-t-1
M[3,3]=b M[3,4]=bt2-2b+1
M[3,5]=bt5-2bt4+bt3+2bt2-3bt+b-t5+2t4-2t2+3t
M[4,1]=cbt4-cbt3-2cbt2+cbt+cb+ct3-ct-c+b2t4+bt5-bt4+2bt2+bt+t4-2t3+3t+2
M[4,2]=-cbt3+cbt+cb-ct-c+bt2+bt-t4+t3+t2+t+1
M[4,3]=cb+bt+t+1 M[4,4]=cbt2-2cb+c-2bt+t3-t2-t-2
M[4,5]=-cb2t3+cbt5-3cbt4+cbt3+2cbt2-4cbt+cb-ct3-ct2+2ct+c+b2t2-bt2+t5-t4-3t+1
\end{verbatim}

This implies that over $\Z$ and, hence, over each $\F_q$
\[
J[k] = \sum^4_{\ell=1} M[\ell,k] \cdot I[\ell]\,, \quad k = 1, \dots, 5.
\]
\begin{verbatim}
matrix N= lift(J,I); N;

N[1,1]=-cb+c-1 N[1,2]=b  N[1,3]=-c2b-1  N[1,4]=cb-c+b+t N[2,1]=cb2-cb+b-1
N[2,2]=-b2+1       N[2,3]=c2b2+c+b-1          N[2,4]=   -cb2+cb-b2-bt+t-1
N[3,1]=cb2t+cbt+2cb+bt+t+2                    N[3,2]=cb+ct-b2t-2bt-2b-t-1
N[3,3]=c2b2t+2c2bt+2c2b+2c2t-cb2+ct-c+bt-b+t+2
N[3,4]=-cb2t-cbt-2cb+c+b2-bt-2b-3t+1            N[4,1]=c       N[4,2]= -1
N[4,3]=0  N[4,4]=0       N[5,1]=-1  N[5,2]=1    N[5,3]=c-1     N[5,4]=t-1
\end{verbatim}
In the same way this implies that
\[
I[k] = \sum^5_{\ell=1} N[\ell,k] \cdot J[\ell ]\,,\quad k = 1, \dots, 4.
\]
We proved that the polynomials $J[1], \dots, J[5]$ generate the ideal $I$.

To show that $J[1], \dots, J[5]$ is a Gr\"obner basis, we use
Buchberger's criterion (cf.\ \cite{GP3}, Theorem 1.7.3).   To see
this for any $q$, we can use the same trick as above.   Let $s =
\spoly(J[i], J[j])$, $i < j$, the $s$--polynomial of $J[i]$ and
$J[j]$.   We have to show that the normal form of $s$ with respect
to $J[1], \dots, J[5]$ is $0$.   We apply \texttt{lift(s,J);} in
\textsc{Singular} and use the result to check by hand that $s$ is
a linear combination of $J[1], \dots, J[5]$ in all
characteristics. As this is similar to above, we dispense with the
output.
\end{proof}

\begin{lemma}\label{lemma3.6}
Let
 \begin{align*}
 f_1 & = t^2 b^4 - t^3(t-2)b^3+(-t^5+3t^4-2t^3+2t+1)b^2
 +t^2(t^2-2t-1)(t-2)b+(t^2-2t-1)^2\,\\
 f_2 & = t(t^2-2t-1)c + t^2b^3 + (-t^4 + 2t^3)b^2+(-t^5+3t^4-2t^3+2t+1)b
 +(t^5-4t^4+3t^3+2t^2)\,,\\
 h & = t(t^2 - 2t-1)\,.
\end{align*}
Then the following holds for any prime power $q$.

\noindent (1)\; $\{f_1, f_2\}$ is a Gr\"obner basis of $IL(t)[b,c]$ with
respect to the lexicographical ordering $c > b$;

\noindent (2)\; $I : h = I$;

\noindent (3)\; $IL(t)[b,c] \cap L[t,b,c] = \langle f_1, f_2\rangle : h^2 = I$.
\end{lemma}

\begin{proof}
  Because $J$ is a Gr\"obner basis of $I$ with respect to the lexicographical
  ordering $c > b > t$, $J$ is a Gr\"obner basis of $IL(t)[b,c]$ with respect
  to the lexicographical ordering $c>b$ (cf.\ \cite{GP3}, Chapter 4.3).  But
  $J[1] = f_1$ and $J[2] = f_2$ and, considered in $IL(t)[b,c]$, the leading
  monomials of $f_1$ and $f_2$ generate already the leading ideal of
  $IL(t)[b,c]$. This shows (1).

(3) is a consequence of (2) because $IL(t)[b,c] \cap L[t,b,c] = \langle
    f_1, f_2\rangle : h^\infty$, see \cite[Prop. 4.3.1]{GP3}, and
    $h^2 I \subset \langle f_1, f_2\rangle$ that we shall see now.
\begin{verbatim}
M=lift(ideal(J[1],J[2]),h^2*I); M;

M[1,1]=(-t5+2t4+t3)*cb2+(-t6+3t5-2t4+3t2+t)*cb+(t6-4t5+2t4+4t3+t2)*c
       +(-t3)*b2+(t3)*b+(t2)
M[1,2]=(-t4+2t3+t2)*cb+(-t5+2t4+t3)*c+(t5-2t4)*b2+(t6-3t5+2t4-t3-3t2-t)*b
       +(-t6+4t5-3t4-2t3-t2)
M[1,3]=(-t5+2t4+t3)*c2b2+(-t6+2t5+2t3+t2)*c2b+(-2t5+4t4+2t3)*c2+(t4-t3-t2)*cb2
       +(-t5+2t4)*cb+(-t6+2t5+t4-t2)*c+(-t3)*b2+(t4-t3-t2)*b+(t2)
M[1,4]=(t5-2t4-t3)*cb2+(t6-3t5+2t4-t3-t2)*cb+(-t6+4t5-3t4-2t3)*c+(t3+t2)*b2
       +(-t3-t2)*b+(-t3+t2)
M[2,1]=(t5-2t4-t3)*cb3+(t6-3t5+2t4-3t2-t)*cb2+(-t6+4t5-2t4-4t3-t2)*cb
       +(-t5+3t4-3t2-t)*c+(t3)*b3+(-t3)*b2+(-t2)*b
M[2,2]=(t4-2t3-t2)*cb2+(t5-2t4-t3)*cb+(t4-2t3-t2)*c+(-t5+2t4)*b3
       +(-t6+3t5-2t4+t3+3t2+t)*b2+(t6-4t5+3t4+2t3+t2)*b+(t5-4t4+2t3+4t2+t)
M[2,3]=(t5-2t4-t3)*c2b3+(t6-2t5-2t3-t2)*c2b2+(2t5-4t4-2t3)*c2b+(t4-2t3-t2)*c2
       +(-t4+t3+t2)*cb3+(t5-2t4)*cb2+(t6-2t5-t4+t2)*cb+(-t4+2t3+t2)*c
       +(t3)*b3+(-t4+t3+t2)*b2+(-t2)*b
M[2,4]=(-t5+2t4+t3)*cb3+(-t6+3t5-2t4+t3+t2)*cb2+(t6-4t5+3t4+2t3)*cb+(t5-3t4+t3
       +t2)*c+(-t3-t2)*b3+(t3+t2)*b2+(t3-t2)*b
\end{verbatim}
This implies
\[
h^2 \cdot n_i = M[1,i] \cdot f_1 + M[2,i] \cdot f_2,\quad i=1,
\dots, 4.
\]
To prove (2) we can use the \textsc{Singular} commands
\begin{verbatim}
poly h=t*(t2-2t-1);
reduce(quotient(I,h),std(I));
_[1]=0
_[2]=0
_[3]=0
_[4]=0
\end{verbatim}
to see that $I : h \subset I$.

If we want to check this by hand, we can use the following method to compute the
quotient (cf.\ \cite{GP3}, 2.8.5):

\noindent If $U = \langle[g_1, 0], \dots, [g_n,0], [h,1]\rangle$ is a
submodule of the free module $L[c,b,t]^2$ and $[0,h_1], \dots, [0,h_r]$ is the
part of the Gr\"obner basis of $U$ (with respect to the ordering $(c,>)$
giving priority to the components (cf.\ \cite[2.3]{GP3})), having the first
component zero, then $\langle g_1, \dots, g_n\rangle : h = \langle h_1, \dots,
h_r\rangle$.
\begin{verbatim}
module N=[J[1],0],[J[2],0],[J[3],0],[J[4],0],[J[5],0],[h,1];
module N1=std(N);
N1;
N1[1]=[0,b4t2-b3t4+2b3t3-b2t5+3b2t4-2b2t3+2b2t+b2+bt5-4bt4+3bt3+2bt2+t4
       -4t3+2t2+4t+1]
N1[2]=[0,ct3-2ct2-ct+b3t2-b2t4+2b2t3-bt5+3bt4-2bt3+2bt+b+t5-4t4+3t3+2t2]
N1[3]=[0,cbt-t2+2t+1]
N1[4]=[0,cb2-ct2+2ct+c-b3t+b2t3-2b2t2+bt4-3bt3+2bt2-bt-t4+4t3-3t2-2t]
N1[5]=[0,c2t-cb+ct-b3t2+b2t4-2b2t3+b2t+bt5-3bt4+bt3+2bt2-2bt-b-t5+3t4-4t2+t]
N1[6]=[t2-2t-1,-cb+t-2]
N1[7]=[b3-b2-bt+2b,cb-4ct2+10ct-c+b5t-b4t-b3t5+3b3t4-4b3t2-4b3t+4b3-b2t6
      +5b2t5-6b2t4+b2t3-3b2t2+4b2t-6b2+bt6-6bt5+14bt4-13bt3+4bt2-7bt-b+t5
      -10t4+27t3-17t2-8t-4]
N1[8]=[cb,-c2+2cb-ct2+2ct-b2+bt2-2bt+t3-3t2+4]
N1[9]=[c2+c+b2-b-t+2,c3b+c2b+4c2-5cb+ct2-3ct+5c+b4t-b3t3+b3t2-b3t+2b2t3
       -4b2t2+b2t+4b2+bt5-bt4-6bt3+5bt2+5bt-5b-t5+2t4+3t3-4t2-2t-1]
\end{verbatim}
We see that in the second component of $N1[1], \dots, N1[5]$ we have exactly
the Gr\"obner basis $J$.

We have to check that $N = N1$ and $N1$ is a Gr\"obner basis.   The last claim
follows again by using Buchberger's criterion (\cite{GP3}, Theorem 1.7.3).
To see that $N = N1$, we compute
\begin{verbatim}
M=lift(N1,N);

M[1,1]=-c+b2t-b2-bt2+2bt+4b-t2+2t  M[1,2]=-b-2t2+3t+3  M[1,3]=0
M[1,4]=b+t2-2t-1 M[1,5]=b+2t2-3t-4 M[1,6]=0
M[2,1]=-b3t-b3-b2t2+b2t+bt2-2bt+8b+t-2
M[2,2]=-b2+b+8 M[2,3]=1 M[2,4]=-4 M[2,5]=b2-b-9 M[2,6]=0
M[3,1]=-b2t-4b2-2bt+2b+t-2        M[3,2]=c-b2t-4b-t+1  M[3,3]=-2
M[3,4]=b2+bt-b+1        M[3,5]=-c2-c+b2t+3b+4          M[3,6]=1
M[4,1]=b3-b2+1 M[4,2]=b2-b        M[4,3]=0             M[4,4]=-2
M[4,5]=-b2+b   M[4,6]=0 M[5,1]=0  M[5,2]=0  M[5,3]=1   M[5,4]=0
M[5,5]=-t-2    M[5,6]=0
M[6,1]=b4-b3t2-b3-b2t3+b2t2-b2t+b2+bt3-2bt2+t2-2t-1
M[6,2]=ct+b3-b2t2-b2-bt3+bt2-bt+b+t3-2t2    M[6,3]=-1
M[6,4]=-c+b2t+bt2-bt-t2+2t    M[6,5]=-b3+b2t2+b2+bt3-bt2-b-t3+t2+t
M[6,6]=t  M[7,1]=2bt+b M[7,2]=2t+1 M[7,3]=0 M[7,4]=-t  M[7,5]=-2t-1
M[7,6]=0  M[8,1]=0     M[8,2]=0    M[8,3]=t M[8,4]=b   M[8,5]=-1
M[8,6]=0  M[9,1]=0     M[9,2]=0    M[9,3]=0 M[9,4]=0   M[9,5]=t
M[9,6]=0
\end{verbatim}
This implies
\[
N[k] = \sum^9_{\ell=1} M[\ell,k] \cdot N1[\ell]\,,\quad k = 1, \dots, 6.
\]
\begin{verbatim}
M=lift(N,N1);

M[1,1]=0 M[1,2]=0 M[1,3]=0 M[1,4]=0 M[1,5]=0 M[1,6]=0 M[1,7]=0 M[1,8]=0
M[1,9]=0 M[2,1]=0 M[2,2]=0 M[2,3]=0 M[2,4]=0 M[2,5]=0 M[2,6]=0 M[2,7]=0
M[2,8]=0 M[2,9]=0 M[3,1]=-b2t3+2b2t2+b2t+t5-4t4+2t3+4t2+t
M[3,2]=-bt3+2bt2+bt        M[3,3]=-t3+2t2+t  M[3,4]=0 M[3,5]=0 M[3,6]=t2-2t-1
M[3,7]=-b2t4+3b2t3-3b2t-b2+bt5-5bt4+7bt3+3bt2-12bt+b+2t6-11t5+16t4+2t3-10t2-3t
M[3,8]=bt2-2bt-b-2t2+5t
M[3,9]=b2t5-5b2t4+5b2t3+2b2t2+2b2t+b2+2bt6-11bt5+18bt4-5bt3-7bt2+3bt-b+t7
       -6t6+13t5-11t4-4t3+14t2-4t+2
M[4,1]=bt4-2bt3-bt2        M[4,2]=t4-2t3-t2  M[4,3]=0 M[4,4]=-t3+2t2+t
M[4,5]=0 M[4,6]=0
M[4,7]=-ct3+2ct2+ct+bt5-4bt4+2bt3+5bt2-bt-b+t5-4t4-2t3+14t2 M[4,8]=-t3+2t2+t
M[4,9]=-ct4+2ct3+2ct+c-b2t4+2b2t3+b2t2+2bt3-4bt2-2bt+t6-4t5+3t4+t3+5t+1
M[5,1]=0 M[5,2]=0 M[5,3]=0 M[5,4]=0 M[5,5]=-t3+2t2+t  M[5,6]=0
M[5,7]=b2t2-2b2t-b2-t4+4t3-2t2-4t-1 M[5,8]=t2-2t-1
M[5,9]=-cbt2+2cbt+cb+b2t3-2b2t2-b2t-t5+4t4-3t3-4t2+5t
M[6,1]=b4t2-b3t4+2b3t3-b2t5+3b2t4-2b2t3+2b2t+b2+bt5-4bt4+3bt3+2bt2+t4
       -4t3+2t2+4t+1
M[6,2]=ct3-2ct2-ct+b3t2-b2t4+2b2t3-bt5+3bt4-2bt3+2bt+b+t5-4t4+3t3+2t2
M[6,3]=cbt-t2+2t+1
M[6,4]=cb2-ct2+2ct+c-b3t+b2t3-2b2t2+bt4-3bt3+2bt2-bt-t4+4t3-3t2-2t
M[6,5]=c2t-cb+ct-b3t2+b2t4-2b2t3+b2t+bt5-3bt4+bt3+2bt2-2bt-b-t5+3t4-4t2+t
M[6,6]=-cb+t-2
M[6,7]=cb-4ct2+10ct-c+b5t-b4t-b3t5+3b3t4-4b3t2-4b3t+4b3-b2t6+5b2t5-6b2t4+b2t3
       -3b2t2+4b2t-6b2+bt6-6bt5+14bt4-13bt3+4bt2-7bt-b+t5-10t4+27t3-17t2-8t-4
M[6,8]=-c2+2cb-ct2+2ct-b2+bt2-2bt+t3-3t2+4
M[6,9]=c3b+c2b+4c2-5cb+ct2-3ct+5c+b4t-b3t3+b3t2-b3t+2b2t3-4b2t2+b2t+4b2+bt5
       -bt4-6bt3+5bt2+5bt-5b-t5+2t4+3t3-4t2-2t-1

\end{verbatim}
This implies
\[
N1[k] = \sum^6_{\ell=1} M[\ell,k] \cdot N[\ell]\,,\quad k=1, \dots,9\,,
\]
and we obtain finally $N = N1$.
\end{proof}

We now continue the proof of Proposition \ref{prop:PSL-irr}. We have
$IL(t)[b,c] \cap L[b,c,t] = \langle f_1, f_2\rangle : h^2 = IL[b,c,t]$.
Therefore, if $IL[b,c,t]$ were reducible, then $IL(t)[b,c]$ would be reducible
too. We are going to prove that this is not the case.

In $L(t)[b,c]$ the polynomial $f_2$ is linear in $c$.  Since $f_1$
does not depend on $c$, we have $L(t)[b,c]/I \cong L(t)[b]/\langle
f_1\rangle$ and, hence, it suffices to prove that the polynomial
$f_1$ is irreducible.

Set $x=bt$, and let $p(x,t) = t^2 f_1(\tfrac{x}{t},t)$, then
\begin{align*}
  p(x,t)  = &\;  x^4 - t^2(t-2) x^3 + (-t^5 + 3t^4 - 2 t^3 + 2t + 1) x^2 +
  t^3(t-2)(t^2-2t-1)x\\
         & +t^2(t^2-2t-1)^2\,.
\end{align*}

To prove that $f_1\in L[t,b]$ is irreducible, it suffices to prove
that $p \in L[x,t]=L[t][x]$ is irreducible.

We have to prove that $p$ has no linear and no quadratic factor
with respect to $x$.

First we prove that $p$ has no linear factor, that is, that $p(x)
= 0$ has no solution in $L[t]$.

Assume that $x(t) \in L[t]$ is a zero of $p(x) = 0$. Then $x(t)
\mid t^2(t^2-2t-1)^2$. If the characteristic of $L$ is not 2, it
is not difficult to see that $x(t)$ cannot contain the square of
an irreducible factor of $t^2(t^2-2t-1)^2$. If the characteristic
of $L$ is 2, it is not possible that $t^2 \mid x(t)$ or $(t+1)^3
\mid x(t)$. Moreover, it is easy to see that the leading
coefficient of $x(t)$ is $(-1)^{\deg (x(t))-1}$.

The following list gives the candidates for a zero of $p(x)$ and
the value of $p(x)$.

If char$(L) > 2$:
\[
\begin{array}{c|c}
x & \text{leading term of } p(x)\\ \hline
-1 & -t^5\\
t & -t^6\\
t-1-\sqrt{2} & \sqrt{2} t^6\\
t-1+\sqrt{2} & \sqrt{2} t^6\\
-t(t-1+\sqrt{2}) & -\sqrt{2} t^8\\
-t(t-1 -\sqrt{2}) & \sqrt{2} t^8\\
-t^2+2t+1 & -t^8\\
t(t^2-2t-1) & -t^{11}.
\end{array}
\]
If char$(L) = 2$:
\[
\begin{array}{c|c}
x & p(x)\\ \hline
1 & t^5+t^3+t^2\\
t & t^6 + t^5 + t^4\\
t+1 & t^5 + t^3\\
t(t+1) & t^7 + t^5\\
(t+1)^2 & t^8 + t^7 + t^6 + t^4 + t^3 + t^2\\
t(t+1)^2 & t^{11} + t^9 + t^8 + t^7 + t^5 + t^4.
\end{array}
\]
This implies that $p(x)$ has no linear factor with respect to $x$
in $L[x,t]$.

Now assume that $p(x) = (x^2 + ax + b) (x^2 + gx + d)$, $a, b, g,
d \in L[t]$.

This implies:

$\begin{array}{crcl}
(1) & bd & = & t^2(t^2-2t-1)^2\\
(2) & ad + bg & = & t^3(t-2)(t^2-2t-1)\\
(3) & d + ag + b & = & -t^5 + 3t^4 - 2t^3 + 2t + 1\\
(4) & a + g & = & -t^2(t-2).
\end{array}$

If $t^2|b$ then, because of (2), we obtain $t^2|a$.   (4) implies
$t^2|g$ and (2) implies $t^3|a$.   (3) implies that $d \equiv 1 +
2t \mod (t^2)$ and (4) implies that $g \equiv 2t^2 \mod(t^3)$.
If char$(L) \not= 2$, we obtain $d = -(t^2 - 2t - 1)$ and $b =
-t^2(t^2 - 2t -1)$, because $(t^2 - 2t - 1)^2 \equiv 1 + 4t \mod
(t^2)$.   If char$(L) = 2$, then $t^3|a$ and $t^3|g$.   (2)
implies that $\tfrac{a}{t^3} \cdot d + \tfrac{g}{t^3} b =
(t+1)^2$.   This implies $(t+1)^2|b$ and $(t+1)^2|d$.   Therefore,
we have in any characteristic $b = -t^2(t^2-2t-1)$ and $d =
-(t^2-2t-1)$.   (3) implies that $ag = -t^3(t-2)^2$. This is a
contradiction to the fact that $t^3|a$ and $t^2|g$.

We showed that $t^2 \nmid b$.  Similarly, we obtain that $t^2\nmid
d$.  This implies that $t|b$ and $t|d$.  If $(t^2 - 2t - 1)^2|b$,
then (2) implies that $t^2-2t - 1|a$.   Let $d =d_1t$ for a
suitable $d_1 \in L$, then (3) implies that $t^2-2t-1 \mid
-t^5+3t^4-2t^3+2t+1 - d_1t$, that is, $d_1 = -1$.   Then $b =
-t(t^2-2t-1)^2$.   Now (3) implies that $ag = -t^4 + 4t^2 + 4t + 1
= -(t^2 -2t -1)(t+1)^2$.

But $t^2-2t-1|a$ and (4) implies that $\deg(a) = 3$ and $\deg(g) =
1$. This implies that $t+1|a$ and $t+1|g$, which is a
contradiction to (4).

Similarly, we obtain that $(t^2-2t-1)^2\nmid d$.  This implies
that $b = b_3 t(t^2-2t-1)$ and $d = \tfrac{1}{b_3} t(t^2-2t-1)$
for a suitable $b_3 \in L$. (3) implies that $\deg(ag) = 5$.
Because of (4), we may assume that $\deg(a) = 3$ and $\deg(g) =
2$.  (4) implies that $a = -t^3 +$ terms of lower degree. (3)
implies that $g = t^2 +$ terms of lower degree.  (4) implies that
$a = -t^3 + t^2 +$ terms of lower degree.  (2) implies that $b_3 =
-1$.  (3) implies that $ag = -t^5 + 3t^4 - 4t^2 + 1$.  Let $a =
t^3 +t^2 + a_1t + a_0$ for suitable $a_1, a_0 \in L$ then, because
of (4), $g = t^2 - a_1 t - a_0$. (3) implies that $a_0^2 = -1$.
Now $-t^5 + 3t^4 - 4t^2+1 = a \cdot g$ implies that $a_0 = 0$,
which is a contradiction.  This proves that $p$ is irreducible, and hence the
proposition is proved.\qed

We can now apply Corollary \ref{cor:HW} to prove Proposition
\ref{prop1.2}.

We compute the Hilbert polynomial $H(t)$ of the projective curve
corresponding to $I_h$, the homogenization of $I$. We obtain $H(t)
= 10t - 11$. The corresponding \textsc{Singular} session is:

\begin{verbatim}
ring S=0,(b,c,t,w),dp;   ideal J=imap(R,J); ideal K=std(J); K;

K[1]=bct-t2+2t+1
K[2]=bt3-ct3+t4-b2t-c2t-2bt2+2ct2-3t3+bc+2t2-t
K[3]=b2c2-bt2+ct2-t3+b2+2bc+c2+2bt-2ct+2t2+2
K[4]=c2t3-ct4+c3t-2c2t2+3ct3-t4-bc2+bt2-2ct2+4t3-2bt+ct-3t2-b-2t
\end{verbatim}
We now compute matrices to represent the generators of $J$ in
terms of the generators of $K$, and vice versa, in order to see
that in any characteristic $KL[b,c,t,w] = JL[b,c,t,w]$. Moreover,
using Buchberger's criterion, it is not difficult to check that
$K$ is a Gr\"obner basis of $IL[b,c,t,w]$ in any characteristic.

\begin{verbatim}
lift(J,K);

_[1,1]=0   _[1,2]=0   _[1,3]=t-1             _[1,4]=-1
_[2,1]=0   _[2,2]=-1  _[2,3]=-bt-t2+2b+3t-2  _[2,4]=b+c-1
_[3,1]=1   _[3,2]=0   _[3,3]=-t+3            _[3,4]=1
_[4,1]=0   _[4,2]=0   _[4,3]=-t3+bt+2t2+c    _[4,4]=-t
_[5,1]=0   _[5,2]=-1  _[5,3]=t-2             _[5,4]=c

lift(K,J);

_[1,1]=-bt3-bct+2bt2+b2-t2+2t+1  _[1,2]=-t3-ct+2t2+b _[1,3]=1
          _[1,4]=t2+c-2t  _[1,5]=t3+ct-2t2-b
_[2,1]=-b2t  _[2,2]=-bt  _[2,3]=0   _[2,4]=b  _[2,5]=bt-1
_[3,1]=0     _[3,2]=0    _[3,3]=0   _[3,4]=0  _[3,5]=0
_[4,1]=0     _[4,2]=0    _[4,3]=0   _[4,4]=0  _[4,5]=0
\end{verbatim}
We homogenise $K$ with respect to $w$ and obtain again a Gr\"obner basis, cf.\
\cite{GP3}, with respect to the lexicographical ordering. Since
the leading ideal is independent of the characteristic, the Hilbert polynomial
is the same in any characteristic.  We compute
\begin{verbatim}
K=homog(K,w); hilbPoly(K);

-11,10
\end{verbatim}
Hence, the Hilbert polynomial is $10t-11$.  From this we obtain the degree $d
= 10$ and the arithmetic genus $p_a = 12$ of the projective closure. Using
Corollary \ref{cor:HW}, we obtain:
\[
N_q \ge q + 1 - 24 \sqrt{q} - 10.
\]
This implies that $C(\F _q)$ is not empty if $q > 593$.

For small $q$, we give a list of points (Tables \ref{fig:1} and \ref{fig:2})
to prove that $C(\F _q)$ is not empty.

Proposition \ref{prop1.2} and, hence, \ref{prop1.1} are proved.

\begin{remark}\label{remark3.6}
  Using the leading terms of $J$, we can even compute the Hilbert polynomial
  without computer.  Hence, once the matrices are computed by the
  \texttt{lift} command and the Gr\"obner bases are given, we can check
  everything by hand, since only simple (although tedious)
  manipulations are necessary.   Therefore, the PSL(2) case can be verified
  without using any computer.   Unfortunately, this will not be the case for
  the Suzuki groups.
\end{remark}


\begin{table}[ht]
\[\footnotesize{
\begin{array}{cccc}
\begin{array}{r|l}
 p & \text{ point in } C(\F_p)\\ \hline
   5   &  (1,2,2)\\
   7   &  (0,1,4)\\
   11   &  (1,9,1)\\
   13   &  (1,1,8)\\
   17   &  (0,7,7)\\
   19   &  (3,2,10)\\
   23   &  (0,11,19)\\
   29   &  (2,12,8)\\
   31   &  (1,18,26)\\
   37   &  (1,25,22)\\
   41   &  (1,4,19)\\
   43   &  (1,15,3)\\
   47   &  (0,2,8)\\
   53   &  (2,16,12)\\
   59   &  (3,33,39)\\
   61   &  (2,21,49)\\
   67   &  (1,11,63)\\
   71   &  (0,18,60)\\
   73   &  (1,44,49)\\
   79   &  (0,17,71)\\
   83   &  (1,54,39)\\
   89   &  (0,19,26)\\
   97   &  (0,10,15)\\
   101   &  (2,1,47)\\
   103   &  (0,23,39)\\
   107   &  (1,61,26)\\
   109   &  (1,69,102)
   \end{array} & \quad
   \begin{array}{r|l}
p & \text{ point in } C(\F_p)\\ \hline
   113   &  (0,37,52)\\
   127   &  (0,10,112)\\
   131   &  (1,14,22)\\
   137   &  (0,5,32)\\
   139   &  (1,19,109)\\
   149   &  (1,87,63)\\
    151   &  (1,99,108)\\
    157   &  (1,22,62)\\
163   &  (1,67,8)\\
   167   &  (0,3,14)\\
   173   &  (1,101,119)\\
   179   &  (1,11,71)\\
   181   &  (1,3,75)\\
   191   &  (0,7,58)\\
   193   &  (0,45,142)\\
   197   &  (1,18,145)\\
   199   &  (0,67,180)\\
   211   &  (1,51,92)\\
   223   &  (5,6,157)\\
   227   &  (1,118,74)\\
   229   &  (3,220,92)\\
   233   &  (0,19,149)\\
   239   &  (1,179,126)\\
   241   &  (0,67,220)\\
   251   &  (3,15,112)\\
   257   &  (3,97,135)\\
\phantom{x} & \phantom{x}
\end{array}
&\quad
\begin{array}{r|l}
p & \text{ point in } C(\F_p)\\ \hline
   263   &  (0,47,154)\\
   269   &  (2,205,73)\\
   271   &  (0,64,97)\\
   277   &  (4,21,7)\\
   281   &  (0,98,150)\\
   283   &  (1,188,250)\\
   293   &  (1,26,270)\\
   307   &  (1,100,10)\\
   311   &  (2,56,162)\\
   313   &  (0,45,194)\\
   317   &  (2,34,146)\\
   331   &  (1,197,323)\\
   337   &  (0,138,312)\\
   347   &  (1,252,267)\\
 349   &  (2,314,255)\\
   353   &  (0,142,187)\\
 359   &  (0,80,20)\\
367   &  (0,28,80)\\
   373   &  (1,82,336)\\
   379   &  (2,9,197)\\
   383   &  (0,149,138)\\
   389   &  (1,27,379)\\
   397   &  (3,271,169)\\
   401   &  (0,48,349)\\
   409   &  (0,50,98)\\
   419   &  (1,121,65)\\
\phantom{x} & \phantom{x}
\end{array} &\quad
\begin{array}{r|l}
p & \text{ point in } C(\F_p)\\ \hline
   421   &  (2,331,151)\\
   431   &  (0,100,189)\\
   433   &  (0,67,228)\\
   439   &  (0,4,22)\\
   443   &  (2,213,143)\\
   449   &  (2,215,286)\\
   457   &  (0,63,378)\\
   461   &  (5,5,267)\\
   463   &  (0,62,204)\\
   467   &  (1,70,461)\\
   479   &  (0,202,293)\\
   487   &  (0,9,92)\\
   491   &  (1,31,439)\\
   499   &  (1,275,40)\\
   503   &  (0,12,158)\\
   509   &  (7,424,256)\\
   521   &  (0,219,250)\\
   523   &  (3,8,369)\\
   541   &  (1,220,80)\\
   547   &  (2,264,122)\\
   557   &  (2,42,261)\\
   563   &  (1,317,485)\\
   569   &  (0,269,369)\\
   571   &  (1,443,422)\\
   577   &  (2,169,514)\\
   587   &  (1,45,229)\\
   593   &  (1,240,5).
\end{array}
\end{array}
}\] \caption{$q=p=5,\dots , 593$}\label{fig:1} 
\end{table}

\begin{table}
\[\footnotesize{
\begin{array}{cc}
\begin{array}{c|c}
n & \text{point in } C(\F_q)\\ \hline
2     &     (a,0,1) \\
3     &     (a,a^2,a^2)\\
4     &     (a^3,a^{12},a^5)\\
5     &     (a^3,a^{20},a^{22})\\
6     &     (a^9,a^9,a^{54})\\
7     &     (a,a^{62},a^{48})\\
8     &     (a,a^{70},a^{200})\\
9     &     (a,a^{191},a^{121}).\\[1.0ex]
\multicolumn{2}{c}{q = 2^n,\; n = 2, \dots, 9}
\end{array}
&\qquad
\begin{array}{c|c}
n & \text{point in } C(\F_q)\\ \hline
2 &          (a,0,a)\\
3 &         (a,a^3,a^{10})\\
4 &         (a,-1,a^{66})\\
5 &         (a^2,a^{10},a^2).\\
  &\\
  &\\
  &\\
  &\\[1.0ex]
\multicolumn{2}{c}{q = 3^n,\; n = 2, \dots, 5}
\end{array}
\end{array}
}\] \caption{}\label{fig:2}
\end{table}

In Table \ref{fig:2},  $a$ denotes a generator of the multiplicative group $\F_q\smallsetminus
\{0\}$.

\bigskip

\section{Suzuki case: details} \label{sec:Suzuki}

\subsection{Universal model $V$} \label{subsec:universal}

We follow the strategy presented in Subsection
\ref{subsec:Suzuki}.

We use the following representation for the Suzuki group Sz$(q)$.
Let $n = 2m+1$ and $q = 2^n$ and consider the automorphism
\[
\pi\colon\F_q \lra \F_q,\quad \pi (a) = a^{2^{m+1}}.
\]
We have $\pi^2(a)=a^2$.

Let
\[
U(a,b) =
\begin{pmatrix}
  1 & 0 & 0 & 0\\
a & 1 & 0 & 0\\
a\pi(a) + b & \pi(a) & 1 & 0\\
a^2\pi(a) + ab + \pi(b) & b & a & 1
\end{pmatrix},
\]

\[
M(c) =
\begin{array}({llll})
  c^{1+2^m} & 0       & 0        & 0\\
  0         & c^{2^m} & 0        & 0\\
  0         & 0       & c^{-2^m} & 0\\
  0         & 0       & 0        & c^{-1-2^m}
\end{array}, \qquad
T =
\begin{pmatrix}
  0 & 0 & 0 & 1\\
0 & 0 & 1 & 0\\
0 & 1 & 0 & 0\\
1 & 0 & 0 & 0
\end{pmatrix}.
\]
Then Sz$(q) = \langle U(a,b), M(c), T \mid a,b,c \in \F_q, c \neq
0\rangle \subset \SL (4,\F_q)$.

To show that $u_1(x,y) = u_2(x,y)$ has a solution with $y \neq
x^{-1}$, we consider the matrices
\[
x = TU(a,b) \text{ and } y = TU(c,d) \text{ in Sz}(q).
\]
Let $V_n$ be the $\F_q$--variety in $\A^4$ defined by the ideal of the
components of the equation
\[
x^{-1} y x^{-1} y ^{-1} x^2 = y x^{-2} y^{-1} xy^{-1}
\]
which is equivalent to $u_1 (x,y) = u_2(x,y)$.   It is easy to see that $y =
x^{-1}$ in $Sz(q)$ if and only if $a,b,c,d$ are all $0$.

Now it is our aim to show that $V_n(\F_q)\smallsetminus \{0\}$ is not empty.

To make the problem independent of $n$, we replace the expressions
$\pi(a)$, $\pi(b)$, $\pi(c)$, $\pi(d)$ by the indeterminates $a_0,
b_0, c_0, d_0$.

Let
\[
S(a,b,a_0, b_0) =
\begin{pmatrix}
  1 & 0 & 0 & 0\\
a & 1 & 0 & 0\\
aa_0 + b & a_0 & 1 & 0\\
a^2a_0 + ab + b_0 & b & a & 1
\end{pmatrix},
\]
then $U(a,b) = S\bigl(a,b,\pi(a), \pi(b)\bigr)$.

Again we consider the matrices $x = TS(a,b,a_0,b_0)$ and $y =
TS(c,d,c_0,d_0)$ and the equation $u_1(x,y) = u_2(x,y)$. This
leads to a system of equations (16 equations in 8 unknowns, see
Subsection \ref{subsec:geom} below) defining an $\F_2$--variety $V
\subset \A^8$ with ideal $I_V$ generated by the 16 polynomials given in
Subsection \ref{subsec:geom}.

Let us describe the relationship between the ``universal'' model
$V$ and the varieties $V_n$ in more detail. Let $R$, resp.\ $R_1$,  be the
polynomial ring over $\F_2$ in the variables $a,b,c,d,a_0,b_0,c_0,d_0$, resp.\
$a, b, c, d$, and
$l_n : R\to R_1$ the ring homomorphism acting as follows:
\begin{align*}
                  a \to a, & \qquad  a_0 \to a^{2^{m+1}}\\
                   b \to b, & \qquad b_0 \to b^{2^{m+1}}\\
                   c \to c, & \qquad c_0 \to c^{2^{m+1}}\\
                   d \to d, & \qquad d_0 \to d^{2^{m+1}}.
\end{align*}
Denote by $I_1=l_n(I_V)\subset R_1$ the image of the ideal $I_V$, and let
$V_1=\Spec (R_1/I_1)$. Since the images of the matrices $S(a,b,a_0,b_0)$ and
$S(c,d,c_0,d_0)$ (viewed over $\F _q$) under the homomorphism $l_n$ take
values in $\Sz (2^q)$, we obtain for the matrices $x$ and $y$ the presentation
from the beginning of this subsection. In other words, the $\F _q$--varieties
$V_1\times _{\F _2}\F _q$ and $V_n$ are isomorphic.

Note here that the existence of an $\F _q$-point on $V$ does not
guarantee the existence of an $\F _q$-point on $V_n$ (the variety
$V_1$ usually has no $\F _2$-points!). However, here are results
of computer experiments.

Computations show that $V_n$ $(n=3, 5, 7, 11, 13)$ are all of
dimension 0; we have $\# V_3(\F _8)=13$, $\# V_5(\F_{2^5})=40$,
$\# V_{11}(\F _{2^{11}})>1000$. So we have a computational phenomenon
of the growing number of solutions for the equation:
$u_1(x,y)=u_2(x,y)$ for the initial word $w(x,y)=x^{-2}y\min x$.
Note that for most other choices of $w$ the situation is quite
different: there is an odd prime $n$ such that $V_n(\F _{2^n})=\{0\}$.

Our next goal is to show that with a ``good'' choice of the
initial word $w$ the variety $V$ carries an additional structure.
To be more precise, define
$$
S_0(a,b,a_0,b_0)=\left( \array{cccc}
1 & 0         & 0    & 0 \\
a_0 & 1         & 0    & 0 \\
a_0a^2+b_0        & a^2 & 1 & 0 \\
{a_0}^2a^2+a_0b_0+b^2 & b_0    & a_0 & 1
\endarray
\right)\,,
$$
$$
S_0(c,d,c_0,d_0)=\left( \array{cccc}
1 & 0         & 0    & 0 \\
c_0 & 1         & 0    & 0 \\
c_0c^2+d_0        & c^2 & 1 & 0 \\
{c_0}^2c^2+c_0d_0+d^2 & d_0    & c_0 & 1
\endarray
\right)\,.
$$

Let $x_0=TS_0(a,b,a_0,b_0)$, $y_0=TS_0(c,d,c_0,d_0)$.

Define $V_0$ to be the affine variety given by the coordinate equations of
$u_1(x_0,y_0)=u_2(x_0,y_0)$ in the form above.  The variety $V_0$ is also
given by 16 equations in 8 unknowns.

Let $W=V\cap V_0$ be the affine variety given by the equations for
$V$ together with those for $V_0$. The variety $W$ is given by 32
equations in 8 unknowns. Denote the corresponding ideal by $I_W$.

Consider the relationship between varieties $V$ and $W$. First of
all, we have

\begin{prop} \label{prop:dim=2}
$$\dim V=\dim W=2\,$$
\end{prop}

\begin{proof} \textsc{Singular} computation by showing
  $\dim\bigl({\tt {std}}(I_V)\bigr) = 2$, $\dim\bigl({\tt {std}}(I_W)\bigr) = 2$.
\end{proof}

Note that with most other choices of $w$ we have $\dim V=0$
(compare with the dimension jump in the $\PSL (2)$ case).

The key point of our approach is the following mysterious
observation:

\begin{prop} \label{prop:V=W}
With the above notation, $V=W$.
\end{prop}

\begin{proof} \textsc{Singular} computation by showing

{\tt reduce}$\bigl(I_W, {\tt std}(I_V)\bigr) = \langle 0 \rangle$,
{\tt reduce} $\bigl(I_V,
  {\tt std}(I_W)\bigr) = \langle 0 \rangle$.
\end{proof}

\begin{remark}
Note that, as above, this statement fails for most other choices of
$w$. A straightforward computer verification shows that this is
the only word of the length smaller than 10 with the property that
$W=V\bigcap V_0$ is equal to $V$. As mentioned above (see Remark
\ref{rem:deep}), more extensive computer experiments show the
following phenomenon: ``good'' words $w$ correspond to ``deep
minima'' of a certain length function.
\end{remark}

Now define an operator $\alpha$ on the affine space $\A^8=\Spec
(R)$ by the rule:
$$ \alpha(a)=a_0, \qquad \alpha(a_0)=a^2, $$
$$ \alpha(b)=b_0, \qquad \alpha(b_0)=b^2, $$
$$ \alpha(c)=c_0, \qquad \alpha(c_0)=c^2, $$
$$ \alpha(d)=d_0, \qquad \alpha(d_0)=d^2. $$

Then $V_0=\alpha(V)$, and the equality $W=V$ means that the
variety $V$ is preserved by the operator $\alpha$. From the
definition of $\alpha$ it follows that
$$\alpha^2 (v)=\alpha (\alpha (v))=(a^2,b^2,c^2,d^2, a_0^2,b_0^2,c_0^2,d_0^2)$$
for $v=(a,b,c,d,a_0,b_0,c_0,d_0)$. This implies that the variety
$V$ not only carries the operator $\alpha$, but the square of this
operator is the (geometric) Frobenius endomorphism $\Fr$. This
very rigid additional condition on $V$ explains the discovered
computational phenomenon on growing number of solutions, and gives
rise to the proof of the existence of rational points on $V_n$.

Let us consider the action of $\alpha$ on $V$ in more detail. We
obtain:
\[
\begin{array}{c}
\alpha^{2m+1}(a,b,c,d,a_0,b_0,c_0,d_0)=\alpha\bigl(\alpha^{2m}(a,b,c,d,a_0,b_0,c_0,d_0)\bigr)\\[1.0ex]
=\alpha\bigl(\Fr^{m}(a,b,c,d,a_0,b_0,c_0,d_0)\bigr)
=\alpha(a^{2^{m}},b^{2^{m}},c^{2^{m}},d^{2^{m}},a_0^{2^{m}},b_0^{2^{m}},c_0^{2^{m}},d_0^{2^{m}})\\[1.0ex]
=(a_0^{2^{m}},b_0^{2^{m}},c_0^{2^{m}},d_0^{2^{m}},
a^{2^{m+1}},b^{2^{m+1}},c^{2^{m+1}},d^{2^{m+1}})\,.
\end{array}
\]

Thus the set of fixed points of $\alpha^{2m+1}$ looks as follows:
\[
\begin{array}{c}
 \Fix \ \alpha^{n}=\Fix \
\alpha^{2m+1}=\{(a,b,c,d,a_0,b_0,c_0,d_0)\mid\\[1.0ex]
 a=a_0^{2^{m}},b=b_0^{2^{m}}, c=c_0^{2^{m}}, d=d_0^{2^{m}}, a_0=a^{2^{m+1}},
b_0=a^{2^{m+1}}, c_0=c^{2^{m+1}}, d=d_0^{2^{m+1}}\}\, .
\end{array}
\]
Therefore, a point $(a,b,c,d, a_0,b_0,c_0,d_0)$ is a fixed point of
$\alpha^{2m+1}$ if and only if $$
a=(a^{2^{m+1}})^{2^{m}}=(a^{2^{2m+1}})=a^{2^n},
$$
$$
a_0={(a_0^{2^{m}})}^{2^{m+1}}=(a_0^{2^{2m+1}})=a_0^{2^n},
$$
and the same formulas are valid for $b$, $c$, $d$, $b_0$, $c_0$,
$d_0$.

We conclude that the set of fixed points of the operator $\alpha^n$ over
$\F_q$ coincides with set of $\F_q$-points of $V_n$. Indeed, a point
$u=(a,b,c,d)$ is an $\F_q$-point of $V_n$, $q=2^n$, if and only if
$u^q=(a^q,b^q,c^q,d^q)=(a,b,c,d)$. Since $q=2^n$, we get $a=a^{2^n},b=b^{2^n},
c= c^{2^n}, d=d^{2^n}$. Moreover, the equality $\alpha(a)=a_0$ implies that
the same formulas are valid for $a_0,b_0,c_0,d_0$. Thus the point
$(a,b,c,d,a_0,b_0,c_0,d_0)$ belongs to the set of fixed points $\Fix \
(\alpha^n)$.  Conversely, if $(a,b,c,d,a_0,b_0,c_0,d_0)\in \Fix \ (\alpha^n)$,
then according to the definition of $l_n$ the point
$l_n(a,b,c,d,a_0,b_0,c_0,d_0)\in V_n$ is a rational point over $\F_q$.

The idea of the observation above is to use $V\subset \A^8$ as a universal
model, and replace the calculations in each specific variety $V_n\subset\A^4$
by calculations in $V$.   By the above, $V_n$ has a non--zero rational point
over $\F_q$, $q = 2^n$, if and only if $\alpha^n$ has an non--zero rational
fixed point on $V$.

To sum up, we obtained the following reduction.

\begin{theorem} \label{th:universal}
  Suppose that for every odd $n>1$ the operator $\alpha ^n$ has a non--zero
  $\F_q$--rational fixed point on the variety $V$. Then the equation $u_1=u_2$
  has a non--trivial solution in $\Sz (q)$ for every $q=2^n$.
\end{theorem}

\subsection{Geometric structure of $V$} \label{subsec:geom}

Denote $a_0=v$, $b_0=w$, $c_0=x$, $d_0=y$. With this notation, the
ideal $I$ defining $V$ is obtained with \textsc{Singular} as
follows:

\begin{verbatim}
ring A=2,(a,b,c,d,v,w,x,y),dp;

matrix S1[4][4] =1           0, 0, 0,
                 a,          1, 0, 0,
                 av+b,       v, 1, 0,
                 a2v+a*b+w,  b, a, 1;

matrix S2[4][4] =1,          0, 0, 0,
                 c,          1, 0, 0,
                 cx+d,       x, 1, 0,
                 c2x+c*d+y,  d, c, 1;

matrix T[4][4] = 0, 0, 0, 1,
                 0, 0, 1, 0,
                 0, 1, 0, 0,
                 1, 0, 0, 0;

matrix X=T*S1; matrix Y=T*S2;
matrix iX = inverse(X); matrix iY = inverse(Y);

matrix M=iX*Y*iX*iY*X*X-Y*iX*iX*iY*X*iY;
ideal I=flatten(M); I;
I[1]=a3cv2+a3cvx+a2bcv+a2bcx+acv2+acvw+acvx+acwx+adv+bcx
I[2]=a5v2+a3c2vx+a3bcv+a4dv+a3cdv+a2bc2x+a3b2+a2b2c+a3bd+a2bcd
     +a3v2+ac2vx+ac2wx+a3vy+abcv+a2dv+acdv+abcw+a2dw+acdw+abcx
     +a2by+a3+ab2+a2c+abd+av2+avw+aw2+cvx+avy+awy+a
I[3]=a4cv3+a5v2x+a4cv2x+a2c3v2x+a3c2vx2+a2c3vx2+a4bv2+a2bc2v2
     +a3cdv2+a2c2dv2+a3bcvx+abc3vx+a4dvx+a2c2dvx+a2bc2x2+abc3x2
     +ab2c2v+a3bdv+abc2dv+a3v3+a3b2x+a3bdx+abc2dx+ac2v2x+c3vwx
     +ac2wx2+c3wx2+a2cv2y+a3vxy+a2cvxy+a2b3+ab3c+a2b2d+ab2cd
     +bc2v2+a2dv2+acdv2+bc2vw+acdvw+c2dvw+a2bvx+abcwx+a2dwx
     +c2dwx+a2bvy+abcvy+a2bxy+abcxy+b2cv+ac2v+abdv+av3+cv3+b2cw
     +abdw+bcdw+av2w+cvw2+a3x+ab2x+a2cx+b2cx+ac2x+abdx+bcdx
     +aw2x+cw2x+ab2y+av2y+cvwy+awxy+cwxy+a2b+b3+abc+b2d+bv2
     +dv2+bvw+bw2+bvy+bwy+ax+b
I[4]=a6v3+a5cv3+a4c2v3+a5cv2x+a3c3v2x+a2c4v2x+a3c3vx2+a2c4vx2
     +a5bv2+a4bcv2+a3bc2v2+a2bc3v2+a5dv2+a2c3dv2+a3bc2vx+abc4vx
     +a4cdvx+a2c3dvx+a2bc3x2+abc4x2+a4b2v+a3b2cv+ab2c3v+a3bcdv
     +abc3dv+a4v3+a3cv3+a4v2w+a3b2cx+a2b2c2x+ab2c3x+a3bcdx
     +abc3dx+a4v2x+a3cv2x+a2c2v2x+ac3v2x+ac3vwx+c4vwx+a2c2vx2
     +ac3vx2+ac3wx2+c4wx2+a4v2y+a3cv2y+a2c2v2y+a3cvxy+a2c2vxy
     +a3b3+a2b3c+a3b2d+a2b2cd+a2bcv2+abc2v2+bc3v2+a3dv2+ac2dv2
     +abc2vw+bc3vw+a2cdvw+c3dvw+abc2vx+a3dvx+ac2dvx+abc2wx
     +bc3wx+a2cdwx+c3dwx+abc2x2+bc3x2+a2bcvy+abc2vy+a2bcxy
     +abc2xy+a2b2v+ab2cv+ac3v+a2bdv+a2v3+c2v3+a2b2w+a2v2w
     +acv2w+a2vw2+acvw2+c2vw2+a2b2x+a3cx+ab2cx+a2c2x+b2c2x
     +ac3x+a2bdx+abcdx+a2v2x+c2v2x+a2vwx+acvwx+acw2x+a2b2y
     +a2v2y+acv2y+acvwy+c2vwy+a2vxy+acvxy+acwxy+c2wxy+b3c+a3d
     +a2cd+b2cd+adv2+cdv2+abvw+abw2+bcw2+adw2+cdw2+abvx+bcvx
     +advx+cdvx+abwx+bcwx+bcvy+abxy+bcxy+a2v+b2v+acv+bdv+acw
     +bdw+w3+a2x+a2y+b2y+w2y+ab+bc+ad+w
I[5]=a4cv2+a3c2v2+a3c2vx+a2c3vx+a3cdv+a2c2dv+a2bc2x+abc3x+a2b2c
     +ab2c2+a2bcd+abc2d+a3v2+ac2v2+ac2vw+a2cvx+ac2vx+ac2wx+c3wx
     +a2cvy+a2bv+abcv+bc2v+acdv+bc2w+acdw+c2dw+abcx+bc2x+abcy
     +a2c+b2c+ac2+bcd+av2+avw+cw2+avx+cwx+cwy+c
I[6]=a4vx+a3cvx+a3bv+a2bcv+abc2v+a3dv+a3bx+a2bcx+abc2x+bc3x
     +a2b2+ab2c+b2c2+a2bd+abcd+bc2d+c2v2+a2vx+acvx+c2vx+a2wx
     +acwx+abv+bcv+adv+cdv+abw+bcw+adw+abx+bcx+bcy+a2+ac+c2
I[7]=a3cv2x+a2c2v2x+a4vx2+a2c2vx2+a2bcv2+a2cdv2+a2bcvx+abc2vx
     +a3dvx+a2cdvx+a3bx2+abc2x2+a2b2v+a3cv+ab2cv+a2c2v+a2bdv
     +abcdv+a2c2x+ac3x+a2bdx+abcdx+a2v2x+acv2x+acvwx+c2vwx
     +a2wx2+c2wx2+ab3+a2bc+abc2+ab2d+a2cd+ac2d+abv2+bcv2+adv2
     +bcvw+cdvw+bcvx+adwx+cdwx+abx2+bcx2+a2v+acv+c2v+bdv+v3+b2w
     +acw+bdw+a2x+b2x+c2x+v2x+acy+ab+bc+cd+v
I[8]=a5v2x+a2c3v2x+a4cvx2+a2c3vx2+a4bv2+a3bcv2+a2bc2v2+a4dv2
     +a3cdv2+a2c2dv2+a2bc2vx+abc3vx+a3cdvx+a2c2dvx+a3bcx2
     +abc3x2+a2b2cv+ab2c2v+a2bcdv+abc2dv+a3b2x+ab2c2x+a2bcdx
     +abc2dx+a3v2x+ac2v2x+a2cvwx+c3vwx+a3vx2+ac2vx2+a2cwx2
     +c3wx2+a2b3+a2b2d+bc2v2+a2dv2+acdv2+abcvw+bc2vw+acdvw
     +c2dvw+bc2vx+a2dvx+acdvx+abcwx+bc2wx+acdwx+c2dwx+a2bx2
     +abcx2+a3v+ab2v+b2cv+bcdv+av3+cv3+a3x+ab2x+c3x+abdx+bcdx
     +aw2x+cw2x+avx2+cvx2+b3+abc+bc2+a2d+b2d+acd+c2d+bv2+bvw
     +bw2+dw2+bvx+bwx+av+cv+aw+cw+cy+b
I[9]=a5v3+a3c2v3+a5v2x+a4cv2x+a3c2v2x+a2c3v2x+a4cvx2+a2c3vx2
     +a3cdv2+a2c2dv2+a3bcvx+abc3vx+a4dvx+a3cdvx+a3bcx2+abc3x2
     +a3b2v+ab2c2v+a3bdv+abc2dv+a3d2v+a2cd2v+a2cv3+ac2v3
     +a2cv2w+ac2v2w+a3b2x+ab2c2x+a3bdx+a2bcdx+a3v2x+ac2vwx
     +c3vwx+a3vx2+a2cvx2+ac2vx2+a2cwx2+c3wx2+a3v2y+a2cv2y
     +a3vxy+a2cvxy+a2b2d+ab2cd+a2bd2+abcd2+a2bv2+bc2v2+a2dv2
     +abcvw+bc2vw+a2dvw+acdvw+c2dvw+a2bvx+abcvx+abcwx+bc2wx
     +a2dwx+acdwx+a2bx2+abcx2+bc2x2+a2bvy+abcvy+a2dvy+a2bxy
     +abcxy+a3v+ab2v+b2cv+ac2v+ad2v+bcdw+ad2w+cd2w+avw2+cvw2
     +a3x+b2cx+ac2x+av2x+avwx+aw2x+cw2x+avx2+awx2+abdy+avwy
     +cvwy+awxy+cwxy+a2d+b2d+acd+bd2+dw2+bvx+bwx+dwx+bx2+dwy
     +cv+cx+d
I[10]=a3bv2+abc2v2+a3bvx+a2bcvx+abc2vx+bc3vx+a2bcx2+bc3x2+a4v
      +a2b2v+a3cv+b2c2v+abcdv+bc2dv+acv3+c2v3+a2b2x+b2c2x+a2bdx
      +abcdx+acvx2+c2vx2+a3b+a2bc+ab2d+b2cd+abd2+bcd2+abv2+bcv2
      +adv2+abvw+bcvw+bcvx+advx+abwx+bcwx+abx2+abvy+bcvy+abxy
      +bcxy+acv+c2v+bdv+d2v+a2w+acw+bdw+a2x+b2x+c2x+bdx+v2x+vx2
      +bdy+ab+ad+cd+x
I[11]=a4v2+a3cv2+a3cvx+ac3vx+a3cx2+ac3x2+a2cdv+ac2dv+a3dx+a2cdx
      +a2b2+ab2c+a2bd+abcd+a2d2+acd2+a2v2+c2v2+a2vw+c2vw+a2vx
      +acwx+c2wx+a2x2+acx2+c2x2+a2vy+acvy+a2xy+acxy+bcv+adv+abw
      +bcw+adw+abx+ady+b2+bd+d2+v2+vx+x2
I[12]=a5v2+a2c3v2+a4cvx+a2c3vx+a3bcv+abc3v+a3bcx+abc3x+a3b2
      +a2b2c+a2cvw+c3vw+a2cvx+c3vx+a2cwx+c3wx+a2cx2+c3x2+bc2v
      +acdv+c2dv+abcx+bc2x+a2dx+acdx+a3+a2c+b2c+abd+bcd+ad2
      +cd2+av2+cv2+avw+cvw+aw2+cw2+cvx+awx+cwx+ax2+avy+cvy
      +axy+cxy+bv+bw+dw+bx+dx+dy+a
I[13]=a6v3+a4c2v3+a6v2x+a5cv2x+a5cvx2+a2c4vx2+a5dv2+a4cdv2
      +a4bcvx+a2bc3vx+a5dvx+a3c2dvx+a4bcx2+abc4x2+a4b2v+a2b2c2v
      +a4bdv+a2bc2dv+a3cd2v+a2c2d2v+a4v3+a2c2v3+a3cv2w+a2c2v2w
      +a4b2x+ab2c3x+a4bdx+a2bc2dx+a3cv2x+ac3v2x+a2c2vwx+a4vx2
      +a3cvx2+a2c2vx2+ac3vx2+a3cwx2+c4wx2+a4vxy+a2b2cd+ab2c2d
      +a2bcd2+abc2d2+a3bv2+abc2v2+a3dv2+ac2dv2+a2bcvw+abc2vw
      +a3dvw+abc2vx+bc3vx+a3dvx+a2cdvx+a2bcwx+abc2wx+bc3wx
      +a3dwx+ac2dwx+a3bx2+a2bcx2+abc2x2+bc3x2+a3bvy+a2bcvy
      +a3dvy+a3bxy+a4v+a3cv+abcdv+bc2dv+a2d2v+acd2v+a2v3+acv3
      +bc2dw+acd2w+c2d2w+a2v2w+acv2w+a2vw2+acvw2+a4x+ab2cx
      +b2c2x+ac3x+a2bdx+abcdx+acv2x+acvwx+a2w2x+acw2x+c2w2x
      +a2vx2+a2wx2+a2b2y+ab2cy+a2bdy+acv2y+a2vwy+acvxy+a2wxy
      +a2vy2+a2bc+abc2+ab2d+a2cd+b2cd+ac2d+abd2+bcd2+abv2+bcv2
      +adv2+abvw+bcvw+cdvw+cdw2+bcvx+advx+abwx+abx2+abvy+bcvy
      +advy+bcwy+adwy+bcxy+aby2+acv+bdv+acw+c2w+bdw+d2w+acx
      +c2x+bdx+a2y+b2y+acy+bdy+w2y+wy2+bc+ad+cd+y
I[14]=a4bv2+a2bc2v2+a4bvx+a3bcvx+a3bcx2+bc4x2+a3b2v+ab2c2v
      +a3bdv+a2bcdv+a2cv3+ac2v3+a3b2x+b2c3x+a3bdx+abc2dx
      +c3v2x+a2cvx2+ac2vx2+c3vx2+ab2cd+b2c2d+abcd2+bc2d2
      +a2dv2+acdv2+c2dv2+a2bvw+abcvw+a2bvx+bc2vx+a2dvx+acdvx
      +a2bwx+abcwx+bc2wx+a2bx2+a2bxy+ab2v+b2cv+ac2v+cd2v+av3
      +cv3+bcdw+a3x+ab2x+c3x+cv2x+avx2+ab2y+b2cy+abdy+cv2y
      +cvxy+a2b+abc+bc2+b2d+acd+c2d+bd2+dv2+dvx+bvy+dvy+bwy
      +by2+aw+ax+ay+cy+d
I[15]=a5v2+a3c2v2+a5vx+a4cvx+a4cx2+ac4x2+a4bv+a2bc2v+a4dv+a3cdv
      +a4bx+abc3x+a4dx+a2c2dx+a2bcd+abc2d+a2cd2+ac2d2+a3v2
      +a2cv2+ac2v2+a3vw+a2cvw+a3vx+c3vx+a3wx+a2cwx+ac2wx+a3x2
      +a2cx2+ac2x2+c3x2+a3xy+a2bv+abcv+a2dv+c2dv+acdw+a2bx+bc2x
      +a2dx+acdx+a2by+abcy+a2dy+ab2+a2c+ac2+abd+bcd+ad2+cd2+cv2
      +cvw+cvx+awx+cwx+ax2+cvy+awy+cxy+ay2+bv+dv+bw+dx+by+dy+c
I[16]=a3cv2+a4vx+a2c2vx+a3cx2+c4x2+a3bv+a2bcv+a3dv+a2cdv+a3bx
      +a2bcx+abc2x+bc3x+a3dx+ac2dx+a2b2+abcd+bc2d+acd2+c2d2
      +a2vw+c2vw+a2vx+a2wx+a2x2+a2xy+abv+cdv+cdw+bcx+aby+bcy
      +ady+a2+b2+ac+c2+bd+d2+w2+wy+y2
\end{verbatim}

To show that $\alpha^n$ has a rational fixed point on $V$, we want to apply
the Lefschetz trace formula, which requires that $V$ is absolutely
irreducible.   This is not the case.
Therefore, we exhibit a subvariety $V^\prime \subset V$ for which we
can show that it is absolutely irreducible.   Then we apply the
Lefschetz trace formula to the non--singular locus of $V^\prime$ which
happens to be affine.

Set $J = I : a^3 x^2$, then $J \supset I$ and $V^\prime := \V (J)
\subset \V (I) = V$.

We shall show that $V'$ is an absolutely irreducible surface.

The ideal $J$ is given as follows:\footnote{Computing $I : a^3x^2$ is not an
  easy task.  However, once $J$ is given, it is much simpler to check $J
  \supset I$, which is all we need.}

\begin{verbatim}
ideal J=quotient(I,a3x2); J;

J[1]=d2+adv+cdv+a2v2+c2v2+abx+bcx+wx+c2x2+vy+xy+c2;
J[2]=a2b+acd+a2cv+aw+a3x+a2cx+ac2x+ay+av+cx;
J[3]=bcw+acvw+w2+a2wx+acwx+b2+bd+d2+abv+bcv+c2v2+bcx+adx+a4
     +a3c+vx+x2+ac+1;
J[4]=adv2+cdv2+d2x+abvx+bcvx+advx+cdvx+vwx+abx2+bcx2+wx2
     +c2x3+v2y+vxy+x2y+ab+cd+acv+c2v+w+a2x+acx+c2x+y;
J[5]=abd+abcv+bc2v+a2dv+dw+avw+cvw+bc2x+c2dx+ac2vx+awx
     +a2cx2+ac2x2+c3x2+by+cxy+dv+av2+cv2+bx+cx2+ac2+a+c;
J[6]=bcd+cd2+a2bv+abcv+a2dv+c2dv+bw+avw+cvw+a2dx+c2dx+c3vx
     +a3x2+a2cx2+ac2x2+by+dy+cvy+axy+bv+dv+cv2+dx+cvx+ax2+
     a3+a+c;
J[7]=a3v2+a2cv2+c2dx+a3vx+ac2vx+a2cx2+ac2x2+c3x2+cxy+cx2;
J[8]=d2v+acv3+c2v3+cdvx+a2vx2+acvx2+a2bc+ac2d+ac3v+acw+a3cx
     +vx2+acy+a2v+acx+v;
J[9]=advx+cdvx+a2v2x+c2v2x+abx2+bcx2+a2vx2+c2vx2+wx2+vxy+c3d
     +a3cv+a2c2v+a3cx+a2c2x+c4x+c2y+cd+a2v+c2v+c2x+y;
J[10]=a2vw+acvw+c2vw+w2+ac2dx+c3dx+a3cvx+ac3vx+acwx+c2wx
     +a3cx2+c4x2+aby+acxy+c2xy+a2v2+acv2+abx+adx+cdx+a2vx
     +acvx+c2vx+a2x2+c2x2+a4+a2c2+v2+1;
\end{verbatim}



We compute a Gr\"obner basis of the ideal
$J3=J\F_2(a,c)[w,y,b,d,x,v]$ with respect to the lexicographical
ordering.

\begin{verbatim}
ring s=(2,a,c),(w,y,b,d,x,v),lp;
ideal J3=std(imap(r,J));J3;
J3[1]=(a8+a6c2+a4c4+a2c6)*v6+(a8+a7c3+a6c2+a5c3+a4c4+a3c7+a2c6
      +ac7)*v4+(a7c3+a6c2+a5c5+a5c3+a3c7+a3c5+a2c6+a2c4+ac9+c6)
      *v2+(ac9+ac5+c8+c4);
J3[2]=(a4c4+a3c7+a3c5+a3c3+a2c8+a2c4+ac7+c4)*x+(a8+a7c+a4c4
      +a3c5)*v5+(a8+a7c+a6c2+a5c3+a4c4+a4c2+a2c6+a2c4)*v3
      +(a4c4+a4c2+a3c7+a3c3+a2c8+a2c6)*v;
J3[3]=(c2+1)*d2+(xc3+xc)*d+(v3xa2+v3xc2+v2a4+v2a3c+v2ac3+v2c4
      +vxa4+vxa3c+vxa2c2+vxc4+x2a2c2+x2a2+x2ac3+x2c2+c4+c2);
J3[4]=(ac5+ac)*b+(v4a2c2+v4a2+v3xac+v2x2c4+v2x2c2+v2a5c+v2a4
      +v2a2c4+v2ac3+vx3ac+vx3+vxa5c+vxa4+vxa3c3+vxa2c4+vxa2c2
      +vxac5+vxac3+vxac+vxc2+vx+x2a3c3+x2a2c2+x2ac5+x2c4+a2c4
      +a2c2+ac3+ac+c4+c2)*d+(v2xa2c3+v2xc5+vx2a3+vx2ac2+va5c2
      +va4c+va3+va2c+vac6+vac4+vac2+vc5+xa5c2+xa4c+xa2c3+xa2c
      +xac6+xac2+xc3);
J3[5]=(c)*y+(va2c+va)*bd+(v2ac+v2c2+x4+c4+1)*b+(v4a3c+v4ac
      +v3xc2+v2x2ac3+v2x2ac+v2a4+v2a3c3+v2ac3+v2c2+vx3c2+vxa4
      +vxa3c3+vxa3c+vxa2c4+vxa2+vxac3+vxac+x2a3c+x2a2c4+x2a2
      +x2ac+x2+ac3+ac+c2+1)*d+(v3x2ac2+v3x2c3+v3a3c4+v3a3c2
      +v3a2c3+v3a2c+v3ac4+v3c3+v2xa3c2+v2xa3+v2xa2c5+v2xc5
      +v2xc+vx4a3+vx4a2c+vx4a+vx4c+vx2a7+vx2a5+vx2a3c2+vx2a2c3
      +vx2a+vx2c3+vx2c+va7c2+va7+va4c+va3c6+va3c4+va3c2+va2c5
      +va2c+vac6+vac4+vc3+vc+x5a+x5c+x3a6c+x3a5+x3a4c+x3a3
      +x3a2c+x3ac2+x3c+xa7c2+xa6c+xa5c2+xa4c3+xa3c6+xa3c4
      +xa3c2+xa3+xa2c5+xa2c3+xa2c+xac6+xac4+xa);
J3[6]=w+(vx+1)*y+(a)*b+(vxc+c)*d+(v3a2+v3ac+v2xa2+v2xc2+vx2a2
      +vx2ac+vx2c2+vx2+va2c2+va2+vac3+v+xa2c2+xa2+xac3+xc2);
\end{verbatim}

\texttt{dim(J3);} returns 0, hence $V^\prime$ is a surface.

Let $f = (a^3 + a^2c^3 + a^2c+ac^4+ac^2+c)(ac+1)(a+c)(c+1)ac$ be the least
common multiple of the leading coefficients of this Gr\"obner basis. Then,
using \textsc{Singular}\footnote{The first equality is a general fact (cf.\
  \cite{GP3}). To see that $\langle J3[1], \dots, J3[6]\rangle : f^\infty =
  J$, it is sufficient to know that $J \supset \langle J3[1], \dots,
  J3[6]\rangle$, $J = J : f$ and that $\langle J3[1], \dots, J3[6]\rangle :
  f^\infty$ is a prime ideal, which we shall see later.  This is,
  computationally,  much easier
  to check than a direct computation.}, we obtain
\[
J3 \cap \F_2 [a,c,w,y,b,d,x,v] = \langle J3[1], \dots,
J3[6]\rangle : f^\infty = \langle J3[1], \dots,
  J3[6]\rangle : f^6 = J.
\]

Since $J:f=J$, no factor of $f$ divides all elements of $J$.
That is why the irreducibility of $J3$ as an ideal of
$\F_2(a,c)[w,y,b,d,x,v]$ implies the irreducibility of $J$.

Furthermore, we compute the vector space dimension over $\F_2(a,c)$ as
\[
\dim_{\F_2(a,c)} \F_2(a,c)[w,y,b,d,x,v]/J3 = 12\,.
\]
Next we show that $J3 \cap \F_2(a,c) [b] = \langle h \rangle$ with
the following polynomial $h$, which we compute directly by elimination (using
\textsc{Singular}).

\begin{verbatim}
poly h=(a18c2+a16+a14c6+a12c4+a10c10+a8c8+a6c14+a4c12)*b12
+(a20c2+a19c5+a18+a17c7+a17c5+a17c3+a16c6+a15c7+a15c5+a15c3
+a14c4+a13c5+a12c10+a11c13+a10c8+a9c15+a9c13+a9c11+a8c14
+a7c15+a7c13+a7c11+a6c12+a5c13)*b10+(a21c5+a20c4+a19c5+a19c3
+a18c2+a17c9+a17c3+a16c6+a16+a15c7+a14c6+a14c4+a14c2+a13c13
+a12c12+a11c13+a11c11+a10c10+a10c6+a9c17+a9c11+a8c14+a8c8
+a7c15+a6c14+a6c12+a6c10+a2c14)*b8+(a24c2+a22c4+a22+a18c6
+a18c4+a17c11+a17c3+a16c8+a15c13+a15c9+a15c7+a15c5+a14c10
+a14c8+a13c15+a13c11+a13c9+a13c7+a12c14+a12c12+a12c8+a11c17
+a11c13+a11c11+a11c5+a10c16+a10c12+a10c10+a10c4+a9c15+a9c13
+a9c11+a9c9+a8c12+a8c6+a7c13+a7c11+a6c14+a6c8+a5c15+a5c13
+a4c10+a3c15+a3c13)*b6+(a26c2+a25c5+a24c4+a24+a23c5+a23c3
+a22c2+a21c9+a21c5+a21c3+a20c8+a20c6+a20+a19c5+a18c10+a18c8
+a18c6+a18c2+a17c13+a17c5+a17c3+a16c6+a16c4+a16+a15c13
+a15c11+a15c9+a15c7+a15c5+a14c12+a14c8+a14c6+a14c4+a13c17
+a13c11+a13c5+a12c14+a12c12+a12c10+a12c8+a12c6+a12c4+a11c11
+a11c9+a11c7+a10c12+a10c10+a9c17+a9c7+a8c12+a8c10+a8c8+a8c4
+a7c11+a6c14+a6c12+a6c6+a5c17+a5c15+a5c13+a5c11+a4c12+a2c10
+c12)*b4+(a27c5+a26c4+a25c7+a25c5+a25c3+a24c6+a24c2+a23c7
+a23c5+a23c3+a22c6+a21c7+a21c5+a21c3+a20c8+a19c13+a19c9+a19c7
+a19c5+a18c12+a18c10+a18c8+a18c6+a18c4+a18+a17c15+a17c13+a17c9
+a17c5+a16c14+a16c12+a16c8+a16c4+a15c15+a15c3+a14c12+a14c10
+a14c6+a14c4+a13c11+a13c5+a12c14+a12c8+a11c13+a11c9+a11c5
+a10c14+a10c12+a10c10+a9c13+a9c11+a9c9+a8c12+a8c10+a7c13+a6c14
+a5c15+a4c14+a4c12+a4c8+a3c15+a3c13+a2c14+a2c10)*b2
+(a26c6+a24c4+a22c6+a20+a18c14+a16c12+a16c4+a16+a14c14+a14c10
+a14c2+a8c12+a8c8+a8c4+a6c14+a6c10+a4c12+a2c14+a2c10+c8);
\end{verbatim}

$h$ is a polynomial of degree 12 with respect to $b$ and therefore
$\dim_{\F_2(a,c)} \F_2(a,c)[b]/(J3\cap \F_2(a,c)[b]) = 12$. Since
$\dim_{\F_2(a,c)} \F_2 (a,c)[w,y,b,d,x,v]/J3$ is also 12, we know that a
lexicographical Gr\"obner basis with respect to $b <v<x<d<y<w$ of $J3$ must
have leading polynomials as follows: $b^{12}, v, x,d,y,w$.\footnote{We do not
  need to compute directly $J3 \cap \F_2(a,c)[b]=\langle h \rangle$ which is
  difficult. Once $h$ is given, it suffices to know that $h$ is irreducible of
  degree 12, $\dim_{\F_2(a,c)} \F_2(a,c)[w,y,b,d,x,v]/J3 = 12$ and $h \in J3$,
  which is much easier to check.}

It follows that the projection
$$
[a,b,c,d,v,w,x,y]\to (a,b,c)
$$
over the field $\F_2(a,c)$ is birational on $\ V(J3)$. The image
of $\V (J3)$ in $\F_2(a,c)[b]$ is defined by the polynomial $h$.

This implies that $J3\overline{\F}_2(a,c)[w,y,d,x,v,b]$ is a prime
ideal if $h$ is absolutely irreducible. In particular, we obtain
that $J$ is absolutely irreducible if $h$ is absolutely
irreducible.

To prove that $h$ is absolutely irreducible, we proceed as
follows:

First we show that the radical of the ideal of the coefficients of
$h$ in $\overline{\F}_2[a,c]$ with respect to $b$ is $\langle
a,c\rangle \cap \langle a + 1, c + 1\rangle$. We do this using the
factorising Gr\"obner basis algorithm.

\begin{verbatim}
ideal JF=coeffs(h,b); facstd(JF);
[1]:
   _[1]=c
   _[2]=a
[2]:
   _[1]=c+1
   _[2]=a+1
\end{verbatim}

This implies that $h$ cannot have a nontrivial factor in
$\F_2[a,c]$. Then we consider $\tilde{h}(b,c) = h(1,b,c)$.

\begin{verbatim}
subst(h,a,1);

(c+1)^14*b12+(c+1)^14*b10+(c+1)^11*(c6+c5+c4+c+1)*b8+(c+1)^11
*(c6+c4+c2+c+1)*b6+(c+1)^8*(c9+c7+c5+c4+c3+c2+1)*b4+(c+1)^10
*b2+(c+1)^10*c2;
\end{verbatim}

It is sufficient to show that $f(x,c) = \tilde{h}
\left(\tfrac{x}{c+1}, c\right)/(c+1)^2$ is absolutely irreducible.
To simplify the situation, we make the transformation $c \mapsto c
+ 1$.

Let $a_4 = c^6 + c^5 + c^4 + c^2 + 1$ and $a_2 = c^9 + c^8 + c^7 +
c^6 + c^4 + c^2 +1$.

\begin{lemma}\label{lemma3.1}
  The polynomial
\[
f = x^{12} + c^2x^{10} + c a_4 x^8 + c^3(c^6 + c+1) x^6 + c^2a_2
x^4 + c^6x^2 + c^8(c+1)^2
\]
is irreducible in $\overline{\F}_2[x,c]$.
\end{lemma}

\begin{proof}
We check that $f(x,c^2)$ is the square of some polynomial $g$, that is, $g$ is
 defined by $g^2(x,c) = f(x,c^2)$. It suffices to prove
that $g$ is irreducible: if $f = f_1 f_2$ is a non--trivial
decomposition, then $g^2=f_1(x,c^2) f_2(x,c^2)$. If $g$ is
irreducible, we obtain $g = f_1(x,c^2) = f_2(x,c^2)$. This implies
$f = f_1^2$, which is obviously not true.

\noindent\textit{First step:\/} $g$ has no linear factor in $x$.

A linear factor of $g$ has to be of the form $x - x_0 c^i(c+1)^j$ for some
$x_0 \in \overline{\F}_2$ and $i \le 8$, $j \le 2$.  Now it is easy to see,
using divisibility by $c$, that $g\bigl(x_0 c^i(c+1)^j,c\bigr) \not= 0$ for $i
= 0, 1, 2, 4, \dots, 8$.  In the case $i = 3$, $g\bigl(x_0 c^3(c+1)^j,c)\bigr)
\not= 0$ because $\bigl(x_0 c^3(c+1)^j\bigr)^6$ has degree $18 + 6j$ with
respect to $c$, which is strictly larger than the degree of the other
summands.

\noindent\textit{Second step:\/} $g$ has no quadratic factor in $x$.

Assume that $g = (x^4 + \alpha x^3 + \beta x^2 + \gamma x +
\delta) (x^2 + \varepsilon x + \mu)$ for $\alpha, \dots, \mu \in
\overline{\F}_2[c]$.   Then we obtain

\renewcommand{\arraystretch}{1.25}
$\begin{array}{lrcl}
(1) & \mu\delta &                              = & c^8(c+1)^2\\
(2) & \mu\gamma + \varepsilon \delta &         = & c^6\\
(3) & \mu\beta + \varepsilon\gamma + \delta &  = & c^2a_2\\
(4) & \mu \alpha + \varepsilon \beta + \gamma& = & c^3(c^6 + c+1)\\
(5) & \mu + \varepsilon \alpha + \beta &       = & c a_4\\
(6) & \varepsilon + \alpha & = & c^2.
\end{array}$

Now $g(x,0) = x^6$ implies that $c|\alpha, \beta, \gamma, \delta, \varepsilon,
\mu$.  Therefore, they all have degree $\le 10$.  Equation (2) implies that
$(c+1)|\mu$ and $(c+1)|\delta$ are not possible.  (3) and (4) imply that
$c^2|\delta$ and $c^2|\gamma$ and, therefore $\deg(\mu) \le 8$.  (4), (5) and
(6) imply that $\deg(\varepsilon) \le 4$ and $\deg (\alpha) \le 4$.

If $\deg(\mu) = 8$, then $\deg(\mu \gamma) \ge 10$ and (2) implies that
$\deg(\varepsilon \delta) \ge 10$.  This implies that $\deg(\varepsilon) \ge
8$, which is not possible, as we already saw.

If $\deg(\varepsilon) = \deg(\alpha) = 4$, then (5) implies
$\deg(\beta) = 8$.   This implies $\deg(\varepsilon\beta) = 12$
and, therefore, by (4), $\deg(\mu\alpha) = 12$.   This contradicts
$\deg(\mu) \le 7$ and $\deg(\alpha) = 4$.   Thus, we have
$\deg(\varepsilon) \le 3$ and $\deg(\alpha) \le 3$.   This
implies, using (5), $\deg(\beta) \le 7$.

If $\deg(\mu) = 6$, then $\deg(\delta) = 4$ implies $\deg(\mu\gamma) \ge 8$
and $\deg(\delta \varepsilon) \le 7$, contradicting (2).

We obtain $\deg(\mu) \le 5$ and, using (5), $\deg(\beta) = 7$.
If $\deg(\mu) \le 3$, then $\deg(\mu\beta) \le 10$ and (3) implies
$\deg(\varepsilon \gamma) = 11$.   We shall see that this is not
possible.

If $\deg(\varepsilon) = 3$, then $\deg(\varepsilon\beta) = 10$ and
(4) implies that $\deg(\gamma) = 10$.   This contradicts $\deg
(\varepsilon\gamma) = 11$.

If $\deg(\varepsilon)=2$ then $\deg(\varepsilon \beta) = 9$ and
$\deg(\gamma) = 9$.   This contradicts (2).

If $\deg(\varepsilon) = 1$, then $\deg(\varepsilon\beta) = 8$ and
$\deg(\gamma) = 10$.   This contradicts (4).

Finally, we obtain $4 \le \deg(\mu) \le 5$.   This implies
$c^2|\mu$ and $c^3|\delta$ and, consequently, $c^3|\mu\beta$.
But we know already that $c^2|\gamma$ and, therefore,
$c^3|\varepsilon\gamma$ and obtain a contradiction to (3).

\noindent\textit{Third step:\/} $g$ has no cubic factor in $x$.

Let $g = (x^3 + \alpha x^2 + \beta x + \gamma) (x^3 + \delta x^2 +
\varepsilon x + \mu)$, then we obtain

\renewcommand{\arraystretch}{1.25}
$\begin{array}{lrcl}
(1) & \gamma\mu                                   & = & c^8(c+1)^2\\
(2) & \gamma\varepsilon + \mu\beta                & = & c^6\\
(3) & \mu\alpha + \varepsilon\beta + \gamma\delta & = & c^2a_2\\
(4) & \mu + \alpha\varepsilon+\beta\delta+\gamma  & = & c^3(c^6+c+1)\\
(5) & \varepsilon + \alpha\delta +\beta           & = & ca_4\\
(6) & \alpha + \delta                             & = & c^2.
\end{array}$

Now, $g(x,0) = x^6$ implies $c |
\alpha,\beta,\gamma,\delta,\varepsilon,\mu$.

As in the previous case, $(c+1)|\gamma,\mu$ is not possible.   If
$c^4 \nmid \gamma$, then $c^5|\mu$ and, by (2), $c^6 |
\varepsilon\gamma$, which implies $c^3|\varepsilon$.   This
contradicts (3) and (4), because (3) implies $c^3 \nmid
\gamma\delta$ and, therefore, $c^2\nmid \gamma$.

We obtain that $c^4|\gamma$ and, by symmetry, $c^4|\mu$.  We may
assume that $\gamma = \gamma_0 c^4(c+1)^2$ and $\mu = \mu_0c^4$
for suitable $\gamma_0, \mu_0 \in \overline{\F}_2$.  This implies
$\deg(\delta) \le 5$, $\deg(\varepsilon) \le 5$ and $\deg(\alpha)
\le 7$, $\deg(\beta) \le7$ by using (3), since $a_2$ is of degree 9.

If $\deg(\alpha) \ge 4$, then $\deg(\delta) \ge 4$ by (6).   This
implies $\deg(\alpha \delta) \ge 8$, which contradicts (5).   We
obtain that $\deg(\alpha) \le 3$, $\deg(\delta) \le 3$.   This
implies, using (5), that $\deg(\beta) = 7$.   Now (4) implies that
$\deg(\delta) = 2$ and we obtain, using (3), that
$\deg(\varepsilon) = 4$. This is a contradiction to (2) and
finishes the third step.

Altogether, we proved now that $V'=\V (J)$ is absolutely
irreducible.
\end{proof}

Next we compute the singular locus of $\V (J)$, using
\textsc{Singular} (with a special procedure).

\begin{lemma} \label{lem:sing}
The singular locus of $\V (J)$ is the union of the following six
smooth curves defined by the ideals $S1, \dots, S6$.

\begin{verbatim}
S1[1]=y; S1[2]=x; S1[3]=v2+vw+w2+1; S1[4]=d+1; S1[5]=c+1;
S1[6]=b+w+1; S1[7]=a+1;

S2[1]=y+1; S2[2]=x+1; S2[3]=v+w+1; S2[4]=d; S2[5]=c;
S2[6]=b2+w2+w+1; S2[7]=a;

S3[1]=y+1; S3[2]=x+1; S3[3]=v+1; S3[4]=d; S3[5]=c; S3[6]=b2+w+1;
S3[7]=a2+ab+w;

S4[1]=y; S4[2]=x; S4[3]=v; S4[4]=d+1; S4[5]=c+1; S4[6]=b2+b+w+1;
S4[7]=a+b+1;

S5[1]=x2+y; S5[2]=wy+x; S5[3]=wx+1; S5[4]=v; S5[5]=d2+xy+x;
S5[6]=c; S5[7]=by+b+dw+d; S5[8]=bx+b+dw; S5[9]=bw+b+dw2;
S5[10]=bd+x+1; S5[11]=b2+w; S5[12]=a+dw;

S6[1]=x; S6[2]=w3y+w2+1; S6[3]=v+w2y; S6[4]=d+wy+1; S6[5]=c+w2y+w;
S6[6]=b+w; S6[7]=a;
\end{verbatim}
\end{lemma}

\begin{corol} \label{cor:S}
  The singular locus of $V'$ is contained in the set $S =V'\cap \V (xc)$.  The
  variety $U=V'\smallsetminus S$ is a smooth irreducible affine surface
  invariant under the morphism $\alpha$. For any odd $n$, $\alp^n$ has no
  fixed points in $S$.
\end{corol}

\begin{proof} The first two assertions are checked directly, looking at
the equations S1--S6 and the equation for the action of $\alp$. To
prove the third, assume that $p=(a,b,\dots, y)$ is a fixed point
of $\alp ^n$ lying on
$S$. Let $x=0$. Then, since $p$ is $\alp ^n$ invariant, we have
$c=0$. Since $p\in V'$, equation J3[1] gives $a^8v^6+a^8v^4=0$.
(The variety defined by the ideal J3 contains $V'$ as a component,
so $p$ must satisfy all the equations of J3.) Hence we have either
$a=0$, or $v=0$, or $v=1$.

In any of the two first cases we have $a=v=c=x=0$, and equation J3[3]
gives $d=0$. Since $p$ is an invariant point, we get $y=0$. Furthermore,
equation J[4] gives $w=0$. Hence $b=0$, contradiction.

If $v=1$, then $a=1$ which, taking into account $a=c=0$, contradicts J[7].
\end{proof}

\subsection{Trace formula}\label{subsec:trace}

\def\tr{\textrm{tr}}

Throughout this subsection $k$ denotes a (fixed) algebraic closure
of $\F _2$. All varieties under consideration, even those defined
over $\F _2$, are viewed as $k$-varieties.

Let $V'$ be the variety defined by equations $J[1], \dots, J[10]$
(see Subsection \ref{subsec:geom}). We have seen that this is an
irreducible affine surface.  Computations in Subsection
\ref{subsec:geom} show that the singular locus of $V'$ is
contained in the set $S =V'\cap \V (xc)$.  By Corollary
\ref{cor:S} the variety $U=V'\smallsetminus S$ is a smooth
irreducible affine surface invariant under the morphism $\alpha$
acting in $\mathbb A^8$ as
\begin{equation}
\alp (a,b,c,d,v,w,x,y)=(v,w,x,y,a^2,b^2,c^2,d^2) \label{eq:action}
\end{equation}
(see Subsection \ref{subsec:universal}).

Our goal is to prove that for $n$ odd and large enough, the set $U$ has an
$\alp^n$-invariant point. In this subsection we prove an estimate
of Lang--Weil type:

\begin{theorem} \label{th:LW}
With the above notation, let $\#\Fix (U,n)$ be the number of fixed
points of $\alp ^n$ $($counted with their multiplicities$)$. Then
for any odd $n>1$ the following inequality holds:
\begin{equation}
|\#\Fix (U,n)-2^n|\leq b^12^{3n/4}+b^22^{n/2}, \label{eq:LW}
\end{equation}
where $b^i=\dim H^i\et (U,\Ql )$ are $\ell$--adic Betti numbers
$(\ell\neq 2)$.
\end{theorem}

The strategy of proof is as follows. The operator $\alp$ and all
its powers act on the \'etale $\ell$-adic cohomology groups
$H^i\comp (U,\Ql )$ of $U$ (with compact support). We are going to
apply Deligne's conjecture (proved by T.~Zink for surfaces
\cite{Zi}, by Pink \cite{Pi} in arbitrary dimension (modulo
resolution of singularities), and by Fujiwara \cite{Fu} in the
general case) saying that the
Lefschetz(--Weil--Grothendieck--Verdier) trace formula is valid
for any operator on $U$ composed with sufficiently large power of
the Frobenius (in our case this means sufficiently large odd power
of $\alp$). We shall show that in our case the trace formula is
already valid after twisting with the first power of the
Frobenius. This fact is a consequence of the above mentioned
results on Deligne's conjecture together with the following
crucial observation: roughly speaking, if we consider the closure
$\Ub$ of $U$ in $\mathbb P^8$, $\alp$ (as well as any of its odd
powers) has no fixed points at the boundary (i.e.\ on
$\Ub\smallsetminus U$). As soon as the trace formula is
established, the proof can be finished by applying Deligne's
estimates of the eigenvalues of the Frobenius.

Let us make all this more precise.

Denote by $\Gam$ (the transpose of) the graph of $\alp$ acting on
$\A ^8$ by formulas (\ref{eq:action}), i.e.\ $\Gam =\{(\alp (M),M)
: M\in \A ^8\}$, and let $\Gam _U=\Gam\cap (U\times U)$.

Consider the natural embedding $\mathbb A^8\subset\mathbb P^8$, and denote by
$\overline\Gam$ (resp. $\overline\Gam _U$) $\subset\mathbb P^8\times\mathbb
P^8$ the closure of $\Gam$ (resp.  $\Gam _U$) with respect to this embedding.
Let $H_0=(\mathbb P^8\times\mathbb P^8)\smallsetminus(\mathbb A^8\times\mathbb
A^8)$, $H_1=(V^\prime\times V^\prime)\smallsetminus (U\times U)$, $H=H_0\cup
H_1$. Let $\Del$ denote the diagonal of $\mathbb A^8\times\mathbb A^8$,
$\overline\Del$ the diagonal of $\mathbb P^8\times\mathbb P^8$, $\Del
_U=\Del\cap\Gam _U$, and $\overline\Del _U=\overline\Del\cap\overline\Gam _U$.
If $n$ is a positive integer, denote the corresponding objects related to
$\alp ^n$ by $\Gam^{(n)}$, $\overline\Gam^{(n)}$, $\Gam _U^{(n)}$,
$\overline\Gam _U^{(n)}$, $\Del _U^{(n)}$, $\overline\Del _U^{(n)}$.

\begin{lemma} \label{lem:inf}
If $n$ is odd, $\overline\Del _U^{(n)}=\Del _U^{(n)}$.
\end{lemma}

\begin{proof}
We have
$$
\overline\Del _U^{(n)}\smallsetminus\Del
_U^{(n)}=\overline\Gam_U^{(n)}\cap\overline\Del\cap H.
$$
We wish to prove that this set is empty. Since
$$
\overline\Gam_U^{(n)}\cap\overline\Del\cap
H\subseteq\overline\Gam^{(n)}\cap (\overline{U\times U})
\cap\overline\Del\cap H,
$$
it is enough to prove that
$$
\overline\Gam^{(n)}\cap\overline\Del\cap H=\emptyset .
$$

First note that
$$
\overline\Gam^{(n)}\cap\overline\Del\cap
H_1=\Gam^{(n)}\cap\Del\cap H_1=\emptyset
$$
(the first equality is obvious since $H_1$ is contained in
$\mathbb A^8\times\mathbb A^8$, and the second one immediately
follows from Corollary \ref{cor:S}). Hence we only have to prove
that $\overline\Gam^{(n)}\cap\overline\Del\cap H_0=\emptyset . $

Let $(a,b,c,d,v,w,x,y),(a',b',\dots ,y')$ be the coordinates in
$\mathbb A^8\times\mathbb A^8$, and let $(a:b:\dots
:t),(a':b':\dots :t')$ be the homogeneous coordinates in $\mathbb
P^8\times\mathbb P^8$. Suppose that
$$
M=((a:b:\dots :t),(a':b':\dots
:t'))\in\overline\Gam^{(n)}\cap\overline\Del\cap H_0.
$$
If $n=2m+1$, denote $s=2^m$. With this notation, since $M\in
\overline\Gam^{(n)}$, formulas (\ref{eq:action}) imply that
$$
\begin{array}{llll}
a'=v^st^s, & b'=w^st^s, & c'=x^st^s, & d'=y^st^s, \\
v'=a^{2s}, & w'=b^{2s}, & x'=c^{2s}, & y'=d^{2s}, \\
t'=t^{2s}. & & &
\end{array}
$$
On the other hand, since $M\in H_0$, we have $t=t'=0$, and hence
$a'=b'=c'=d'=0$. Furthermore, since $M\in\overline\Del$, we have
$a'=\lam a$, $b'=\lam b$, $c'=\lam c$, $d'=\lam d$ for some
$\lam\in k$, and hence $a=b=c=d=0$. This implies $v'=w'=x'=y'=0$,
contradiction.
\end{proof}

The next goal is to show that the Lefschetz trace formula holds
for all odd $n$th powers of $\alp$ ($n>1$). We shall do it using
the above mentioned results on Deligne's conjecture. First we
briefly recall the general approach (\cite{SGA5}, \cite{Zi},
\cite{Pi}, \cite{Fu}); we mainly use the notation of \cite{Pi} and
refer the reader to that paper for more details.

\noindent(i) {\it Global term}. We can (and shall) view our operator $\alp$
as a particular case of the correspondence $a$:
$$
U\stackrel{a_1}\longleftarrow \Gam _U\stackrel{a_2}\longrightarrow
U
$$
(here $a_1$ and $a_2$ stand for the first and second projections,
respectively). We regard an odd power $\alp^{2m+1}$ as a
``twisted'' correspondence $b=\Fr^m\circ a$ with $b_1=\Fr ^m\circ
a_1$, $b_2=a_2$.

Let $\Lambda$ denote a finite field extension of $\Q_\ell$, $L$ a
constructible $\Lambda$-sheaf (in our situation it suffices to
consider the constant sheaf $L=\Ql$). Then a cohomological
correspondence $u$ on $L$ with support in $b$ is a morphism
$u\colon b_1^*L\to b_2^!L$, where ${}^*$ stands for the inverse
image functor, and ${}^!$ for the extraordinary inverse image
functor (cf. \cite[Section 1]{Pi} and references therein); in our
situation $b_2=\id$ and hence $b_2^!L=L$. Since $b_1$ is a proper
morphism, $u$ induces an endomorphism $u_!\colon H^{\bullet }\comp
(U,L)\to H^{\bullet }\comp (U,L)$ which possesses a well-defined
trace $\tr (u_!)\in\Lambda$; this is the global term in the
desired trace formula. In down-to-earth terms, in our situation we
have
\begin{equation}
\tr (u_!)=\sum_{i=0}^4(-1)^i\tr (\alp^n| H^i\comp (U,\Ql ).
\label{eq:global}
\end{equation}

\noindent(ii) {\it Compactification}. Furthermore, since $b_1$ is
proper, our correspondence $b$ can be extended to a
compactification $\bar b$
$$
\CD
U @<b_1<< \Gam _U @>b_2>> U \\
@VjVV @VVV @VVV \\
\bar U @<\bar b_1<< \overline\Gam _U @>\bar b_2>> \bar U
\endCD
$$
where the vertical arrows are open embeddings and the bottom line
is proper. This gives rise to a cohomological correspondence $\bar
u_!$ on the sheaf $j_!L$ with support in $\bar b$; here ${}_!$
stands for the direct image functor with compact support
(extension by 0), cf. \cite[2.3]{Pi}.

The global term does not change after compactification:
\begin{equation}
\tr (\bar u_!)=\tr (u_!). \label{eq:comp}
\end{equation}
(see \cite[Lemma 2.3.1]{Pi}).

For a compactified correspondence the Lefschetz--Verdier trace
formula is known (cf. \cite[2.2.1]{Pi}):
\begin{equation}
\tr (\bar u_!)=\sum_DLT_D(\bar u) \label{eq:Verdier}
\end{equation}
where $D$ runs over all the connected components of $\Fix (\bar
b)$, and the local terms $LT_D(\bar u)$ are defined as in
\cite[2.1]{Pi}. In our case $\Fix (\bar b)$ consists of isolated
points (since this is true for the Frobenius), and all these points
are contained in $U$ (because of Lemma \ref{lem:inf} there are no
fixed points at the boundary, neither on the singular locus, nor
at infinity).

\noindent(iii) {\it Local terms}. Suppose that $b_2$ is quasifinite and $y$
is a point not at infinity. Let $x=b_2(y)$, then
$$
d(y)=[k(y)/k(x)]_i\cdot\textrm{length}\, O_{\Gam
_U,y}/b_2^*(m_{U,x}O_{U,x}),
$$
where $[k(y)/k(x)]_i$ denotes the inseparable degree of the
residue field extension. Clearly, in our case $b_2=\id$ implies
$d(y)=1$.

By \cite[Th.\ 5.2.1]{Fu}, for an isolated fixed point $y$ at finite
distance we have
\begin{equation}
LT_y(u)=\tr_y(u) \label{eq:fd}
\end{equation}
provided $2^m>d(y)$. In our setting,
\begin{equation}
\tr_y(u) \textrm{ equals the multiplicity of } y \label{eq:naive}
\end{equation}
(cf. \cite[p.~338]{Zi}, \cite[8.3.1]{Pi}).

\noindent(iv) Summing up, (i) -- (iii) (or, more precisely, formulas
(\ref{eq:global}), (\ref{eq:comp}), (\ref{eq:Verdier}),
(\ref{eq:fd}), (\ref{eq:naive}), together with Lemma \ref{lem:inf})
imply

\begin{proposition} \label{prop:trace}
If $n>1$ is an odd integer, then
\begin{equation}
\#\Fix (U,n)=\sum_{i=0}^4(-1)^i\tr \bigl(\alp^n \mid H^i\comp (U,\Ql )\bigr)\,.
\label{eq:trace}
\end{equation}
\end{proposition}

We are now ready to prove Theorem \ref{th:LW}. Since $U$ is
non--singular, the ordinary and compact Betti numbers of $U$ are
related by Poincar\'e duality \cite[p.~6]{Ka2}, and we have
$b^i\comp =b^{4-i}$. Since $U$ is affine, $b^i=0$ for $i>2$
\cite[loc.~cit.]{Ka2}. Since $U$ is geometrically integral,
$b^0=1$ and $\Fr$ acts on the one-dimensional vector space
$H^0(U,\Ql )$ as multiplication by $4$ \cite[loc.~cit.]{Ka2}.
Hence $\alp$ acts on the same space as multiplication by $2$.
(Indeed, if it were multiplication by $(-2)$, for a sufficiently
big power of $\alp$ the right-hand side of (\ref{eq:trace}) would
be negative.) Hence $\alp^n$ acts as multiplication by $2^n$. Thus
$\tr (\alp^n|H^4\comp (U,\Ql ))=2^n$.

On the other hand, according to Deligne \cite[Th.~1]{De} for every
eigenvalue $\alp _{ij}$ of $\Fr$ acting on $H^i\comp (U,\Ql )$ we
have $|\alp _{ij}|\leq 2^{i/2}$. This yields similar inequalities
for the eigenvalues $\bet_{ij}$ of $\alp$: $|\bet _{ij}|\leq
2^{i/4}$ and the eigenvalues $\bet_{ij,n}$ of $\alp^n$:
$|\bet_{ij,n}|\leq 2^{ni/4}$. We thus get
$$
|\tr \bigl(\alp^n|H^3\comp (U,\Ql )\bigr)|\leq b^12^{3n/4}\,,
$$
$$
|\tr \bigl(\alp^n|H^2\comp (U,\Ql )\bigr)|\leq b^22^{n/2}\,.
$$

This proves the theorem. \qed

\begin{remark}
  Probably one can get another proof of Proposition \ref{prop:trace} (and
  hence Theorem \ref{th:LW}) using an approach of \cite{DL}. In that paper the
  Lefschetz trace formula is established for any endomorphism of finite order.
  A remark in Section 11 of the above cited paper (see also \cite[Sommes
  trig., 8.2, p.~231]{SGA41/2}) says that the results of the paper can be
  extended to the case of an endomorphism $\alp$ with the property $\alp
  ^2=\Fr$.
\end{remark}

\subsection{Estimates of Betti numbers} \label{subsec: Betti}

\newcommand{\ef}{\end{equation}}
\chardef\bslash=`\\ 
\newcommand{\ntt}{\normalfont\ttfamily}
\newcommand{\cn}[1]{{\protect\ntt\bslash#1}}
\newcommand{\pkg}[1]{{\protect\ntt#1}}
\newcommand{\fn}[1]{{\protect\ntt#1}}
\newcommand{\env}[1]{{\protect\ntt#1}}
\hfuzz1pc 

\newcommand{\cA}{\mathcal{A}}
\newcommand{\B}{\mathcal{B}}
\newcommand{\st}{\sigma}
\renewcommand{\k}{\varkappa}

\newcommand{\X}{\mathcal{X}}
\newcommand{\wt}{\widetilde}
\newcommand{\wh}{\widehat}
\newcommand{\mk}{\medskip}
 \renewcommand{\sectionmark}[1]{}
\renewcommand{\Im}{\operatorname{Im}}
\renewcommand{\Re}{\operatorname{Re}}
\newcommand{\la}{\langle}
\newcommand{\ra}{\rangle}

 \renewcommand{\th}{\theta}
\newcommand{\ve}{\varepsilon}
\newcommand{\1}{^{-1}}
\newcommand{\iy}{\infty}
\newcommand{\iintl}{\iint\limits}
 \newcommand{\cupl}{\operatornamewithlimits{\bigcup}\limits}
\newcommand{\suml}{\sum\limits}
\newcommand{\ord}{\operatorname{ord}}
\newcommand{\bk}{\bigskip}
\newcommand{\fc}{\frac}
\newcommand{\g}{\gamma}
 \newcommand{\dl}{\delta}
\newcommand{\Dl}{\Delta}
\newcommand{\lm}{\lambda}
\newcommand{\Lm}{\Lambda}
\newcommand{\ov}{\overline}
\newcommand{\vp}{\varphi}
\newcommand{\BC}{\field{C}}
\newcommand{\C}{\field{C}}
\newcommand{\BM}{\field{M}}
\newcommand{\BR}{\field{R}}
\newcommand{\BZ}{\field{Z}}
\renewcommand{\Im}{\operatorname{Im}}
\newcommand{\ep}{\endproclaim}

\newcommand{\doe}{\overset{\text{def}}{=}}
 \newcommand{\supp} {\operatorname{supp}}

\newcommand{\z}{\zeta}
\renewcommand{\a}{\alpha}

As in the previous subsection, we assume that the ground field is
$k=\bar\F _2$.

Recall that we consider the variety $V'$ defined by equations
$J$[1--10] (see Subsection \ref{subsec:geom}) whose singular locus
is contained in the set $S =V'\cap \V (xc)$. As before, we denote
$U=V'\smallsetminus S$; it is a smooth irreducible affine variety
invariant under the morphism $\alpha$. Our aim is to estimate
$b^1(U)$ and $b^2(U)$.

First we deal with $b^1(U)$. We want to use the Lefschetz Theorem
on hyperplane sections. For technical reasons we want to use
hyperplanes of special type, namely those defined by equations
$\alpha a+\beta c+\gamma=0$. These hyperplane sections are not
general, and in order to apply the Lefschetz Theorem, we have to
provide a quasifinite map of the surface $V'$ onto $\mathbb A^2$
with coordinates $a,c$.

The next step is to estimate the Euler characteristic of $U$. To
do this, we represent $U$ as the union of an open subset $U'$ and
a finite number of curves. We estimate the Euler characteristics
of these curves and of $U'$ separately, using the fact that $U'$
is a double cover of a simpler variety. Having in hand bounds for
$b^1(U)$  and $\chi(U),$ we estimate $b^2(U)$.

\begin{prop}\label{prop1}
A regular map $\pi\colon U\to \mathbb A^2$ defined as
$\pi(a,b,c,d,v,w,x,y) =(a,c)$ is quasifinite.
\end{prop}

\begin{proof}
Consider the variety $\widetilde W$ defined in $\mathbb A^8$ by
equations $J3[1-6]$ (see Subsection \ref{subsec:geom}).

We have $\widetilde W\supset V'$ and $\widetilde W\smallsetminus
V'\subset \V (f) \subseteq \A^8$,  where
$$f(a,c)=c(ac+1)a(a+c)(c+1)(a^3+a^2c^2c^3+a^2c+ac^4+ac^2+c)$$ (see
Subsection \ref{subsec:geom}). Thus, if $f(a,c)\ne0,$ the equation
$J3[1]$ provides at most six different possible values for $v$.
The equation $J3[2]$ implies that for each of these six values
only one value of $x$ is possible. The equation $J3[3]$ gives at
most two values for $d,$ and all the proceeding equations provide
one value for $b,$ $y$ and $w$. Hence, for any point $(a,c)\in
\mathbb A^2$, the preimage $\pi^{-1}(a,c)$ is finite if
$f(a,c)\ne0$.

Let now $A=\V (f) \subset \A^2$. Then $\pi^{-1}(A)\cap U=\cup
A_i,$\ $i=1,\dots ,6$ which may be described as follows.

\begin{enumerate}
\item[1.] $A_1= U\cap \V (c-1)$.
\end{enumerate}

According to calculations, $A_1=A_1^1\cup A_1^2,$ where
$A_1^1=U\cap \V (c-1)\smallsetminus \V
\bigl(a(a+1)(a^2+a+1)\bigr)$ and $A_1^2=U\cap \V\bigl(c-1,
a(a+1)(a^2+a+1)\bigr)$.

The set $A_1^1$ is defined by the ideal L.
\begin{verbatim}
L[1]=c+1;
L[2]=(a5+a4+a3+a2)*v4+(a5+a)*v2+(a4+a2+1);
L[3]=x+(a3+a2)*v3+(a3+a2+a)*v;
L[4]=(a+1)*d2+(a4+a2)*dv3+(a4+a)*dv+(a8+a6+a5+a4+a3+a2+1)*v2+(a8+a5+a+1);
L[5]=(a4+a2+1)*b+(a5+a4+a2+a)*dv2+(a4+a)*d+(a6+a4)*v3+(a6+a2+a+1)*v;
L[6]=(a2+a+1)*y+(a2+a+1)*d+(a7+a6+a5+a2)*v3+(a7+a6+a4+a3+a2+a)*v;
L[7]=(a4+a2+1)*w+(a6+a5+a3+a2)*dv2+(a5+a2)*d+(a7+a4+a2+a)*v3+(a7+a6+a5+a)*v;
\end{verbatim}
These equations show that for a fixed value of $a,$ if
$a(a+1)(a^2+a+1)\ne 0,$ there are at most six points in $U\cap
\pi^{-1}(a,1)$. The set $A_1^2$ is defined by the ideal L1.
\begin{verbatim}

L1[1]=a2+a+1
L1[2]=c+1
L1[3]=v+a+1
L1[4]=x+a
L1[5]=d2+da+1
L1[6]=b+da+d
L1[7]=y+d+a
L1[8]=w+d+1
\end{verbatim}
It follows that $\pi^{-1}(1,1)=\emptyset$;
$\pi^{-1}(0,1)=\emptyset;$ \ $\pi^{-1}(a_0,1),$ where $a_0$  is  a
root of $a^2+a+1$, consists of two points.

\begin{enumerate}
\item[2.] $A_2=U\cap \V (c)=\emptyset$.

\item[3.] $A_3=U\cap \V (a)=\emptyset$.

\item[4.] $A_4=U\cap \V (a+c)=\emptyset$.

\item[5.] $A_5=U\cap \V (ac+1)=A_5^1\cup A_5^2$.

\item[5.1]  $A_5^1=U\cap \V (ac+1) \cap \D\bigl((a^2+a+1)(a-1)a\bigr)$. This set
is defined by the ideal D.
\end{enumerate}
\begin{verbatim}
D[1]=(a)*c+1;
D[2]=(a3+a2)*v2+(a3+a2+a)*v+1;
D[3]=x+(a4)*v3+(a4)*v;
D[4]=(a6+a2)*d2+(a9+a5)*dv3+(a9+a5)*dv+(a10+a6+a4+a2+1)*v2+(a10+a2+1);
D[5]=(a3+a)*b+(a5+a3)*dv2+(a)*d+(a4+a2)*v3+(a4+a2+1)*v;
D[6]=(a)*y+d+(a7+a5)*v3+(a7+a5+a3)*v;
D[7]=(a2+1)*w+(a5+a3)*dv2+(a)*d+(a4+a2)*v3+v;
\end{verbatim}
which show that for any point $a\ne0,1,$ or $a_0$ (a root of $a^2+a+1),$ the
set $\pi^{-1}(a,\frac{1}{a})$ contains at most four points.

\begin{enumerate}
\item[5.2]  $A_5^2=U\cap \V\bigl((a-1)(a-a_0)a\bigr)$.
\end{enumerate}

This set consists of four points defined by the ideal D1.
\begin{verbatim}
D1[1]=a2+a+1
D1[2]=c+a+1
D1[3]=v+a
D1[4]=x+1
D1[5]=d2+da+d+1
D1[6]=b+d+a
D1[7]=y+da+d+a+1
D1[8]=w+da+a
\end{verbatim}

\begin{enumerate}
\item[6.] $A_6=U\cap \V (h_1)$, where
  $h_1(a,c)=a^3+a^2c^3+a^2c+ac^4+ac^2+c$.

\noindent $A_6=A_6^1\cup A_6^2\cup A_6^3$, where

\noindent $A_6^1=U\cap \V (h_1)\cap \D\bigl(v^2+ac^3+c^2+a^2,
a(a+1) (a^2+a+1)\bigr)$;

\noindent $A_6^2=U\cap \V (h_1, v^2+ac^3+c^2+a^2)\cap
\D\bigl(a(a+1)(a^2+a+1)\bigr)$;

\noindent $A_6^3=U\cap \V\bigl(a(a+1)(a^2+a+1)\bigr)$.
\end{enumerate}

\noindent The set $A_6^1$ is defined by the ideal K1.

\begin{verbatim}
K1[1]=(a)*c4+(a2)*c3+(a)*c2+(a2+1)*c+(a3);
K1[2]=(a3+a)*v4+(a6+a4+a2)*v2c3+(a5+a)*v2c2+v2c+(a7+a3+a)*v2+(a10+a8+1)*c3
     +(a9+a5)*c2+(a4+a2+1)*c+(a11+a7+a5);
K1[3]=(a2+1)*xv2+(a3+a)*xc3+(a2+1)*xc2+(a4+a2)*x+(a)*v3c3+(a)*v3c+v3+(a5)*vc3
     +(a4+a2)*vc2+(a)*vc+(a6+a4+a2)*v;
K1[4]=(a6+a4+a2+1)*x2+(a7+a5)*xvc3+(a8+a4)*xvc2+(a7+a)*xvc+(a8+a2)*xv
     +(a3)*v2c3+(a8+a6)*v2c2+(a7+a5+a3)*v2c+(a6)*v2+(a11+a3+a)*c3+(a10+a8)*c2
     +(a5+a3+a)*c+(a12+a10+a6+a4+1);
K1[5]=(a10+a6+a4+1)*d2+(a10+a6+a4+1)*dxc+(a13+a11+a9+a7+a3+a)*xvc3
     +(a14+a12+a4+a2)*xvc2+(a13+a9+a7+a5)*xvc+(a6+a4)*xv+(a13+a11+a9+a7+a3)
     *v2c3+(a14+a12+a8+a6+a2+1)*v2c2+(a13+a11+a5+a3+a)*v2c+(a8+a4+a2)*v2
     +(a13+a9+a5)*c3+(a12+a6)*c2+(a5+a3+a)*c+(a14+a12+a8+a4+a2);
K1[6]=(a12+a10+a8+a6+a4+a2)*b+(a3+a)*dxvc3+(a8+a6)*dxvc2+(a5+a)*dxvc
     +(a10+a4+a2+1)*dxv+(a5+a3)*dv2c3+(a10+a8)*dv2c2+(a7+a3)*dv2c
     +(a12+a6+a4+a2)*dv2+(a7+a5+a3)*dc3+(a12+a8+a6+a2)*dc2+(a11+a9+a7)*dc
     +(a12+a10+a8)*d+(a8+a2)*xc3+(a7+a)*xc2+(a10+a8+a4+a2)*xc+(a11+a9+a5+a3)
     *x+(a13+a11+a9+a7+a5+a3)*v3c2+(a14+a12+a10+a8+a6+a4)*v3c+(a13+a11+a9+a7
     +a5+a3)*v3+(a6+a4+a2)*vc3+(a13+a11+a7+a5)*vc2+(a14+a10+a8+a6+a2)*vc
     +(a11+a9+a7)*v;
K1[7]=y+bdx+bdv3+bdv+bc3+bc+d3c+(a)*d3+dx2c3+dx2c+(a)*dxvc12+(a2)*dxvc11
     +(a2+1)*dxvc9+(a)*dxvc6+(a)*dxvc4+(a2)*dxvc3+(a)*dxvc2+(a2)*dxvc+(a)*dxv
     +dv4c+(a)*dv2c12+(a2)*dv2c11+(a2+1)*dv2c9+dv2c7+(a)*dv2c6+dv2c5+(a2+1)
     *dv2c3+dv2c+(a)*dc10+(a)*dc8+(a2)*dc7+(a2)*dc5+(a2)*dc3+(a)*dc2+(a)
     *d+x3c2+x3+x2vc2+x2v+(a)*xc21+(a2)*xc20+(a2+1)*xc18+(a)*xc17+(a2+1)
     *xc16+(a2)*xc14+xc12+(a)*xc11+(a2+1)*xc10+(a)*xc9+(a)*xc7+(a2)*xc6+(a)
     *xc5+xc4+(a)*xc3+(a)*xc+(a)*v3c19+(a2)*v3c18+(a)*v3c17+v3c16+(a2)*v3c14
     +(a)*v3c11+(a2+1)*v3c10+(a2)*v3c8+(a)*v3c7+(a2)*v3c6+v3c4+v3c2+(a)*v3c
     +(a2)*v3+(a)*vc21+(a2)*vc20+(a)*vc19+vc18+(a2)*vc16+(a)*vc15+(a2)*vc14
     +(a)*vc13+vc10+(a)*vc9+vc6+(a)*vc5+(a2)*vc4+(a2+1)*vc2+(a2+1)*v;
K1[8]=w+bdx+bdv3+bdvc2+(a)*b+(a)*d3+dx2c3+dx2c+(a)*dxvc12+(a2)*dxvc1
     1+(a)*dxvc10+dxvc9+(a)*dxvc8+(a2)*dxvc5+(a2+1)*dxvc3+(a)*dxvc2+(a2)*dxvc
     +(a)*dxv+(a)*dv2c12+(a2)*dv2c11+(a)*dv2c10+dv2c9+(a)*dv2c8+(a)*dv2c6
     +(a2)*dv2c5+dv2c3+dv2c+(a)*dc10+(a2)*dc7+dc5+(a2+1)*dc3+(a2)*dc
     +(a)*d+x3c2+(a)*xc21+(a2)*xc20+(a)*xc19+xc18+(a2)*xc16+(a)*xc15+(a2)
     *xc14+(a)*xc13+(a)*xc11+xc10+(a2)*xc6+(a2)*xc2+x+v5+(a)*v3c19+(a2)*v3c18
     +(a2+1)*v3c16+v3c14+(a)*v3c13+(a2+1)*v3c12+(a)*v3c11+v3c10+(a)*v3c9+v3c8
     +(a)*v3c7+v3c6+(a)*v3c5+(a)*v3c3+(a2)*v3c2+(a)*v3c+(a2+1)*v3+(a)*vc21
     +(a2)*vc20+(a2+1)*vc18+vc16+vc14+(a2+1)*vc10+(a2)*vc8+vc6+(a)*vc5
     +(a)*vc3+vc2+(a)*vc+v;
\end{verbatim}

This shows that each point $(a,c)$ such that $f(a,c)=0,$\
$v^2+ac^3+c^2+a^2\ne0,$\ $a(a+1)(a^2+a+1)\ne0$ has at most four
preimages in $U_1$. The set $A_6^2$ is defined by the ideal K2.

\begin{verbatim}
K2[1]=(a)*c4+(a2)*c3+(a)*c2+(a2+1)*c+(a3);
K2[2]=v2+(a)*c3+c2+(a2);
K2[3]=(a4+1)*x2+(a5)*xvc3+(a6+a4)*xvc2+(a5+a3+a)*xvc+(a6+a4+a2)*xv
     +(a3+a)*c3+(a6+a2)*c2+(a3+a)*c+(a2+1);
K2[4]=(a9+a7+a3+a)*d2+(a9+a7+a3+a)*dxc+(a12+a8+a2)*xvc3+(a13+a3)*xvc2
     +(a12+a10+a6)*xvc+(a5)*xv+(a12+a10+a8+a6+a4+a2)*c3+(a13+a11+a7+a5)*c2
     +(a6+1)*c+(a11+a5);
K2[5]=(a12+a10+a8+a6+a4+a2)*b+(a3+a)*dxvc3+(a8+a6)*dxvc2+(a5+a)*dxvc
     +(a10+a4+a2+1)*dxv+(a11+a7+a3)*dc3+(a12+a10+a8+a6+a4+a2)*dc2
     +(a11+a7+a3)*dc+(a12+a8+a4)*d+(a8+a2)*xc3+(a7+a)*xc2+(a10
     +a8+a4+a2)*xc+(a11+a9+a5+a3)*x+(a6+a4+a2)*vc3+(a9+a3)*vc2
     +(a10+a8+a6)*vc+(a11+a9+a7)*v;
K2[6]=(a4+a2+1)*y+(a4+a2+1)*dc+(a)*xc3+(a4+a2)*xc2+(a5)*xc+(a6+a4+a2+1)*x
     +(a5+a3+a)*vc+(a6+1)*v;
K2[7]=(a11+a9+a7+a5+a3+a)*w+(a3+a)*dxvc3+(a8+a6)*dxvc2+(a5+a)*dxvc+(a10
     +a4+a2+1)*dxv+(a11+a7+a3)*dc3+(a12+a10+a8+a6+a4+a2)*dc2+(a11+a7+a3)*dc
     +(a12+a8+a4)*d+(a10+a8+a6+a4+a2+1)*xc+(a11+a9+a7+a5+a3+a)*x+(a6+a4+a2)
     *vc3+(a9+a3)*vc2+(a10+a8+a6)*vc+(a13+a5+a3)*v;
\end{verbatim}

It follows that in this case the preimage of each point is finite.

The set $A_6^3$ consists of 54 points defined by the ideals W1, W2, W3, W4,
H1, H2 and H3.
\begin{verbatim}
W1[1]=c3+c+1                    W2[1]=c3+c+1
W1[2]=a+1                       W2[2]=a+1
W1[3]=v+1                       W2[3]=v+1
W1[4]=x+c2+c                    W2[4]=x+c2+c
W1[5]=d+c2+1                    W2[5]=d+c
W1[6]=b+c2+c                    W2[6]=b+c2+c+1
W1[7]=y+c+1                     W2[7]=y+c2+c
W1[8]=w+c2                      W2[8]=w+c2+1

W3[1]=c3+c+1                    W4[1]=c3+c+1
W3[2]=a+1                       W4[2]=a+1
W3[3]=v+c2                      W4[3]=v+c2
W3[4]=x+c                       W4[4]=x+1
W3[5]=d2+dc2+1                  W4[5]=d2+dc+c
W3[6]=b+d+c                     W4[6]=b+dc
W3[7]=y+dc+1                    W4[7]=y+dc+c
W3[8]=w+d+c2+1                  W4[8]=w+dc

H1[1]=a2+a+1                    H2[1]=a2+a+1
H1[2]=c3+c2a+c2+ca+a+1          H2[2]=c3+c2a+c2+ca+a+1
H1[3]=v+c2a+c2+1                H2[3]=v+c2a+c2+1
H1[4]=x+c2a+c2+a+1              H2[4]=x+ca+1
H1[5]=d2+dc2a+dca+da+c2a+c2     H2[5]=d2+dc2a+dc+c
H1[6]=b+dc2+dc+da+d+c2a         H2[6]=b+dc2a+dc2+dca+da+c+a+1
H1[7]=y+dc+a                    H2[7]=y+dc+c2+ca+a+1
H1[8]=w+dc2a+dca+d+ca+c+1       H2[8]=w+dc2+dca+dc+da+d+c2a+c2+ca+c+a

H3[1]=a2+a+1
H3[2]=c3+c2a+c2+ca+a+1
H3[3]=v+c2a+ca+c+1
H3[4]=x+c2a+c2
H3[5]=d2+dc2a+dc+da+c2a+c2+ca+c
H3[6]=b+dc2a+dc2+dca+c2a
H3[7]=y+dc+c2a+ca+c+a
H3[8]=w+dc2+dca+dc+ca+c+a
\end{verbatim}

Thus, any point in $\mathbb A^2$ has a finite (maybe, empty)
preimage. Hence $\pi$ is quasifinite.
\end{proof}

Further on we shall consider the following sets:

$V'\subset \mathbb A^8$, defined by the ideal J;

$\widetilde W\subset \mathbb A^8$, defined by the ideal J3;

$W\subset \mathbb A^4$ with coordinates $(a,c,v,x),$ defined by
the ideal $\langle$J3(1), J3(2)$\rangle$.

$L=W\cap \V (f)\subset \mathbb A^4;$

$U=V'\smallsetminus \V (xc)\subset \mathbb A^8$;

$U'=V'\smallsetminus \V (f)\subset \mathbb A^8$;

$Y= \V (J3[1]) \cap \D (f)\subset A^3$ with coordinates $(a,c,v)$.

$Z=W \smallsetminus L\subset \mathbb A^4;$

These affine sets are included in the following diagram:
\begin{alignat*}{3}
&\quad\widetilde W\supset V' && \supset U\supset &&\ U'\\
&{}^{\pi_1} \downarrow && &&\downarrow{}^{\pi_1}\\
&\quad W &&\supset && Z\\
& && &&\downarrow{}^{\pi_2}\\
& && && Y
\end{alignat*}
The inclusion $U\supset U'$ follows from computations: we have
$V'\cap \V (x)\subset \V (f)\cap V'$. The map $\pi_1$ is a double
unramified cover. This follows from the structure of equations
$J3[1], \dots, J3[6]$: all the branch points are contained in the
set $\V (f)$. The map $\pi_2$ is an isomorphism since $x$ appears
linearly in the equation $J3[2]$ and its coefficient does not
vanish in $U'$.

\begin{prop}\label{prop2}
$b^1(U)\le 675$.
\end{prop}

\begin{proof}
This estimate follows from the Weak Lefschetz Theorem proved by
N.~Katz (\cite[Cor.\ 3.4.1]{Ka1}). Indeed, we have:
\begin{itemize}
\item an algebraically closed field of characteristic $2\le \ell$.
\item $U$, a separated $k$-scheme of finite type which is a local
complete intersection, purely of dimension $2>0$.
\item $U\to \mathbb A^2,$ a quasifinite morphism (see
Proposition \ref{prop1}).
\end{itemize}

Then, for a constant $\mathbb Q_\ell$-sheaf $\mathcal F$ on $U$,
there exists a dense open set $\mathcal U\subset \A^3$ such that
for any $(\alpha,\beta,\gamma)\in \mathcal U$ the restriction map
\[
H^1(U,\mathcal F)\to H^1\bigl(U\cap\{\alpha a+\beta c+\gamma=0\}, i^\ast
\kf\bigr)
\]
is injective ($i$ denotes the embedding of the hyperplane section
into $U$).

Denote:
\begin{align*}
S_1 & = U \cap \V (\alpha a + \beta c + \gamma);\\
  \widetilde{S} & = S_1 \cap U^\prime = U^\prime \cap \V (\alpha a + \beta c +
  \gamma) \subset S_1;\\
S & = Y \cap \V (\alpha a + \beta c + \gamma) \subset Y\,.
\end{align*}

Since $U'$ is a double unramified cover of $Y,$\ $\widetilde S$ is
a double unramified cover of $S$. The curve $S$ is defined in
$\A^3$ with coordinates $(a,c,v)$ by $\V (J3[1],\alpha a+\beta
c+\gamma) \cap \D (f)$ with
\begin{verbatim}
J3[1]=(a8+a6c2+a4c4+a2c6)*v6+(a8+a7c3+a6c2+a5c3+a4c4+a3c7+a2c6
      +ac7)*v4+(a7c3+a6c2+a5c5+a5c3+a3c7+a3c5+a2c6+a2c4+ac9+c6)
      *v2+(ac9+ac5+c8+c4)=0;
\end{verbatim}

Let $\overline S$ be the projectivization of $S$ in $\mathbb P^3$.
For a general triple $(\alpha,\beta,\gamma)$ it is an irreducible
complete intersection of degree $d=14$. By \cite[Cor.~7.4]{GL}, we
have
\[
b^1(\overline S)\le (d-1)(d-2)\le 156\,.
\]
Let $B$ be the union of the plane at infinity with the closure of
the set $\V\bigl((\alp a+\beta c+\gamma)f(a,c)\bigr)$. Since $\deg
f=11,$ we have $\deg B=13$. Thus $\overline S\cap B$ contains at
most $14\cdot 13=182$ points. Hence $b^1(S)\le 156+182=338$. Since
$\widetilde S$ is a double unramified cover of $S,$ \
$b^1(\widetilde S)=2b^1(S)-1\le 675$. Since $\widetilde S\subset
S_1,$\ $b^1(S_1)\le b^1(\widetilde S)\le 675$.
\end{proof}

\begin{prop}\label{prop3}
The Euler characteristic of $L$ can be estimated as follows:
$\chi(L)\le 71430 < 2^{17}$.
\end{prop}

\begin{proof}
The set $L=W\cap \V (f)$ consists of several components. According
to computations, the list of components is as follows:
\[
\begin{array}{lcllcllcl}
F_1 & = & \V (a,c);    & \dim F_1 & = & 2, & \chi(F_1) & = &\phantom{-}1\\
F_2 & = & \V (v,c);    & \dim F_2 & = & 2, & \chi(F_2) & = &\phantom{-}1\\
F_3 & = & \V (v-1,c);  & \dim F_3 & = & 2, & \chi(F_3) & = &\phantom{-}1\\
F_4 & = & \V (a-1,c-1);& \dim F_4 & = & 2, & \chi(F_4) & = &\phantom{-}1\\
E   & = & \V (ac-1,v); & \dim E   & = & 2, & \chi(E)   & = &\phantom{-}0\\
G   & = & \V (ac-1, av^2 + c^2 + av + cv + v^2 + v),& \dim G&  =&
2,& \chi(G) &
= & -3\\
C_1 & =  & \V (x,a,c^2+cv+1),& \dim C_1 & = & 1,& \chi(C_1)& = &\phantom{-}0\\
C_2 & =  & \V (c-1,v,x),     & \dim C_2 & = & 1,& \chi(C_2)& = &\phantom{-}1\\
C_3 & =  & \V (I_3),         & \dim C_3 & = & 1,&                           \\
H_1 & =  & \V (I_1),         & \dim H_1 & = & 2,&                           \\
H_2 & =  & \V (I_2),         & \dim H_2 & = & 2,&
\end{array}
\]
where $I_3 = \langle c-1,a^2v^2x+a^2v+v^2x+av+ax+v+x,
a^4x^2+a^2vx^3+a^3v^2+a^3x^2+a^4
+a^2vx+vx^3+avx+ax^2+a^2+vx+1,av^2x^4+v5x+v^4x^2+v^2x^4+av^4+avx^3+v^4
+a^2vx+a^2x^2+vx^3+x^4+avx+vx+x^2\rangle$,

$I_1 = \langle c^3+c^2v+c^2+av+cv+v^2,
acv+ac+c^2+av+v^2+a+c+v, a^2v+a^2+ac+cv+v^2+c\rangle$, and

$I_2 = \langle ac^2v+c^3v+c^3+c^2v+av^2+cv^2+cv+a,
c^4+acv+c^2v+ac+cv+v+1, a^3v^2+a^2v^3+acv^3+c^3v+a^2v^2+acv^2+cv^2\rangle$.

Let us explain how the Euler characteristics were computed. We
have $\chi(F_i)=1, i=1,\dots ,4$ because the $F_i$'s are just
affine spaces. $E$ is isomorphic to $\mathbb A^1$ with coordinate
$a$ punctured at the point $a=0,$ so $\chi(E)=1-1=0$. The
component $G$ is the direct product of $\mathbb A^1$ with
coordinate $x$ and a curve $T$ which is a ramified covering of
$\mathbb A^1$ with coordinate $a$. For a fixed point $(a,c,v)$ in
$T$ we have $c=\frac{1}{a}$, and $v$ is defined by the quadratic
equation
\[
v^2(a^3+a)+v(a^3+a^2+a)+1=0\,.
\]
It follows that if $a\ne 0, a\ne 1,a^2+a+1\ne 0$, there are
precisely two points in $T$ with this value of $a$. There are no
points with $a=0$ and precisely one point for each value $a=1$ or
$a^2+a+1=0$. Since the Euler characteristics of $\mathbb A^1$
without 4 points is $-3$, we have $\chi(G)=2(-3)+3=-3$.

In order to estimate the Euler characteristics of $C_3,$ $H_1,$
$H_2,$ we use the following theorem of Adolphson and Sperber:

\begin{prop} \textrm{\cite[Th.~5.27]{AS}, \cite{Ka2}} \label{prop:AS}
If an affine variety $V$ is defined in $\A^N$ by $r$ polynomial
equations all of degree $\le d,$ then
\begin{equation}
|\chi(V)|\le 2^r D_{N,r}\underbrace{(1,1+d,\dots,1+d)}_{r+1},
\label{eq:AS}
\end{equation}
where $D_{N,r}(x_0,\dots,x_r)=\sum\limits_{|W|=N} X^W$ is the
homogeneous form of degree $N$ in $x_0,\dots ,x_r$ all of whose
coefficients equal $1$.
\end{prop}

According to formula (\ref{eq:AS}),

$|\chi(C_3)|\le 2^3D_{3,3}(1,8,8,8)\le 44232 < 2^{16}$

$|\chi(H_1)|\le 2^3D_{3,3}(1,4,4,4)\le 5992  <  2^{13}$

$|\chi(H_2)|\le 2^3D_{3,3}(1,6,6,6)\le 19160 <  2^{15}$.

The pairwise intersection of these components is a union of 16
lines and 10 points. The triple intersections contain 3 lines and
3 points. No four of these components intersect. Thus,
$|\chi(L)|\le 5-3+44232+5992+21224-26+6<2^{17}$.
\end{proof}

\begin{prop}\label{prop4}
$b^2(U)\le 2^{22}$.
\end{prop}

\begin{proof}
We consider two cases:

\noindent \phantom{I}I.\ $\chi(U)\le 0$. Then $1-b^1(U)+b^2(U)\le
0$ and $b^2(U)\le b^1(U) < 675$.

\noindent II.\ $\chi(U)>0$. We first find $ |(\chi(U')|$. Since
$U'$ is a double cover of $Z,$ we have $|\chi(U')|=2|\chi(Z)|$.
Since $Z=W\smallsetminus L,$ we have $\chi(Z)=\chi(W)-\chi(L)$. By
formula (\ref{eq:AS}), we get $|\chi(W)|\le 2^2D_{4,2}(1,15,15)\le
1128908$. In view of Proposition \ref{prop3}, we have
$|\chi(L)|\le 71430$. Hence $|\chi(Z)|\le |\chi(W)|+|\chi(L)|\le
1200338,$ and therefore $|\chi(U')|\le 2400676 < 2^{22}$. On the
other hand, $\chi(U)=\chi(U')+\chi(U\smallsetminus U')$. In order
to find $\chi(U),$ we have to evaluate $\chi(U\smallsetminus U')$.
Let $N=U\setminus \V (f)$. Since $N$ is the intersection of the
smooth affine surface $U$ with the hypersurface $\V (f)$, all of
its irreducible components $N_i$ are curves (i.e., $\dim N_i=1$).
This follows from \cite[Th.5, p.74]{Sh}, and is confirmed by
calculations. Since by Proposition \ref{prop1} the projection
$\pi\colon U\to \mathbb A^2$ is quasifinite, none of $N_i$ is
mapped into a point. Hence $\pi(N_i)\subset \mathbb A^2$ is a
curve. This curve does not meet the lines $\V (c)$ and $\V (a)$
because $\V (a)\cap U_1=\emptyset $. This means that the ring
$O(\pi(N_i))$ contains the nonvanishing function $ac$. If
$ac=\const$ on $\pi(N_i)$, then $\pi(N_i)$ has two punctures at
infinity. If $ac\ne\const,$ then the normalization of $\pi(N_i)$
has at least two punctures, as does any curve having a nonconstant
and nonvanishing regular function. Thus $\chi (\pi(N_i))\le 0$.
But then for any $N_i$ we have
$$\chi(N_i)\le \chi \bigl(\pi(N_i)\bigr)\le 0\,.$$
Hence
$$\chi\left(\bigcup N_i\right)=\chi(N)\le 0;$$
$$\chi(U)=\chi(U')+\chi(N)\le \chi(U')\le 2400676,$$
and, therefore,
$$b^2(U)=\chi(U)+b^1(U)\le 2401351 < 2^{22}.$$
\end{proof}

\begin{corol} \label{cor:big}
Let $n>48$, $q=2^n$. Then $V_n$ has an $\F _q$-point.
\end{corol}

\begin{proof} On plugging the estimates of Propositions
\ref{prop2} and \ref{prop4} into formula (\ref{eq:LW}), we see
that $\#\Fix (U,n)>0$ as soon as $n>48$. This proves the
corollary.
\end{proof}

\subsection{Small fields} \label{subsec:small}

\def\BF{\mathbb F}

The purpose of this section is to study the fixpoints and also
numbers of fixpoints of the operator $\alpha^n$ on the variety
$V'$ given by the equations $J[1],\ldots, J[10]$. Let $k$ denote
the algebraic closure of $\BF_2$ and $N_n$ the number of fixpoints
of $\alpha^n$ on $V'(k)$. As explained before, if $n$ is even
($n=2k$) then $N_n$ is just the number of points of $V'$ in the
field $\BF_{2^k}$. We are interested here in the numbers $N_p$ for
odd primes $p$.

We first give a table of the numbers $N_n$ for $1\le n\le 23$.
\enlargethispage{1cm}
\begin{table}[ht]
\[
\begin{array}{ccc}
\begin{array}{r|r}
n & N_n\\\hline
1 & 0 \\
2 & 8 \\
3 & 12\\
4 & 16\\
5 & 20\\
6 & 56 \\
7 & 140\\
8 & 240
\end{array} &\quad
\begin{array}{r|r}
n & N_n\\\hline
9 & 516\\
10 & 1088\\
11 & 2332 \\
12 & 3904\\
13 & 8372\\
14 & 16416\\
15 & 32012\\
16 & 65360
\end{array} &\quad
\begin{array}{r|r}
n & N_n\\\hline
17 & 130084\\
18 & 263504\\
19 & 523260\\
20 & 1050016\\
21 & 2102420\\
22 & 4198752\\
23 & 8378348\\
&
\end{array}
\end{array}
\]
\caption{Fixpoint numbers}\label{fig:4}
\end{table}

The zeta-function $Z(\alpha,T)$ of the operator $\alpha$ is
defined by
\[
Z(\alpha,T):=\exp\left(-\sum_{n=1}^\infty {\frac{N_n}{n}} T^n\right)\,.
\]
>From Table \ref{fig:4} we find
\begin{align*}
Z(\alpha,T)= &\; 1-4T^2-4T^3+4T^4+12T^5+4T^6-20T^7-22T^8+12T^9+32T^{10}
-12T^{11}+68T^{13}\\
& +32T^{14}-76T^{15}-179T^{16}-28T^{17}+172T^{18}+96T^{19}
+92T^{20}+32T^{21}\\
& -196T^{22}-112T^{23}+{\bf O}(T^{24})\,.
\end{align*}

In the next table we give for every prime $p$ with $3\le p\le 47$
a fixpoint of $\alpha^p$. This point has coordinates $a,\, b,\,
c,\, d,\, v,\, w,\, x,\, y$ in the field $\BF_{2^p}$. Here $t$ is
a primitive element of $\BF_{2^p}$ and MP is its minimal
polynomial over $\BF_2$.

These data close the gap between Corollary \ref{cor:big} and
Theorem \ref{conj:main} thus finishing the proof of Theorem
\ref{conj:P}. Although it is quite difficult to find fixpoints of
$\alpha^p$, it is easily checked, given the coordinates of point,
whether it is a fixpoint or not.

\begin{small}
\noindent $p=3$,  \qquad {\rm MP} $= t^3 + t + 1$\\
$a=1,\quad b=   t^4,\quad c=  t,\quad d=   t^6,\quad
v=   1,\quad w=   t^2,\quad x= t^4,\quad y= t^3$\,.
\medskip

\noindent $p=5$,\qquad  {\rm MP} $=t^5 + t^2 + 1$\\
$a=t^5,\quad b= 0  ,\quad c=  t^{14},\quad d=
1,\quad v=   t^9,\quad w=   0,\quad x= t^{19},\quad y= 1$.
\medskip

\noindent $p=7$,\qquad {\rm MP} $=t^7 + t + 1$\\
$a=t^5,\quad b= t^3  ,\quad c=  t^{56},\quad d=
t^{91},\quad v=   t^{80},\quad w=   t^{48},\quad x= t^{7},\quad y=
t^{59}$.
\medskip

\noindent $p=11,\qquad  {\rm MP}=t^{11} + t^2 + 1$\\
$a=t^3,\quad b= t^{228}  ,\quad c=  t^{151},\quad d=
t^{8},\quad v=   t^{192},\quad w=   t^{263},\quad x=
t^{1476},\quad y= t^{512}$.
\medskip

\noindent $p=13,\qquad  {\rm MP}=t^{13} + t^4+t^3+t + 1$\\
$a=t^9,\ b= t^{2129}  ,\ c=  t^{6077},\ d=
t^{7814},\ v=   t^{1152},\ w=   t^{2209},\ x= t^{7902},\ y=
t^{890}$.
\medskip

\noindent $p=17,\qquad  {\rm MP}=t^{17} +t^3+ 1$\\
$a=t^5,\quad b= t^{39028}  ,\quad c=  t^{30333},\quad
d=   t^{16060}, \quad v=   t^{2560},\quad w=   t^{59544},\quad x=
t^{64118},\quad y= t^{96318}$.
\medskip

\noindent $p=19,\qquad  {\rm MP}=t^{19} +t^5+t^2+t+ 1$\\
$a=t,\quad b= t^{45681}  ,\quad c=  t^{503015},\quad
d=   t^{8107},\quad v=   t^{1024},\quad w=   t^{115801},\quad x=
t^{237526},\quad y= t^{437263}$.
\medskip

\noindent $p=23,\qquad  {\rm MP}=t^{23} +t^5+ 1$\\
$a= t,$\\
$b= t^{22} + t^{21} + t^{19} + t^{16} + t^{15} + t^{12} + t^{10} +
t^7 + t^4 + t,$\\
$c= t^{21} + t^{20} + t^{17} + t^{16} + t^{13} + t^{12} + t^{11} +
t^{10} + t^6 + t^5 + t^3 + t,$\\
$d= t^{22} + t^{21} + t^{20} + t^{19} + t^{16} + t^{14} + t^{13} +
t^{12} + t^{10} + t^9 + t^6 + t^5 + t^2,$\\
$v=t^{22} + t^{21} + t^{20} + t^{19} + t^{17} + t^{15} + t^{13} +
t^{12} + t^{10} + t^7 + t^5 + t^4,$\\
$w=t^{22} + t^{18} + t^{15} + t^{13} + t^{12} + t^8 + t^7 + t^5 +
t,$\\
$x=  t^{21} + t^{19} + t^{18} + t^{15} + t^{13} + t^{11} + t^9 +
t^8 + t^6 + t^4 + t^2,$\\
$y=t^{19} + t^{18} + t^{15} + t^9 + t^7 + t^6 + t^3 + t^2 + t$.
\medskip

\noindent $p=29,\qquad  {\rm MP}=t^{29} +t^2+ 1$\\
$a=   t^2 + t,$\\
$b=t^{28} + t^{27} + t^{26} + t^{25} + t^{24} + t^{23} + t^{22} + t^{19} +
t^{16} + t^{13} + t^{12} + t^9 + t^8 + t^7 + t^6 + t^5 +
t^4 + 1,$\\
$c= t^{25} + t^{23} + t^{20} + t^{19} + t^{17} + t^{16} + t^{15} + t^{14}
+ t^{13} + t^9 + t + 1,$\\
$d=  t^{26} + t^{24} + t^{23} + t^{22} + t^{20} + t^{17} + t^{16}
+ t^{13} + t^{12} + t^6,$\\
$v= t^{28} + t^{26} + t^{20} + t^{18} + t^{17} + t^{11} + t^{10} +
t^8 + t^4 + t^3 + t,$\\
$w= t^{23} + t^{19} + t^{17} + t^{16} + t^{15} + t^{14} + t^{13} +
t^{12} + t^{11} + t^8 + t^6 + t^3 + t^2,$\\
$x= t^{26} + t^{25} + t^{19} + t^{16} + t^{13} + t^{11} + t^9 +
t^8 + t^7 + t^6 + t^4 + 1,$\\
$y= t^{27} + t^{24} + t^{23} + t^{22} + t^{21} + t^{15} + t^{13} +
t^9 + t^8 + t^6 + t^5 + t^4 + 1,$
\medskip

\noindent $p=31,\qquad  {\rm MP}=t^{31} +t^3+ 1$\\
$a=   t^3,$\\
$b= t^{30} + t^{27} + t^{25} + t^{23} + t^{21} + t^{20} + t^{18} + t^{17}
+ t^{16} + t^{15} + t^{14} + t^{12} + t^{11} + t^9 + t^8 +  t^5 + t^3 + t,$\\
$c= t^{27} + t^{24} + t^{20} + t^{19} + t^{18} + t^{16} + t^{14} + t^{12}
+ t^{10} + t^9 + t^8 + t^7 + t^5 + t^4,$\\
$d= t^{23} + t^{21} + t^{20} + t^{19} + t^{13} + t^8 + t^7 + t^6,$\\
$v= t^{24} + t^{10} + t^9 + t^2 + t,$\\
$w= t^{29} + t^{24} + t^{23} + t^{21} + t^{20} + t^{17} + t^{16} +
t^{15} + t^{11} + t^9 + t^4 + t^3 + t^2 + t,$\\
$x= t^{30} + t^{28} + t^{26} + t^{25} + t^{24} + t^{20} + t^{18} +
t^{17} + t^{16} + t^{15} + t^{11} + t^9 + t^8 + t^7 + t^5 + t,$\\
$y= t^{30} + t^{29} + t^{28} + t^{19} + t^{17} + t^{16} + t^{13} +
t^{11} + t^{10} + t^{8} + t^6 + t^5 + t^4 + t^3 + t$.
\medskip

\noindent $p=37,\qquad  {\rm MP}=t^{37} +t^5+t^4+t^3+t^2+t+ 1$\\
$a= t^2 + t,$\\
$b= t^{36} + t^{34} + t^{30} + t^{29} + t^{28} + t^{27} + t^{26} + t^{25} +
t^{23} + t^{21} + t^{17} + t^{14} + t^{11} + t^{10} + t^7
+ t^4 + t^3 + t^2 + t + 1,$\\
$c= t^{36} + t^{35} + t^{33} + t^{32} + t^{31} + t^{29} + t^{19} +
t^{14} + t^{11} + t^{10} + t^9 + t^7 + t^2,$\\
$d= t^{36} + t^{34} + t^{32} + t^{31} + t^{27} + t^{26} + t^{24} + t^{23} +
t^{22} + t^{20} + t^{19} + t^{18} + t^{16} + t^{12} +
t^{10} + t^9 + t^8 + t^7 + t^6 + t^4 + t^3 + t + 1,$\\
$v= t^{34} + t^{33} + t^{31} +t^{30} + t^{29} + t^{28} + t^{27} + t^{26} +
t^{19} + t^{18} + t^{17} + t^{16} + t^{14} + t^{13} + t^{10} + t^8
+ t^7 + t^6 + t^5 + t^4 + t^3 + t^2 + t + 1,$\\
$w= t^{36} + t^{35} + t^{34} + t^{33} + t^{32} + t^{29} + t^{26} + t^{25} +
t^{20} + t^{19} + t^{18} + t^{17} + t^{15} + t^{14} +
t^{12} + t^{10} + t^8 + t^7 + t^6 + t^5 + t^4 + t^2 + t + 1,$\\
$x=t^{36} + t^{31} + t^{29} + t^{28} + t^{27} + t^{24} + t^{21} + t^{19}
+ t^{18} + t^{16} + t^{15} + t^{14} + t^{12} + t^{10} + t^9 +  t^8 + t^3 + t^2
+ t,$\\
$y= t^{35} + t^{33} + t^{32} + t^{27} + t^{26} + t^{25} + t^{24} +
t^{23} + t^{19} + t^{18} + t^{15} + t^{11} + t^{10} + t^9 + t^7 +
t^4 + t^3 + t + 1$.
\medskip

\noindent $p=41$,\qquad  {\rm MP}$=t^{41} +t^3+ 1$\\
$a= t,$\\
$b= t^{39} + t^{37} + t^{36} + t^{35} + t^{34} + t^{33} + t^{30} + t^{28} +
t^{26} + t^{24} + t^{23} + t^{20} + t^{19} + t^{18} + t^{17} + t^{16} + t^{15}
+ t^{13} + t^{11} + t^{10} +$\\
\phantom{b= } $t^9 + t^6 + t^3 + t^2 + t +1,$\\
$c= t^{40} + t^{38} + t^{36} + t^{33} + t^{32} + t^{31} + t^{29} +
t^{28} + t^{27} + t^{26} + t^{24} + t^{23} + t^{22} + t^{19} +
t^{18} + t^{14} + t^{12} + t^{10} + t^7 + t^3 + 1,$\\
$d= t^{40} + t^{39} + t^{37} + t^{35} + t^{34} + t^{30} + t^{28} +
t^{27} + t^{23} + t^{21} + t^{19} + t^{18} + t^{15} + t^{14} +
t^{13}  + t^{11} + t^9 + t^7 + t^4 + t^2,$\\
$v= t^{38} + t^{37} + t^{36} + t^{35} + t^{34} + t^{32} + t^{31} +
t^{29} + t^{26} + t^{23} + t^{21} + t^{20} + t^{18} + t^{17} +
t^{15} + t^{14} + t^{13} + t^{12} + t^{11} + t^{10} +$\\
\phantom{v= } $t^7 + t^6 + t^4 + t^3 + t^2 + t,$\\
$w= t^{40} +t^{37} + t^{35} + t^{33} + t^{32} + t^{31} + t^{30} +
t^{29} + t^{28} + t^{27} + t^{26} + t^{25} + t^{24} + t^{20} +
t^{19} + t^{18} + t^{16} + t^{14} + t^{13} + t^9 +$\\
\phantom{v= } $t^8 + t^7 + t^6 + t^4 + t^2 + 1,$\\
$x= t^{38} + t^{36} + t^{34} + t^{33} + t^{32} + t^{29} + t^{21} +
t^{20} + t^{18} + t^{13} + t^{12} + t^7 + t^6 + t^2 + t + 1,$\\
$y=t^{40} + t^{39} + t^{35} + t^{26} + t^{23} + t^{22} + t^{19} +
t^{17} + t^{16} + t^{15} + t^{14} + t^{13} + t^7 + t^5 + t^4
+t^2$.
\medskip

\noindent $p=43$,\qquad  {\rm MP}$=t^{43} +t^6+t^4+t^3+ 1$\\
$a=   t^3,$\\
$b= t^{42} + t^{39} + t^{37} + t^{32} + t^{31} + t^{30} + t^{29} + t^{26} +
t^{23} + t^{19} + t^{18} + t^{17} + t^{10} + t^7 + t^6 + t^5 + t^4
+   1,$\\
$c=t^{40} + t^{39} + t^{38} + t^{36} + t^{35} + t^{34} + t^{33} + t^{32} +
t^{31} + t^{30} + t^{26} + t^{23} + t^{22} + t^{21} + t^{19} + t^{16} +
t^{14} + t^{13} + t^{12} + t^9 + t^4 + t^3 + t^2 + 1,$\\
$d= t^{39} + t^{38} + t^{37} + t^{36} + t^{35} + t^{34} + t^{29} + t^{28} +
t^{27} + t^{22} + t^{21} + t^{20} + t^{19} + t^{18} + t^{17} + t^{16} +
t^{13} + t^9 + t^7 + t^6 + t^4 + t^3 + t^2 +  + 1,$\\
$v= t^{40} + t^{37} + t^{35} + t^{33} + t^{31} + t^{30} + t^{26} +
t^{17} + t^{16} + t^{15} + t^{11} + t^8 + t^4 + t + 1,$\\
$w=t^{41} + t^{39} + t^{38} + t^{37} + t^{36} + t^{31} + t^{29} + t^{28} +
t^{24} + t^{21} + t^{18} + t^{17} + t^{14} + t^9 + t^8 + t^7 + t^5
+t^4 + t^3 + 1,$\\
$x=t^{41} + t^{39} + t^{34} + t^{33} + t^{32} + t^{29} + t^{25} + t^{24} +
t^{20} + t^{18} + t^{17} + t^{16} + t^{15} + t^{12} + t^{11} + t^{10} +
t^9 + t^7 + t^6 + t^5 + t^4 + t + 1,$\\
$y= t^{41} + t^{40} + t^{38} + t^{37} + t^{35} + t^{34} + t^{33} + t^{30} +
t^{29} + t^{28} + t^{27} + t^{26} + t^{24} + t^{19} + t^{18} + t^{16} + t^{15}
+ t^{12} + t^{11} + t^9 + t^6 + t^4 + t^3 + t^2 + t,$
\medskip

\noindent $p=47$,\qquad  {\rm MP}$=t^{47} +t^5+ 1$\\
$a=   t^2 + t + 1,$\\
$b= t^{46} + t^{44} + t^{43} + t^{41} + t^{40} + t^{39} + t^{37} + t^{36} +
t^{35} + t^{34} + t^{27} + t^{25} + t^{24} + t^{23} + t^{22} + t^{20} + t^{19}
+ t^{17} + t^{16} + t^{14} +$\\
\phantom{b= }$t^{13} + t^{12} + t^{11} + t^9 + t^8+ t^7 + t^6 + t^5 + t^2,$\\
$c=t^{40} + t^{33} + t^{31} + t^{30} + t^{29} + t^{26} + t^{25} + t^{24} +
t^{23} + t^{21} + t^{19} + t^{16} + t^{15} + t^{14} + t^{11} + t^{10} +
t^9 + t^8 + t^7 + t^6,$\\
$d=t^{44} + t^{42} + t^{41} + t^{39} + t^{36} + t^{35} + t^{31} + t^{27} +
t^{26} + t^{24} + t^{20} + t^{19} + t^{17} + t^{16} +
t^{15} + t^{13} + t^9 + t^7 + t^5 + t^4 + t^2 + t,$\\
$v=t^{44} + t^{41} + t^{40} + t^{38} + t^{37} + t^{34} + t^{33} + t^{32} +
t^{31} + t^{28} + t^{27} + t^{26} + t^{25} + t^{23} + t^{22} + t^{18} + t^{17}
+ t^{16} + t^{12} + t^{10} +$\\
\phantom{v= } $t^9 + t^7 + t^6 + t^5 + t^2 +t + 1,$\\
$w= t^{45} + t^{43} + t^{42} + t^{41} + t^{39} + t^{38} + t^{36} + t^{35} +
t^{34} + t^{32} + t^{30} + t^{24} + t^{23} + t^{21} + t^{20} + t^{16} +
t^{12} + t^9 + t^5 + t^4 + t^3 + t^2 + t,$\\
$x=t^{46} + t^{45} + t^{44} + t^{41} + t^{39} + t^{37} + t^{35} +
t^{33} + t^{31} + t^{30} + t^{29} + t^{27} + t^{26} + t^{25} +
t^{21} + t^{20} + t^{18} + t^{17} + t^{11} + t^9 +$\\
\phantom{x= } $t^8 + t^7 + t^6 + t^5 + t^3 + t^2,$\\
$y= t^{46} + t^{42} + t^{41} + t^{40} + t^{39} + t^{38} + t^{37} +
t^{35} +
 t^{34} + t^{33} + t^{32} + t^{31} + t^{30} + t^{28} + t^{26} + t^{25} +
t^{22} + t^{15} + t^{13} + t^{11} +$\\
\phantom{y= } $t^{10} + t^9 + t^6 + t^5 + t^4+t^2$.
\end{small}

\begin{table}[h]
\caption{Fixpoints in $\BF_{2^p}$}
\end{table}

\section{Appendix}

\subsection{Variations of proofs} \label{app:var}

\def\Z{{\bf Z}}
\def\I{{\bf I}}
\def\BZ{{\mathbb Z}}
\def\BP{{\mathbb P}}
\def\BF{\mathbb F}
\def\BG{{\mathbb G}}
\def\BN{\mathbb N}
\def\BR{{\mathbb R}}
\def\BC{{\mathbb C}}
\def\BQ{{\mathbb Q}}
\def\BT{\mathbb T}
\def\Z{{\bf Z}}
\def\I{{\bf I}}
\def\g{{\bf g}}
\def\S{{\mathcal S}}
\def\A{{\rm Aut}}
\def\a{{\lie Aut}}
\def\r{{\rm rk}}
\def\R{{\bf R}}
\def\Ad{{\rm Ad}}
\def\g{{\Gamma}}
\def\s{{{\sl R}_{\rm u}(G_1)}}
\def\t{{{\sl R}_{\rm u}(G_2)}}
\def\IQ{{\hbox{\Bbbb Q}}}
\def\IZ{{\hbox{\Bbbb Z}}}
\def\GL{{\bf GL}}

\def\lie{\mathfrak{a}}
\def \Box{\lower .1 em
          \vbox{\hrule \hbox{\vrule \hskip .6 em \vrule height .6 em} \hrule}}
\def \Mid{\quad \vrule \quad}

\def \rfish{\mathbin {\hbox {\msbm \char'157}}}

In this section we present another proof of the main theorem in
the $\PSL (2)$ case (Proposition \ref{prop1.1}).
We use the notations of Section \ref{sec:PSL} and consider the curve defined
by the ideal I.

The difference to the proof in Section \ref{sec:PSL} is the proof
of the absolute irreducibility of the polynomial $f_1 = J[1]$,
which uses here the analysis of the singularities.   Furthermore,
the Hasse--Weil Theorem is applied here to the normalisation of
the plane curve defined by $f_1$, while in Section \ref{sec:PSL}
it was applied to the curve defined by I and not to its projection
defined by $f_1$.

\begin{lemma}\label{lemmaA1}
  With the notations of Lemma $\ref{lemma3.6}$  we obtain, substituting
$ c= \tfrac{t^2-2t-1}{tb}$
\[
  \begin{array}{rcl}
bJ[2]\left(t,b,\dfrac{t^2-2t-1}{tb}\right) & = &J[1]\\[1.5ex]
J[3]\left(t,b,\dfrac{t^2-2t-1}{tb}\right) &= &0\\[1.5ex]
tbJ[4]\left(t,b,\dfrac{t^2-2t-1}{tb}\right) &= &-J[1]\\[1.5ex]
tb^2J[5]\left(t,b,\dfrac{t^2-2t-1}{tb}\right) &= &(1-tb)J[1]\,.
  \end{array}.
\]
\end{lemma}

\begin{proof}
  This is an easy computation.
\end{proof}

\begin{corollary}\label{coroA2}
  A point $(t,b)$ of the plane curve defined by $J[1] = 0$ with $tb \not= 0$
  defines a point $\left(t,b,\tfrac{t^2-2t-1}{tb}\right)$ of the curve defined
  by the ideal $I$.
\end{corollary}

\begin{proof}
  Just note that $J[1], \dots, J[5]$ is a Gr\"obner basis of $I$ and use Lemma \ref{lemmaA1}.
\end{proof}

\begin{remark}\label{remarkA2}
  In Section \ref{sec:PSL} we did not use this reduction to the case of a
  plane curve since this allowed a verification without computer.  We used the
  Hasse--Weil theorem involving the arithmetic genus which avoids an analysis
  of the singularities.  The arithmetic genus is 12 for the curve defined by
  $I$ and 15 for its projection to the plane defined by $f_1$.  The analysis
  of singularities allows us to use the geometric genus, which is 8.  In
  principle, this does not make a big difference because we are using a
  computer for small fields $\# q$, anyway.  For genus 15, resp.\ 12, resp.\
  8, Hasse--Weil guarantees rational points if $q \ge 977$, resp.\ $q \ge
  593$, resp.\ $q \ge 277$.  Hence, the analysis of the singularities reduces
  the number of small fields which have to be treated by computer.  On the
  other hand, we shall see that analysing the singularities, we have the
  disadvantage of treating the field $\F_{864007}$.  That such a
  large prime will play a special role in the analysis of singularities was
  unexpected for us.
\end{remark}

We reduced the problem to find a point $(t,b) \in \F_q$ on the
plane curve $\V (f_1)$ with $tb \not= 0$.   Note that $f_1(0,b) =
b^2+1$ and $f_1(t,0) = (t^2-2t-1)^2$. Hence there are at most four
points on the curve with $t=0$ or $b=0$. We shall show that there
are at least five points on such a curve. We did the calculations
in \textsc{Singular} and MAGMA to work with independent computer
algebra systems.

>From now on, we denote $P(t,b) = f_1(t,b)$.

We shall analyse the plane algebraic curve given by the polynomial
$P(t,b)$ over various (finite) fields. We put ${C}, \overline{C}$
for this curve (over the complex numbers $\BC$) and its projective
closure, respectively. We use the coordinate system $(t:b:z)$ in
the projective plane. The projective curve $\overline {C}$ is then
given by the homogeneous polynomial:
\begin{align*}
  \overline{P}(t,b,z)=& -b^3t^4-b^2t^5 + b^4t^2z +2b^3t^3z + 3b^2t^4z+
  bt^5z- 2b^2t^3z^2-4bt^4z^2+3bt^3z^3 + t^4z^3 \\& + 2b^2tz^4
  + 2bt^2z^4 - 4t^3z^4 + b^2z^5 +2t^2z^5 + 4tz^6 + z^7\,.
\end{align*}
If $p$ is a prime number, we put ${C}_p, \overline {C}_p$ for
these curves over $\overline \BF_p$. The curves ${C}_p$ and
$\overline {C}_p$ are then defined by the reductions of the
polynomials $P(t,b)$ and $\overline P(t,b,z)$ modulo $p$
respectively.

We use the standard formula:
\[
g(C)=\dfrac{\bigl({\rm degree}(\overline{C})-1)({\rm
degree}(\overline{C})-2\bigr)}{2} -\sum_{Q\in
\overline{C}(\bar{k})}\, \delta_Q
\]
where the local contributions of singular points $\delta _Q$ are
defined as $\dim_{\bar k}\widetilde{\mathcal O}_Q/\mathcal O_Q$;
here $\mathcal O_Q$ is the local ring of $Q$ and
$\widetilde{\mathcal O}_Q$ is the integral closure of $\mathcal
O_Q$ in the function field of $\overline{C}$.

We also define
$$A(\BF_q)=\#\{\, (t,b)\in \BF_q^2\Mid P(t,b)=0\, \}=\#{C}_p(\BF_q)$$
for a prime power $q$. The following tables contain the solution
numbers $A(L)$ for various finite fields $L$.
\begin{table}[ht]
\[
\begin{array}{lcrlcrlcrlcr}
A(\F_2)     & = &2  &\quad  A(\F_{2^8}) & = & 218 & \quad A(\F_3) & = & 0 &
\quad A(\F_{3^8}) & = & 6806 \\
A(\F_{2^2}) & = & 6 & \quad  A(\F_{2^9}) & = & 551 & \quad A(\F_{3^2}) & = &
14 & & &\\
A(\F_{2^3}) & = & 11 & \quad  A(\F_{2^{10}}) & = & 1026 & \quad A(\F_{3^3}) &
= &36 & & &\\
A(\F_{2^4}) & = & 10 & \quad A(\F_{2^{11}}) & = & 2048 & \quad A(\F_{3^4}) & = &78 & & &\\
 A(\F_{2^5}) & = & 32 & \quad A(\F_{2^{12}}) & = & 4279 & \quad A(\F_{3^5}) & = &190 & & &\\
A(\F_{2^6}) & = & 39 & \quad A(\F_{2^{13}}) & = & 7880 & \quad A(\F_{3^6}) & = &734 & & &\\
 A(\F_{2^7}) & = & 128 & \quad A(\F_{2^{14}}) & = & 16722 & \quad A(\F_{3^7})
 & = &2380 & & &
\end{array}
\]
\caption{Number of points on $C(\BF _q)$, $q=2^n$ or
$3^n$}\label{fig:5}
\end{table}

The numbers contained in  Table \ref{fig:5} and also those in
Table \ref{fig:6} can be obtained in microseconds on a computer.
We have, in fact, used MAGMA and verified this with
\textsc{Singular}.

\begin{table}[ht]
\[
\begin{array}{lcrlcrlcrlcrlcr}
A(\F_5) & = & 11 & \quad A(\F_{53}) & = & 53 & \quad A(\F_{109}) & = & 121 &
\quad A(\F_{181}) & = & 210 & \quad A(\F_{257}) & = & \\
A(\F_7) & = & 5 & \quad A(\F_{59}) & = & 36  & \quad A(\F_{113}) &
= & 122 & \quad A(\F_{191}) & = & 233 & \quad A(\F_{263}) & = & \\
A(\F_{11}) & = & 8 & \quad A(\F_{61}) & = & 72 & \quad A(\F_{127}) & = & 136
& \quad A(\F_{193}) & = & 223 & \quad A(\F_{269}) & = & \\
A(\F_{13}) & = & 16 & \quad A(\F_{67}) & = & 57 & \quad A(\F_{131}) & = &
121 & \quad A(\F_{197}) & = & 201 & \quad A(\F_{271}) & = & \\
A(\F_{17}) & = & 19 & \quad A(\F_{71}) & = & 76 & \quad A(\F_{137}) & = &
121 & \quad A(\F_{199}) & = & 167 & \quad A(\F_{277}) & = & \\
A(\F_{19}) & = & 15 & \quad A(\F_{73}) & = & 78 & \quad A(\F_{139}) & = &
134 & \quad A(\F_{211}) & = & 229 &  &  & \\
A(\F_{23}) & = & 9 & \quad A(\F_{79}) & = & 89 & \quad A(\F_{149}) & = &
140 & \quad A(\F_{223}) & = & 203 & & & \\
A(\F_{29}) & = & 45 & \quad A(\F_{83}) & = & 76 & \quad A(\F_{151}) & = &
164 & \quad A(\F_{227}) & = & 230 & & & \\
A(\F_{31}) & = & 33 & \quad A(\F_{89}) & = & 82 & \quad A(\F_{157}) & = &
161 & \quad A(\F_{229}) & = & 220 & & &\\
A(\F_{37}) & = & 36 & \quad A(\F_{97}) & = & 92 & \quad A(\F_{163}) & = &
170 & \quad A(\F_{233}) & = & 250 & & &\\
A(\F_{41}) & = & 61 & \quad A(\F_{101}) & = & 98 & \quad A(\F_{167}) & = &
136 & \quad A(\F_{239}) & = & 272 & & &\\
A(\F_{43}) & = & 32 & \quad A(\F_{103}) & = & 97 & \quad A(\F_{173}) & = &
167 & \quad A(\F_{241}) & = & 277 & & &\\
A(\F_{47}) & = & 42 & \quad A(\F_{107}) & = & 98 & \quad
A(\F_{179}) & = & 128 & \quad A(\F_{251}) & = & 233 & & &
\end{array}
\]
\caption{Number of points on $\kc(\F_p)$, $p$ prime} \label{fig:6}
\end{table}

We add the information:
$$A(\BF_{523})=474,\qquad A(\BF_{864007})\ge 3000 \leqno (*)$$
which can be also obtained by a simple computer calculation.

We shall show:

\begin{prop} \label{prop:F3}
If $q=p^k$ for a prime $p$ and $q\ne 2,\, 3$ then $A(\BF_q)\ge 5$.
\end{prop}


Our theorems would be much easier to prove if a rational point on
$C$ could be found. Unfortunately, even an extensive computer
search has not revealed such a point.

\medskip

We proceed by our analysis of the curve $C$. Consider affine
charts ${C}^1=\D (z)\cap\overline{C}$, ${C}^2=\D
(b)\cap\overline{C}$, and ${C}^3=\D (t)\cap\overline{C}$. The part
at infinity $\V (z)\cap\overline{C}$ is denoted by ${C}^\infty$.
Putting a prime $p$ as an index to $C$ stands then for the
analogous construction over $\overline \BF_p$. We have:

\begin{lemma} \label{lem:F2}
The part at infinity ${C}^\infty(\BC)$ consists exactly of the
points $(0:1:0)$, $(1:0:0)$, $(1:1:0)$. Also
${C}_p^\infty(\overline\BF_p)$ consists exactly of the points
$(0:1:0)$, $(1:0:0)$, $(1:1:0)$ for every prime $p$.
\end{lemma}

\begin{proof}
We find
$$\overline P(t,b,0)=-b^2t^4(b+t)$$
and the statement follows.
\end{proof}

We shall later prove Proposition \ref{prop:F3} by an application
of a Hasse--Weil estimate for the number of points on
$\overline{C}_p$. To do this, we have to understand the
singularities of $\overline{C}_p$ and also prove the absolute
irreducibility as the prime $p$ varies. The following contains a
description of the singularities of $\overline{C}$.

\begin{lemma} \label{lem:F3}
The projective curve $\overline{C}$ has the $4$ singular points
$$Q_1=(\omega +1:0:1),\ Q_2=(-\omega +1:0:1),\ Q_3=(1:0:0),\ Q_4=(0:1:0)$$
where $\omega=\sqrt {2}$. The points $Q_1$, $Q_2$, $Q_3$ are
ordinary double points whereas $Q_4$ is a singularity of type
$D_6$, that is $Q_4$ is a triple point with $3$ branches, two of
which are simply tangent. The projective curve $\overline{C}_\BQ$
is absolutely irreducible and $g(\overline{C}_\BQ) =8$.
\end{lemma}

\begin{proof} Most of this statement is computed by MAGMA and
\textsc{Singular}, the absolute irreducibility follows from Bezout. We
shall not carry this out here since we shall give the same
argument over the finite fields $\BF_p$ later.
\end{proof}

>From general theory it is clear that Lemma \ref{lem:F3} also holds
for the curve $\overline{C}_p$ for almost all primes $p$. To get
later explicit estimates, we have to find the exceptional set of
primes. We put:
$$S=\{\, 2,\, 23,\, 37,\, 523,\, 864007\,\}$$
and prove:

\begin{prop} \label{prop:F4}
Let $p$ be a prime with $p\notin S$. Then the projective curve
$\overline {C}_{p}$ has the $4$ singular points
$$Q_1=(\omega +1:0:1),\ Q_2=(-\omega +1:0:1),\ Q_3=(1:0:0),\ Q_4=(0:1:0)$$
where $\omega$ is a root of $x^2-2$ in $\overline\BF_p$. The
points $Q_1$, $Q_2$, $Q_3$ are ordinary double points whereas
$Q_4$ is a triple point with $3$ branches, two of which meet in
$Q_4$ of order $2$ and the third intersects them transversally
($D_6$-configuration). The projective curve $\overline{C}_{p}$ is
absolutely irreducible and $g(\overline {C}_p)=8$.
\end{prop}

\begin{proof}
We shall first find the singularities of $\overline{C}_{p}$. The
description of the singularities is obtained by looking at the
blow ups of $\overline{C}_{p}$ in the four singular points. These
can be computed by \textsc{Singular} or MAGMA.

We shall now analyze the singularities on the first affine patch
${C}^1_p$. Let ${\lie }_1$ be the ideal in $\BZ[t,b]$ generated by
$P$ and its derivatives with respect to $t,\, b$. A Gr\"obner
basis computation over $\BZ$ carried out in \textsc{Singular} or
MAGMA shows that $sb\in {\lie }_1$ where
$$s=35378249251012=4\cdot 23^2\cdot 37\cdot 523\cdot 864007.$$
We have $\left<{\lie }_1,b\right>=\left<b,t^2-2t-1\right>$. This
shows that the affine patch ${C}_p^1$ contains only the (distinct)
singular points $Q_1$, $Q_2$.

Let $\mathcal O$ be the ring of integers in $\BQ[\sqrt 2]$. Note
that ${\mathcal O}=\BZ[\sqrt 2]$. The points $Q_1$, $Q_2$ have
their coordinates in ${\mathcal O}/\mathfrak{p}$ where
$\mathfrak{p}$ is a prime ideal of $\mathcal O$ containing $p$.
The polynomial $P(v+\sqrt{2}+1,b)$ has
$H_2(b,v)=-(\sqrt{2}+1)b^2+2(\sqrt{2}+2)bv+8v^2$ as its
homogeneous part of lowest degree. A simple computation shows that
the only prime ideals $\mathfrak{p}$ of $\mathcal O$ with the
property that $H_2(b,v)$ is a square modulo $\mathfrak{p}$ are
$\mathfrak{p}_1=\sqrt{2} {\mathcal O}$ and
$\mathfrak{p}_2=(-3+4\sqrt{2}) {\mathcal O}$. Note that
$\mathfrak{p}_2$ contains $23$. The point $Q_2$ is analyzed
similarly.

This shows that for $p\notin S$ the affine patch ${C}_p^1$ only
contains the ordinary double points $Q_1$, $Q_2$ as singularities.
Note that $\delta_{Q_1}=\delta_{Q_2}=1$.

\medskip

We shall now analyze the singularities on the second affine patch
${C}^2_p$. Put $P_2(t,z)=\overline P(t,1,z)$. Let ${\lie }_2$ be
the ideal in $\BZ[t,z]$ generated by $P_2$ and its derivatives
with respect to $t,\, z$. A Gr\"obner basis computation over $\BZ$
carried out in MAGMA shows that $s_2b\in {\lie }_2$ where
$$s_2=66877597828=4\cdot 37\cdot 523\cdot 864007 .$$
We have $\left<{\lie }_2,z\right>=\left<z,t^2\right>$. This shows
that this affine patch contains only the singular point $Q_4$. The
polynomial $P_2(t,z)$ has $t^2z$ as its homogeneous part of lowest
degree, hence $Q_4$ is a triple point. Let $C$ be the affine curve
over $\overline\BF_p$ given by $P_2$ and $C_1$ be the curve given
by the polynomial $T_1(t,z)=P_2(t,zt)/t^3$. The polynomial $T_1$
has $t+z$ as its homogeneous part of lowest degree. This shows
that the blown up curve $C_1$ has only a simple point lying over
$(0,0)\in C$. Let $C_2$ be the curve given by the polynomial
$T_2(t,z)=P_2(tz,z)/z^3$. The polynomial $T_2$ has $t^2+z^2$ as
its degree $2$ homogeneous part. This shows that the blown up
curve $C_2$ has an ordinary double point lying over $(0,0)\in C$.

This shows that $3$ branches meet in $Q_4$. Two of them intersect
of order 2 and the third intersects these transversally. By M.~
Noether's formula ($\delta_{Q_4}$ equals the sum of $m_Q(m_Q-1)/2$
where $Q$ runs over all points in all blow-ups lying over $Q_4$
and $m_Q$ is the multiplicity of $Q$) we find $\delta_{Q_4}=4$.

\medskip

We shall now analyze the singularities on the third affine patch
${C}^3_p$. Put $P_3(b,z)=\overline P(1,b,z)$. Let ${\lie }_3$ be
the ideal in $\BZ[b,z]$ generated by $P_3$ and its derivatives
with respect to $b,\, z$. A Gr\"obner basis computation over $\BZ$
carried out in MAGMA shows that $sb\in {\lie }_3$. We have
$\left<{\lie }_3,b\right>=\left<b,z(z^2+2z-1)\right>$. This shows
that $Q_3$ is the only singular point on this patch which was not
found on the previous affine patches. The polynomial $P_3(b,z)$
has $-(b+z)z$ as its homogeneous part of lowest degree.

This shows that $Q_3$ is an ordinary double point and
$\delta_{Q_3}=1$.

\medskip

So far we have described the singularities of $\overline{C}_p$.
Also, the degree of ${C}_p$ being $7$, we find
$g(\overline{C}_p)=15-1-1-1-4=8$.

\medskip

It remains to prove the absolute irreducibility of
$\overline{C}_p$. Suppose $\overline{C}_{p}$ had $2$ components
$C_1$, $C_2$. From the description of the singularities we infer
the following possibilities for the intersection numbers:
$$I(C_1,C_2;Q_1)=0,\, 1, \quad I(C_1,C_2;Q_2)=0, 1,\quad I(C_1,C_2;Q_3)=0,\,
1,$$
$$I(C_1,C_2;Q_4)=0,\,2,\,  3.$$
Note that $Q_1,\ldots, Q_4$ do not lie on a common line. The
degree of $\overline{C}_{p}$ being $7$, Bezout's theorem shows
that $\overline{C}_{p}$ is absolutely irreducible.
\end{proof}

Although we shall not need all of it, we shall also describe the
situation for the exceptional primes $p$  in $S$.   We start with $p = 2$.

\begin{prop} \label{prop:F5}
The projective curve $\overline{C}_{2}$ has the $4$ singular
points
$$Q_1=(1:0:1),\ Q_2=(0:1:1),\ Q_3=(1:0:0),\ Q_4=(0:1:0).$$
The points $Q_1$, $Q_2$, $Q_3$ are double points whereas $Q_4$ is
a triple point. The point $Q_3$ is ordinary, at $Q_1$ two branches
with a common tangent touch of order $2$, $Q_2$ is an ordinary
cusp, at $Q_4$ two branches with distinct tangents meet, one of
them behaves like a third order cusp, the other is smooth in $Q_4$
(a $D_9$-configuration). The projective curve $\overline{C}_{2}$
is absolutely irreducible and $g(\overline{C}_{2})=6.$
\end{prop}

\begin{proof}
We shall first find the singularities of $\overline{C}_{2}$. The
description of the singularities is obtained by looking at the
blow ups of $\overline{C}_{2}$ in the four singular points. These
can be computed by MAGMA or \textsc{Singular}.

\medskip

We shall now analyze the singularities on the first affine patch
${C}^1_2$ ($z=1$). The Jacobian ideal of $P(t,b)$ is generated by
$b^2(b^2+1)$ and $t^2+b^2+1$. This shows that $Q_1,\, Q_2$ are the
only singularities on this affine patch.

Put $L_1(t,b)=P(t+1,b)$. The polynomial $L_1(t,b)$ has $b^2$ as
homogeneous component of lowest degree. Let $C$ be the affine
curve over $\overline\BF_2$ given by $L_1$ and $C_1$ be the curve
given by the polynomial $T_1(t,b)=L_1(tb,b)/b^2$. A look at $T_1$
shows that there is no point of $C_1$ lying over $(0,0)\in C$. Let
$C_2$ be the curve given by the polynomial
$T_2(t,b)=L_1(t,tb)/t^2$. The polynomial has $b(b+t)$ as its
homogeneous component of lowest degree. This shows that there is
an ordinary double point over $(0,0)\in C$ on $C_2$. Altogether we
find that two branches with a common tangent touch of order $2$ in
$Q_1$. This implies $\delta_{Q_1}=2$.

Put $L_2(t,b)=P(t,b+1)$. The polynomial $L_2(t,b)$ has $(t+b)^2$
as homogeneous component of lowest degree. Let $C$ be the affine
curve over $\overline\BF_2$ given by $L_2$. Both blow-ups of
$(0,0)\in C$ contain (the same) smooth point over $(0,0)\in C$.
This shows that $Q_2$ is a cusp (one branch passing through $Q_2$)
and $\delta_{Q_1}=1$.

\medskip

We shall now analyze the singularities on the second affine patch
${C}^2_2$ ($b=1$). The points $Q_2$ and $Q_4$ are the only
singularities on this affine patch. To analyze $Q_4$, put
$P_2(t,z)=\overline P(t,1,z)$. The polynomial $P_2$ has $t^2z$ as
its homogeneous component of lowest degree. Hence $Q_4$ is a
triple point with two distinct tangents. Let $C$ be the affine
curve over $\overline\BF_2$ given by $P_2$. In the first blow-up
($z=zt$) we find a simple point over $(0,0)\in C$, In the second
blow-up we find a point $Q_5$ of multiplicity $2$ with a double
tangent over $(0,0)\in C$. The blow-ups of $Q_5$ give one double
point with a double tangent $Q_6$. The blow-ups of $Q_6$ give one
simple point $Q_7$. This shows that at $Q_4$ two branches with
distinct tangents meet, one of them behaves like a third order
cusp, the other is smooth in $Q_4$. By M.~Noether's formula we
find $\delta_{Q_4}=3+1+1=5$.

\medskip

We shall now analyze the singularities on the third affine patch
${C}^3_2$ ($t=1$). The points $Q_1$ and $Q_3$ are the only
singularities on this affine patch. To analyze $Q_3$ put
$P_3(b,z)=\overline P(1,b,z)$. The polynomial $P_3$ has $bz$ as
its homogeneous component of lowest degree. This shows that $Q_3$
is an ordinary double point and $\delta_{Q_3}=1$.

\medskip

The analysis of the singularities being completed, we have found
$g(\overline {C}_2)=15-2-1-1-5=6$.

\medskip

Suppose $\overline{C}_{2}$ had $2$ components $C_1$, $C_2$. From
the description of the singularities we infer the following
possibilities for the intersection numbers:
$$I(C_1,C_2;Q_1)=0,\, 2, \ I(C_1,C_2;Q_2)=0,\ I(C_1,C_2;Q_3)=0,\,
1,$$ $$ I(C_1,C_2;Q_4)=0,\, 2.$$ These numbers cannot add up to
$6$ or more. The degree of $\overline{C}_{\BF_2}$ being $7$,
Bezout's theorem shows that $\overline{C}_{\BF_2}$ is absolutely
irreducible.
\end{proof}

\begin{prop} \label{Prop:F6}
The projective curve $\overline{C}_{{23}}$ has the $4$ singular
points
$$Q_1=(19:0:1),\ Q_2=(6:0:1),\ Q_3=(1:0:0),\ Q_4=(0:1:0).$$
The points $Q_2$, $Q_3$ are ordinary double points, $Q_1$ is an
ordinary cusp, whereas $Q_4$ is a triple point with $3$ branches,
two of which meet in $Q_4$ of order $2$ and the third intersects
them transversally ($D_6$-configuration). $Q_1$ is a cusp
singularity. The projective curve $\overline{C}_{{23}}$ is
absolutely irreducible and $g(\overline{C}_{{23}})=8$.
\end{prop}

\begin{proof} The singular points and their types were computed by
MAGMA. To complete the Bezout-argument notice that $Q_4$ does not
lie on a line with at the three double points.
\end{proof}

\begin{prop} \label{prop:F7}
The projective curve $\overline{C}_{{37}}$ has the $8$ singular
points
$$Q_1=(\omega+1:0:1),\ Q_2=(-\omega+1:0:1),\ Q_3=(1:0:0),\ Q_4=(0:1:0),$$
$$Q_5=(27:17:1),\ Q_6=(10:10:1),\ Q_7=(10:24:1),\ Q_8=(27:34:1)$$
where $\omega$ is a root of $x^2-2$ in $\overline\BF_{37}$. The
points $Q_1$, $Q_2$, $Q_3$, $Q_5$, $Q_6$, $Q_7$, $Q_8$ are
ordinary double points whereas $Q_4$ is a triple point with $3$
branches, two of which meet in $Q_4$ of order $2$ and the third
intersects them transversally ($D_6$-configuration). The
projective curve $\overline {C}_{{37}}$ is absolutely irreducible
and $g(\overline {C}_{{37}})=4$.
\end{prop}

\begin{proof} The singular points and their types were computed by
MAGMA. To complete the Bezout-argument notice that $Q_4$ does not
lie on a line with at least three of the double points, and also
that the points $Q_1$,\dots,$Q_8$ do not lie on a quadric.
\end{proof}

\begin{prop} \label{prop:F8}
The projective curve $\overline{C}_{{523}}$ has the $5$ singular
points
$$Q_1=(\omega+1:0:1),\ Q_2=(-\omega+1:0:1),\ Q_3=(1:0:0),$$
$$Q_4=(0:1:0), Q_5=(479:463:1)
$$ where $\omega$ is a root of $x^2-2$ in $\overline
\BF_{523}$. The points $Q_1$, $Q_2$, $Q_3$, $Q_5$ are ordinary
double points whereas $Q_4$ is a triple point with $3$ branches,
two of which meet in $Q_4$ of order $2$ and the third intersects
them transversally ($D_6$-configuration). The projective curve
$\overline {C}_{{523}}$ is absolutely irreducible and $g(\overline
{C}_{{523}})=7$.
\end{prop}

\begin{proof} The singular points and their types were computed by
MAGMA. To complete the Bezout-argument notice that $Q_4$ does not
lie on a line with at least three of the double points.
\end{proof}

\begin{prop} \label{prop:F9}
The projective curve $\overline{C}_{{864007}}$ has the $5$
singular points
$$Q_1=(767405:0:1),\ Q_2=(96604:0:1),\ Q_3=(1:0:0),$$
$$Q_4=(0:1:0),\ Q_5=(395579:564628:1).$$ The points $Q_1$, $Q_2$,
$Q_3$, $Q_5$ are ordinary double points whereas $Q_4$ is a triple
point with $3$ branches, two of which meet in $Q_4$ of order $2$
and the third intersects them transversally ($D_6$-configuration).
The projective curve $\overline {C}_{{864007}}$ is absolutely
irreducible and $g(\overline {C}_{864007} )=7$.
\end{prop}

\begin{proof} The singular points and their types were computed by
MAGMA. To complete the Bezout-argument notice that $Q_4$ does not
lie on a line with at least three of the double points.
\end{proof}

We are now ready for the

\noindent {\it Proof of Proposition \ref{prop:F3}:} We first
assume that the prime $p$ satisfies $p\notin S$ and also $p\ne 3$.
We shall then show that the statement of Proposition \ref{prop:F3}
is already true for $q=p$. We have to show that
\begin{equation}
\#{C}_{{p}}(\BF_p)\ge 5, \label{eq:F3-1}
\end{equation}
which is by Lemma \ref{lem:F2} equivalent to
\begin{equation}
\#\overline{C}_{{p}}(\BF_p)>7. \label{eq:F3-2}
\end{equation}
We write ${\mathcal D}_{{p}}$ for a nonsingular model of
$\overline{C}_{{p}}$ and
$$\pi_p: {\mathcal D}_{{p}}\to\overline{C}_{{p}}$$
for the birational projection. The map $\pi_p$ is defined over
$\BF_p$. Let $M\subset \overline{C}_{{p}}(\BF_p)$ be the set of
singular points. The map $\pi_p$ defines a bijection
$$\pi_p: {\mathcal D}_{{p}}(\BF_p)\smallsetminus \pi_p^{-1}(M)\to
\overline{C}_{{p}}(\BF_p)\smallsetminus M.$$ Since the
singularities of $\overline{C}_{{p}}$ are three double and a
triple point, we find:
$$\#{C}_{{p}}(\BF_p)\ge\#{\mathcal D}_{{p}}(\BF_p)-5.$$
Hence it is sufficient to show that
\begin{equation}
\#{\mathcal D}_{{p}}(\BF_p)>12.\label{eq:F3-3}
\end{equation}
By the Hasse--Weil estimate we know that
\begin{equation}
p-16\sqrt{p}\le \#{\mathcal D}_{{p}}(\BF_p).\label{eq:F3-4}
\end{equation}
If $p\ge 280$ the estimate (\ref{eq:F3-4}) implies
(\ref{eq:F3-3}). If $p\le 280$ the estimate (\ref{eq:F3-1}) is
already contained in our Table \ref{fig:6}.

For the primes $p\in S$, $p\ne 2,\, 3$ we also have
$\#{C}_{{p}}(\BF_p)\ge 5$ by Table \ref{fig:6} or the addition
(*).

For $q=2^k$, $k\ge 2$, we use Proposition \ref{prop:F5} and an
argument similar to the above to show $A(\BF_{2^k})\ge 5$ for
$k\ge 8$. The remaining values can be found in Table \ref{fig:5}.

For $q=3^k$, $k\ge 2$, we use Proposition \ref{prop:F4} and an
argument similar to the above to show $A(\BF_{3^k})\ge 5$ for
$k\ge 6$. The remaining values can be found in Table \ref{fig:5}.
\bigskip

\subsection{A variant of Zorn's theorem}
\label{app:Engel}

In this appendix we prove

\begin{prop} \label{prop:Engel}
Let $G$ be a finite group, and let $w=w(x,y)$ be a word in two
variables such that: 1) if $w(x,y)\equiv 1$ in $G$ then $G=\{1\}$;
2) the words $x$ and $w(x,y)$ generate the free group
$F_2=\left<x,y\right>$. Then $G$ is nilpotent if and only if it
satisfies one of the identities $[w(x,y),x,x,\dots ,x]=1$.
\end{prop}

\begin{proof} {\it Necessity}. Let $G$ be a nilpotent group of class $n$. Since
the element $e_n(x,y)=[w(x,y),x,\dots ,x]$ lies in the $n$th term
of the invariant series, $e_n(x,y)$ is an identity.

\noindent{\it Sufficiency}. We want to prove that any $G$
satisfying the identity $e_n(x,y)\equiv 1$ for some $n$ is
nilpotent. Assume the contrary.

Suppose that $n=1$. Then according to assumption (1) of the
proposition, the group $G$ is trivial. Let $n>1$. Let $\Gamma$
denote a minimal counterexample, i.e. a non-nilpotent group of the
smallest order satisfying the identity $e_n(x,y)\equiv 1$.
Obviously, all subgroups of $\Gamma$ are nilpotent. Then $\Gamma$
is a Schmidt group, i.e. a non-nilpotent group all of whose proper
subgroups are nilpotent (see \cite{Sch}, \cite{Re} for the
description of these groups). In particular, the commutator
subgroup $\Gamma'$ is the unique maximal Sylow subgroup in
$\Gamma$. Since $\Gamma'$ is nilpotent, it contains a non-trivial
center $Z(\Gamma')$. Take a nontrivial $a\in Z(\Gamma')$. For any
element $x\notin \Gamma'$ there exists $y\in G$ such that
$w(x,y)=a$ (condition (2)). Consider the sequence $[a,x,x,\ldots
,x]=[w(x,y),x,x,\ldots ,x]=e_n(x,y)$. There exists $n$ such that
$[w(x,y),x,x,\ldots, x]\equiv 1$. Let $n$ denote the smallest
number satisfying this equality, and let $b=[w(x,y),x,x,\ldots,
x]=e_{n-1}(x,y)$. Clearly, $b\in Z(\Gamma')$. Moreover,
$[b,x]=e_n(x,y)=1$ and hence $b$ is a nontrivial element from
$Z(\Gamma )$. Take $\bar\Gamma=\Gamma/Z(\Gamma )$. Then the order
of $\bar\Gamma$ is less than the order of $\Gamma$, hence
$\bar\Gamma$ is nilpotent. Therefore $\Gamma$ is nilpotent. Since
$e_n(x,y)$ is an identity in $\Gamma$, we get a contradiction.

The proposition is proved.
\end{proof}

\subsection{Profinite setting} \label{app:pro}

\subsubsection{Pseudovarieties of finite groups} \label{app:ps}

A {\it variety} of groups is a class $C$ of groups defined by some
set of identities $T$ (i.e. $G\in C$ if and only if for every
$u\in T$ the identity $u$ holds in $G$). Birkhof's theorem says
that $C$ is a variety if and only if $C$ is closed under taking
subgroups, homomorphic images, and direct products. To work with
classes of finite groups (which cannot be closed under taking
infinite direct products), one needs a more general notion.

\begin{defn}
A {\it pseudovariety} of groups is a class of groups closed under
taking subgroups, homomorphic images, and {\bf finite} direct
products.
\end{defn}

By Birkhoff's theorem every variety of groups is a pseudovariety.
We will be interested in pseudovarieties of all finite groups, all
finite solvable groups, and all finite nilpotent groups.

Let $F=F(X^0)$ be a free group with countable set of generators
$X^0$. Consider a sequence of words $u=u_1, u_2,\ldots , u_n,\dots
$ in $F$. The sequence $u$ determines a class of groups $V_u$ by
the rule: a group $G$ belongs to $V_u$ if and only if almost all
elements $u$ are identities in $G$. The class $V_u$ is a
pseudovariety. It turns out that this construction is universal:

\begin{theorem} \textrm{\cite{ES}}
For every pseudovariety of finite groups $V$ there exists a
sequence of elements $u\colon\mathbb N\to F$, $u=u_1, u_2,\ldots ,
u_n,\dots$ such that $V=V_u$.
\end{theorem}

We will consider a special class of sequences.

\begin{defn}
Let $X$ be a finite set. We say that a sequence of elements (not
necessarily distinct) $u=u_1, u_2,\ldots ,u_n,\dots $ of the free
group $F(X)$ is {\it correct} if given any group $G$, as soon as
an identity $u_n\equiv 1$ holds in $G$, for all $m>n$ the
identities $u_m\equiv 1$ hold in $G$, too.
\end{defn}

As above, a correct sequence $u$ defines a pseudovariety of groups
$V$ by the rule: $G\in V$ if and only if some identity $u_n\equiv
1$, $u_n\in u$, holds in $G$.

\begin{remark}\label{subseq}
If $u$ is a correct sequence defining a pseudovariety $V$ and $v$
is a subsequence of $u$, then $v$ is also correct and defines the
same pseudovariety $V$.
\end{remark}

Let $F=F(x,y)$ and
\begin{equation}
\begin{array}{ccl}
{}e_1 & = & [\ x,y\ ], \\
{}e_{n+1} & = & [e_n ,y\ ], \dots
\end{array}
\label{seq:gen0}
\end{equation}

This sequence is correct and defines the pseudovariety of all
finite Engel groups. According to Zorn's theorem \cite{Zo}, this
pseudovariety coincides with the pseudovariety of all finite
nilpotent groups.

Our main sequence of quasi-Engel words
\begin{equation}
\begin{array}{ccl}
u_1 & = & w=x^{-2}y^{-1}x, \\
u_{n+1} & = & [x\,u_n\,x\min,y\,u_n\,y\min], \dots
\end{array}
\label{seq:gen1}
\end{equation}
is also correct, and according to Theorem \ref{conj:P} it defines
the pseudovariety of all finite solvable groups.

\subsubsection{Residually finite groups} \label{app:res}

\begin{defn} We say that a group $G$ is {\it residually finite} if
the intersection of all its normal subgroups of finite index
$H_\alpha$, $\alpha \in I$, is trivial.
\end{defn}

Define a partial order on the set $I$ by: $\alpha <\beta$ if and
only if $H_\beta \subset H_\alpha$. The intersection of two normal
subgroups of finite index is also of finite index, and therefore
for every $\alpha, \beta \in I$ there is $\gamma \in I$ such that
$\alpha <\gamma$, $\beta <\gamma$. Thus the set $I$ is directed.

Denote $G_\alpha=G/H_\alpha$. If $ \alpha<\beta$ then there is a
natural homomorphism $\varphi_\alpha^\beta: G_\beta\to G_\alpha$.
If $gH_\beta$ is an element of $G_\beta$ then its image in
$G_\alpha$ is $gH_\alpha$. Let $\bar G$ be the direct product of
all $G_\alpha$. Then there is an embedding $G\to \bar G$ which
associates to each $g\in G$ the element $\bar g=
(gH_{\alpha})_{\alpha\in I}$. Hence $G$ can be approximated by
finite groups $G_\alpha$, i.e. if $f$, $g$ are distinct elements
of $G$ then there is $\alpha$ such that $\bar f_\alpha$ and $\bar
g_\alpha$ are distinct elements of $G_\alpha$.

A group $G$ is regarded as a topological group, with the topology
defined by the system of neighbourhoods of 1 consisting of all
normal subgroups of finite index $H_\alpha$. The system of
neighbourhoods of an element $g\in G$ is given by the cosets
$gH_\alpha$. The group $\bar G$ is also a topological group. To
define the topology, consider the projections
$\pi_\alpha\colon\bar G\to G_\alpha$. Let $\ker
\pi_\alpha=U_\alpha$. Then $\bar G/U_\alpha$ is isomorphic to
$G_\alpha=G/H_\alpha$. For every $g\in G$ the element $\bar g$
lies in $U_\alpha$ if and only if $g\in H_\alpha$. The system of
neighbourhoods of 1 in $\bar G$ consists of all finite
intersections of normal subgroups $U_\alpha$. This defines the
Tikhonov topology on $\bar G$. Since all groups $G_\alpha$ are
finite, the group $\bar G$ is compact.

Let $g_1, \ldots, g_n,\ldots$ be a sequence of elements of $G$. As
usual, we say that this sequence tends to 1 if for every
neighbourhood $H_\alpha$ there exists a natural number
$N=N(\alpha)$ such that for all $ n>N$ the element $g_n$ lies in
$H_\alpha$.

\begin{defn} \label{def:id}
Let $F=F(X)$ be a free group. We say that a sequence $u=u_1,
\ldots, u_n, \ldots$ of elements of $F$ {\it identically converges
to 1} in a group $G$ if for any homomorphism $\mu\colon F\to G$
the sequence $\mu(u)=\mu(u_1), \ldots,\mu(u_n), \ldots$ tends to 1
in $G$. In this case we write $u\equiv 1$ in $G$.
\end{defn}

\begin{prop}\label{ident}
Let $X$ be a finite set. If a sequence  $u_1, \ldots, u_n, \ldots$
identically converges to 1 in $G$ then for every neighbourhood
$H_\alpha$ there exists $N=N(\alpha)$ such that all $u_n, n>N$,
are identities of the group $G/H_\alpha$.
\end{prop}

\begin{proof}
Take a homomorphism $\mu\colon F\to G$, and let $\mu^0\colon G\to
G/H_\alpha$ be the natural projection. Then $\nu=\mu^0\mu$ is a
homomorphism $F\to G/H_\alpha$, and every homomorphism $\nu\colon
F\to G/H_\alpha$ can be represented in this way. Since both
$G/H_\alpha$ and $X$ are finite, the set of different $\nu$'s is
also finite. Denote them $\{\nu_1, \ldots, \nu_k\}$.

Define an equivalence relation on the set of all homomorphisms
$\mu\colon F\to G$ by: $\mu_1\equiv\mu_2$ if
$\mu^0\mu_1=\mu^0\mu_2$. For an arbitrary $u\in F$ we have
$\mu(u)\in H_\alpha$ if and only if $\mu^0\mu(u)=1$. Thus, if
$\mu_1\equiv\mu_2$ then for every $u\in F$ we have $\mu_1(u)\in
H_\alpha$ if and only if $\mu_2(u)\in H_\alpha. $ Indeed, let
$\mu_1(u)\in H_\alpha$. Then $\mu^0\mu_1(u)=1=\mu^0\mu_2(u)=1$,
and $\mu_2(u)\in H_\alpha$.

For every $\nu_i$, $i=1, \ldots, k$ take $\mu_i$ such that
$\mu^0\mu_i=\nu_i$. Consider the equivalence classes
$[\mu_1],\ldots,[\mu_k]$. Each $\mu\colon F\to G$ belongs to one
of these classes. Since the sequence $u_1, \ldots, u_n, \ldots$
identically converges to 1 in $G$, for every $\mu\colon F\to G$
there exists $N=N(\alpha,\mu)$ such that $\mu(u_n)\in H_\alpha$
for $n>N$. Let $N_0$ be the maximum of $N(\alpha, \mu_i)$,
$i=1,\ldots, k$. If $n>N_0$ then $\mu_i(u_n)\in H_\alpha$ for
every $\mu_i$. Since every $\mu$ is equivalent to some $\mu_i$, we
have $\mu(u_n)\in H_\alpha$  for every $\mu$. This means that
$\nu(u_n)=1$ for every $\nu: F\to G/H_\alpha$. Thus the element
$u_n$ defines an identity of the group $G/H_\alpha$.
\end{proof}

\subsubsection{Profinite groups} \label{app:pro-gr}

We now focus on profinite groups, with a goal to establish a
relationship with pseudovarieties and give another reformulation
of our main result. Generalities on profinite groups can be found
in \cite{RZ}, \cite{Al1}, etc. We recall here some basic notions.

Let $V$ be a pseudovariety of finite groups. Given a group $G$, consider all
its normal subgroups of finite index $H_\alpha$ such that
$G/H_\alpha=G_\alpha\in V$. If the intersection of all these $H_\alpha$ is
trivial, we say that $G$ is a {\it residually $V$-group}. This is a
topological group with $V$-topology (the subgroups $H_\alpha$ as above are
taken as the neighbourhoods of 1).

Let $\bar G$ be the direct product of all $G_\alpha$. Denote by
$\widehat G$ a subgroup in $\bar G$ defined as follows: an element
$f\in \bar G$ belongs to $\widehat G$ if and only if for every
$\alpha$ and $\beta$ such that $H_\beta\subset H_\alpha$ the
equality $\varphi_\alpha^\beta(f_\beta )=f_\alpha$ holds. Denote
$f_\alpha= g_\alpha H_\alpha$. Then
$$
\varphi_\alpha^\beta(g_\beta H_\beta)=g_\alpha H_\alpha= g_\beta
H_\alpha.
$$
Recall that $\varphi_\alpha^\beta$ are natural homomorphisms.

The group $\widehat G$ turns out to be the completion of $G$ in
its $V$-topology \cite{ESt}.

Such a group $\widehat G$ is called a {\it pro-$V$-group}. If $V$
is the pseudovariety of all finite groups, $\widehat G$ is called
{\it a profinite group}. Thus in the class of all profinite groups
one can distinguish subclasses related to particular
pseudovarieties $V$.

A free group $F=F(X)$ is residually finite. Take all normal
subgroups of finite index in $F$. They define the profinite
topology in $F$. Denote by $\widehat F$ the completion of $F$ in
this topology. This group is a free profinite group (see, for
example, \cite{RZ}).

Indeed, if $\widehat G$ is the profinite completion of an
arbitrary residually finite group $G$, then every map $\mu\colon
X\to \widehat G$ induces a homomorphism $\mu\colon F\to \widehat
G$ which turns out to be a continuous homomorphism of topological
groups and therefore induces a continuous homomorphism
$\hat{\mu}\colon\widehat F\to \widehat G$.

Another approach to free profinite groups is based on the idea of
implicit operations (cf. \cite{Al1}, \cite{Al2}, \cite{AV},
\cite{MSW}, \cite{We}, etc.). This approach has a lot of
advantages but we do not use it since it needs additional notions
which are not necessary for our aims.

\begin{defn}
Let $f\in \widehat F$. The expression $f\equiv 1$ is called a {\it
profinite identity} of a profinite group $\widehat G$ if for every
continuous homomorphism $\widehat \mu\colon \widehat F\to \widehat
G$ we have $\widehat\mu(f)=1$.
\end{defn}

\begin{defn}(see also \cite{AV}, \cite{Al1})
A variety of profinite groups (for brevity, a {\it provariety}) is
a class of profinite groups defined by some set of profinite
identities.
\end{defn}

An analogue of Birkhoff's theorem for profinite groups says that a
class of profinite groups is a provariety if and only if it is
closed under taking closed subgroups, images under continuous
homomorphisms, and direct products. This implies that for an
arbitrary pseudovariety $V$ of finite groups, the class of all
pro-$V$-groups is a provariety. The converse statement is also
true. For any provariety $C$ there exists a pseudovariety of
finite groups $V$ such that the class of all pro-$V$-groups
coincides with $C$. In the case where $V$ is a correct
pseudovariety of finite groups (i.e., is defined by a correct
sequence), one can construct identities defining the provariety of
pro-$V$-groups in an explicit form.

Let $X$ be a finite set. Let $u=u_1, \ldots, u_n, \ldots$ be a
sequence of elements of a free group $F=F(X)$. Since $\widehat F$
is a compact group, there exists a convergent subsequence $v=v_1,
\ldots, v_m, \ldots$ of $u$.

\begin{proposition} \label{pro}
Let  $v=v_1, v_2, \ldots, v_n , \ldots$ be a convergent sequence
of elements of $F$ with $\lim \bar v_n=f$. Let $\widehat G$ be a
profinite group. Then the identity $f\equiv 1$ holds in $\widehat
G$ if and only if $v\equiv 1$ in $G$ $($i.e. $v$ identically
converges to 1 in $G$, see Definition $\ref{def:id})$.
\end{proposition}

\begin{proof}
First of all the sequence $\bar g_1, \ldots, \bar g_n, \ldots $
converges to 1 in $\widehat G$ if and only if $ g_1, \ldots, g_n,
\ldots $ converges to 1 in $G$.

Let the identity $f\equiv 1$ be fulfilled in $\widehat G$. Then
$$\widehat \mu(f)=\lim\widehat\mu(\bar v_n) =\lim
\overline{\mu(v_n)}=1.$$ Thus, $\lim \mu(v_n)=1 $ in $G$. This
means that $v\equiv 1$ in $G$. Conversely, let $v\equiv 1$ in $G$.
Then for every $\mu\colon F\to G$ the sequence $\mu(v)$ converges
to 1 in $G$. The sequence $\overline{\mu(v)}$ converges to 1 in
$\widehat G$. Using
$$\lim\widehat{\mu}(\bar v_n) =\lim
{\overline{\mu( v_n)}}=1=\widehat{\mu} (f),$$ we conclude that
$\widehat{\mu}(f)=1$ for arbitrary $\mu$. This means that the
identity $f\equiv 1 $ holds in $\widehat G$.
\end{proof}

Let $V$ be a pseudovariety of finite groups defined by a correct sequence
\mbox{$u=u_1$,} \mbox{$u_2, \ldots, u_n , \ldots$,} and let \mbox{$v=v_1, v_2,
  \ldots, v_n , \ldots$} be a convergent subsequence of $u$. Denote the limit
of $v$ by $f$. Since $u$ is a correct sequence, $v$ determines the same class
$V$ as $u$.

\begin{theorem}\label{var}
With the above notation, the class of all pro-$V$-groups is the
provariety defined by the profinite identity $f\equiv 1$.
\end{theorem}
\begin{proof}
Let the profinite identity $f\equiv 1$ hold in a profinite group
$\widehat G$. Then by Proposition \ref{pro}, $v\equiv 1$ in $G$.
Proposition \ref{ident} implies that for every neighbourhood $H$
in $G$ and all sufficiently large $n$ the identity $v_n\equiv 1$
holds in $G/H$. This means that $G/H$ lies in $V$ and $\widehat G$
is a pro-$V$-group.

Conversely, let $G/H$ lie in $V$. By the definition of $V$, this
means that $v$ identically converges to 1 in $G$. Therefore, the
identity $f \equiv 1$ holds in $\widehat G$.
\end{proof}

\begin{remark} Although all convergent subsequences of a correct
sequence define the same pseudovariety, their limits may be
different. For example, consider a correct sequence of the form
$u=v_1,av_1a^{-1}, v_2, av_2a^{-1}, \ldots,v_n, av_n a^{-1} ,
\ldots,$ where $a\in F$ and $v=v_1, v_2, \ldots, v_n , \ldots$ is
a correct convergent sequence. If the limit of the subsequence $v$
is $f$, we get a new convergent subsequence $v'=av_1a^{-1},
av_2a^{-1}, \ldots, av_n a^{-1} , \ldots$ with limit $afa^{-1}$.
However, the elements $f$ and $afa^{-1}$ define the same variety.
\end{remark}

\begin{corollary}\label{e} Let $F=F(x,y)$, and let $u_n$ be defined by
\begin{equation}
\begin{array}{ccl}
u_1 & = & w, \\
u_{n+1} & = & [u_n,y], \dots
\end{array}
\label{seq:nil}
\end{equation}
where $w=[x,y]$ or $w$ is any word satisfying the conditions the
hypotheses of Proposition $\ref{prop:Engel}$.



Let $v_1, v_2, \ldots, v_m , \ldots$ be any convergent subsequence
of $(\ref{seq:nil})$ with limit $f$ from $\widehat F$. Then the
identity $f\equiv 1$ defines the profinite variety of pronilpotent
groups.
\end{corollary}

\begin{proof}
The corollary immediately follows from Proposition
\ref{prop:Engel}, Zorn's theorem, and Theorem \ref{var}.
\end{proof}

\begin{theorem} \label{qe}
Let $F=F(x,y)$, let
\begin{equation}
\begin{array}{ccl}
u_1 & = & w=x^{-2}yx^{-1}, \\
u_{n+1} & = & [xu_nx\min ,yu_ny\min ], \dots
\end{array}
\label{seq:solv}
\end{equation}
be our main sequence, and let $v_1, v_2, \ldots, v_m , \ldots$ be
any convergent subsequence of $(\ref{seq:solv})$ with limit $f$
from $\widehat F$. Then the identity $f\equiv 1$ defines the
profinite variety of prosolvable groups.
\end{theorem}

\begin{proof}
The theorem immediately follows from Theorems \ref{conj:P} and
\ref{var}.
\end{proof}


We can now state the profinite analogue of the Thompson--Flavell
theorem.

\begin{corol} \label{cor:two-gen}
A profinite group $G$ is prosolvable if and only if every closed
two-generator subgroup of $G$ is prosolvable.
\end{corol}

\begin{proof}
Let every two-generator subgroup of $\widehat G$ be prosolvable.
Take an element $f\in \widehat F(x,y)$ which is the limit of a
convergent subsequence of our sequence $u$. Let $\mu$ be an
arbitrary continuous homomorphism  $\widehat F(x,y)\to \widehat
G$. Then $\mu(f)=1$ since $\mu(f)$ lies in a two-generator
subgroup of $\widehat G$. This is true for arbitrary $\mu$ and,
therefore, $f\equiv 1$. According to Theorem \ref{qe}, $\widehat
G$ is prosolvable.
\end{proof}

Corollary \ref{e} and Theorem \ref{qe} should be compared with
results of J.~Almeida \cite{Al2}. He used the language of implicit
operations and the notion of $n!$-type convergent subsequence to
get nice proofs of theorems of similar type. He also noticed that
if our main theorem about solvable groups is true for the sequence
${{}^w\!u_n}$ with $w=[x,y]$, the $n!$ version of the
corresponding statement is also true. We were now able to
formulate the theorem for our sequence (\ref{seq:gen1}).

\begin{remark}
 It is still not clear whether one can take $w=[x,y]$. We believe
that the answer is ``yes".




\end{remark}



\end{document}